\documentclass[a4paper,12pt]{article}
\usepackage{amsmath}
\usepackage{amssymb}
\usepackage{tabularx}
\usepackage{enumerate}
\usepackage{graphicx}
\usepackage{subfigure}
\usepackage{color}
\usepackage[pagewise]{lineno}
\usepackage{longtable}
\usepackage{epstopdf}
\usepackage{float}
\usepackage{setspace}

\usepackage[justification=centering]{caption}
\usepackage{multirow}
\usepackage{cases}
\usepackage[compress]{cite}
\usepackage{lineno}
\usepackage{color}
\usepackage{setspace}
\usepackage{multirow}
 \usepackage{rotating}
\usepackage{geometry}
\newtheorem{thm}{Theorem}[section]
\newtheorem{lem}[thm]{Lemma}
\newtheorem{cor}{Corollary}[section]
\newtheorem{exm}{Example}[section]

\newtheorem{rmk}{Remark}[section]

\newtheorem{proof}{Proof}

\newcommand{\qed}{\hspace{1em}\mbox{\raisebox{0.65ex}{\fbox{}}}}

\numberwithin{equation}{section}

\newcommand{\be}{\begin{equation}}
\newcommand{\ee}{\end{equation}}
\newcommand\bes{\begin{eqnarray}} \newcommand\ees{\end{eqnarray}}
\newcommand{\bess}{\begin{eqnarray*}}
\newcommand{\eess}{\end{eqnarray*}}

\newcommand{\bpf}{{\bf Proof.\ \ }}
\newcommand{\epf}{\mbox{}\hfill $\Box$}

\begin{document}

\thispagestyle{empty}
\title{On an age-structured juvenile-adult model with harvesting pulse in moving and heterogeneous environment\thanks{ Xu 
is supported by the Natural Science Foundation of Jiangsu Province, PR China (No. BK20220553);
Lin is supported by the National Natural Science Foundation of China (No. 12271470); Zhu is supported by Natural Sciences and Engineering Research Council of Canada.}}
\date{\empty}
\author{\thanks{Corresponding author. Email: zglin@yzu.edu.cn (Z. Lin).}  \\
{\small $^1$ School of Mathematical Science, Yangzhou University, Yangzhou 225002, China}\\
}

\author{Haiyan Xu$^1$, Zhigui Lin$^1$ and Huaiping Zhu$^2$\thanks{Corresponding author. }  \\
{\small $^1$ School of Mathematical Science, Yangzhou University, Yangzhou 225002, China}\\
{\small $^2$ Laboratory of Mathematical Parallel Systems (LAMPS),}\\
{\small Department of Mathematics and Statistics,York University, Toronto, ON, M3J 1P3, Canada}
}

\maketitle
\begin{quote}
\noindent
{\bf Abstract} { 
 \small This paper concerns an age-structured juvenile-adult model with harvesting pulse and moving boundaries in a heterogeneous environment, in which
 the moving boundaries describe the natural expanding front of species and
human periodic pulse intervention  is carried on the adults.
 The principal eigenvalue is firstly defined and its properties involving the intensity of harvesting and length of habitat sizes are analysed. Then the criteria to determine whether the species spread or vanish is discussed, and some relevant sufficient conditions characterized by pulse are established. Our results reveal that the co-extinction or coexistence of species is influenced by internal expanding capacities from species itself and external harvesting pulse from human intervention, in which the intensity and timing of harvesting play key roles. The final numerical approximations indicate that the larger the harvesting rate and the shorter the harvesting period, the worse the survival of all individuals due to the cooperation among juveniles and adults, and such harvesting pulse can even alter the situation of species, from persistence to extinction. In addition, expanding capacities also affect or alter the outcomes of spreading-vanishing}.

\end{quote}

\section{Introduction}

As some important natural phenomena, cooperation, competition and predation can promote the virtuous cycle between organisms and organisms or between organisms and the environment\cite{CC,HB}. Take symbiotic and cooperative relationship as an example, leguminous plants provide organic matter for the rhizobium, which in turn convert nitrogen into nitrogen-containing material that plants can absorb \cite{PA}. Whereafter, reaction-diffusion equations have been extensively proposed as credible models to describe such phenomena in biology, and many theoretical developments can be observed, see \cite{BL,LZ,WKS,WLL} for mutualistic model, \cite{TZ,ZZ} for competition model, and \cite{DH,WW} for predator-prey model.

The above-mentioned reaction-diffusion systems and many following works 
have modeled a spatially structured but physiologically unstructured population. Actually, population often have highly structured life cycle by age or other attributes, which lead to  evident differences in dispersal ability, reproduction rate and mortality. Moreover, the dispersal of species is usually influenced by multitudinous environmental factors, such as weather conditions, temperature and food availability etc., which may vary with time $t$. Therefore, with the inspirations by variable coefficients and physiologically structured models, see Alqawasmeh and Lutscher \cite{AL} for model with spatially heterogeneous but temporally homogeneous environment, Cantrell, Cosner and Mart\'{i}nez in \cite{CCM} with spatially variable parameters except diffusion rates, and Brown, Fang, Huang, etal. \cite{HZ,BZ,F} for age-stage reaction-diffusion models, we here consider an age-structured juvenile-adult model
\begin{eqnarray*}
\left\{
\begin{array}{lll}
J_t-d_1J_{xx}=b(t)A-a(t)J-m_1(t)J-\alpha_1(t)J^2,\\[1mm]
A_t-d_2A_{xx}=a(t)J-m_2(t)A-\alpha_2(t)A^2\\
\end{array} \right.
\end{eqnarray*}
in a spatially homogeneous but temporally heterogeneous landscape. All definitions of model parameters are expressed in \eqref{a01}.

Previous works on reaction-diffusion models have mainly focused on the fixed domain, however, the biological habitat is constantly shifting by the movement of species. Since the work of Du and Lin in 2010 \cite{DL}, there have been extensive theoretical developments involving moving boundaries. For example, Li and Lin \cite{LL} proposed a mutualistic ecological model with advection and two free boundaries, in which the asymptotic behaviors of global solution and criteria governing spreading or vanishing were investigated. This model was later extended to a mutualistic model with two different moving boundaries, which was set up by Zhang and Wang in \cite{ZW} and the spreading speeds of two free boundaries were analyzed. See also \cite{W,WZ} for two-species competition-diffusion model, \cite{WZD,TN,YMB} for prey-predator model.
Recently, nonlocal diffusion as a kind of long-distant diffusion has been widely considered in population and epidemic model  with free boundaries. For instance, mutualistic models with different nonlocal diffusion and free boundary were investigated, see \cite{DN,LW} and references therein. Additionally, Wang \cite{ZLW} studied a general and cooperative free boundaries problem with local-nonlocal diffusion, and Pu, Lin and Lou \cite{PLL} proposed a WNv nonlocal model with free boundaries and seasonal session.  Motivated by the previous works aforementioned, to describe the spatial spreading, which means the initial living domain $\Omega_0(=(-h_0,h_0))$ with $t=0$ will evolve into the habitat $\Omega(t)(=(g(t),h(t)))$ as time $t$ increases, we take the Stefan type free boundary
conditions (see \cite{B} for detailed deduction)
\begin{eqnarray*}
\left\{
\begin{array}{lll}
h'(t)=-\mu_1J_x(t,h(t))-\mu_2A_x(t,h(t)), \\[1mm]
g'(t)=-\mu_1J_x(t,g(t))-\mu_2A_x(t,g(t))\\[1mm]
\end{array} \right.
\end{eqnarray*}
into consideration. It implies that expanding fronts $g(t)$ and $h(t)$ expand as a rate that is proportional to species gradient of juveniles and adults. The constants $\mu_1$ and $\mu_2$ measure the spreading capacities of juveniles and adults in new habitat, respectively. Of many other works on different moving boundaries, let us list, \cite{BD,CL} with an extra force ($\alpha(t)$) against the boundary expansion such that
$$h'(t)=-\mu u_x(t,h(t))-\alpha(t),$$
and \cite{FL} with habitat prefers a form of nonlocal boundaries
$$
h'(t)=\mu\int_0^{h(t)}u(t,x)w(h(t)-x)dx
$$
on the assumption that the shifting rate of boundary $h(t)$ is determined by the weight function $w$, and it extended previous dichotomy to the expanding-balancing-shrinking trichotomy.

Apart from moving environment, transient disturbances from human intervention (harvesting, hunting, pesticide and vaccination, etc.) also can have enormous impact on the persistence or extinction of species. For instance, spraying pesticides or insecticides on crops in agriculture \cite{TC} and vaccinating during the COVID-19 pandemic \cite{V,TM} can be regarded as the human efforts for sustainable developments of ecology, economy and society. In recent years, pulses, as a kind of transient and effective control methods, have been frequently introduced into mathematical models to characterize the discrete-continuous phenomenon in nature, and for investigating the spreading of species and propagation of infectious diseases, see \cite{L,FLW} for one-single population model with seasonal pulse, \cite{M} for a two-species competition model, \cite{XLS} for a juvenile-adult mutualistic model with harvesting pulse and \cite{LZZ} for predator-prey model. Except fot population problems, pulses have also been introduced to epidemic models, we refer the readers to \cite{D,DA} and references therein for details. In this paper 
we consider an age-structured juvenile-adult model with moving boundaries and harvesting pulse exerting on adults in heterogeneous time-periodic environment, which reads as
\begin{eqnarray}
\left\{
\begin{array}{lll}
J_t-d_1J_{xx}=b(t)A-{a(t)}J-m_1(t)J-\alpha_1(t)J^2,  & t\in((n\tau)^+,(n+1)\tau],\,x \in (g(t),h(t)),\\[1mm]
A_t-d_2A_{xx}={a(t)}J-m_2(t)A-\alpha_2(t)A^2, & t\in((n\tau)^+,(n+1)\tau],\,x \in (g(t),h(t)),\\[1mm]
J((n\tau)^+,x)=J(n\tau,x), & x \in (g(n\tau),h(n\tau)), \\[1mm]
A((n\tau)^+,x)=H(A(n\tau,x)), & x \in (g(n\tau),h(n\tau)), \\[1mm]
J(t,x)=A(t,x)=0, & t\in(0,+\infty),\,x \in \{g(t),h(t)\}, \\[1mm]
h'(t)=-\mu_1J_x(t,h(t))-\mu_2A_x(t,h(t)), &  t\in(n\tau,(n+1)\tau], \\[1mm]
g'(t)=-\mu_1J_x(t,g(t))-\mu_2A_x(t,g(t)), &  t\in(n\tau,(n+1)\tau], \\[1mm]
h(0)=h_0, g(0)=-h_0, \\[1mm]
J(0,x)=J_0(x), A(0,x)=A_0(x),& x\in[-h_0,h_0],n=0,1,2,\dots.\\[1mm]
\end{array} \right.
\label{a01}
\end{eqnarray}

Problem (1.1) can be regarded as 
a discrete-continuous impulsive population model describing the spreading of juveniles and adults over moving range $[g(t),h(t)]$, where $J(t,x)$ and $A(t,x)$ represent densities of juveniles and adults at time $t$ and space $x$, respectively. It shows a dispersal stage and a harvesting stage imposed on the adult $A(t,x)$ with harvesting pulse occurs at every time $t=n\tau,n=0,1,2,\dots$. $H(A(t,x))$ represents the impulsive harvesting function, the right limit $A((n\tau)^+,x)$ denotes the density of adults after harvesting at time $t=n\tau$ and $A((n\tau)^+,x)=H(A(n\tau,x))$ indicates that the population density of adults after harvesting can be regarded as a function of the population density before harvesting at $t=n\tau$. Moreover, different parameters in $H(A)$ will cause different harvesting pulse, so $1-H(A)/A$ is naturally used to characterize the harvesting intensity on adults and we call it harvesting rate in the following.

For the impulsive harvesting function $H(A)\in C^2([0,+\infty))$, we make the following assumptions:\\

($\mathcal{H}$1)\,\,$H(A)>0$ for $A>0$ and $H(0)=0$; $H'(A)>0$ for $A\geq0$. Further, $H(A)/A$ is nonincreasing in $A$ and $0<H(A)/A<1$. Thus harvesting rate $1-H(A)/A \in (0,1)$.

($\mathcal{H}$2)\,\,There exists positive constant $D,\,\gamma>1$ such that $H(A)\geq H'(0)A-DA^\gamma$ for $0\leq A\leq \sigma$ with some small $\sigma$.\\

It is worth noting that ($\mathcal{H}$1) is a natural assumption about impulsive harvesting function, while ($\mathcal{H}$2) is the mathematical condition required for constructing the upper and lower solutions in the sequel, see also in \cite{L}. Definitions of other model parameters in \eqref{a01} are as follows:

$\bullet$ $d_1$ and $d_2$: the diffusive rates of juveniles and adults, respectively;

$\bullet$ $b(t)$: reproduction rate of adults;

$\bullet$ $h(t)$: the rate at which juveniles mature into adults;

$\bullet$ $m_1(t)$ and $m_2(t)$: the death rates of juveniles and adults, respectively;

$\bullet$ $\alpha_1(t)$ and $\alpha_2(t)$: competition coefficients in juvenile and adult individuals.

Moreover, the initial value and positive functions $b(t), a(t), m_1(t), m_2(t), \alpha_1(t), \alpha_2(t)$ meet the following assumptions:\\

($\mathcal{A}$1)\,\, $b(t), a(t), m_1(t), m_2(t), \alpha_1(t),\alpha_2(t)\in C^{1+\alpha/2}([0,+\infty))$ and are $\tau-$periodic in time $t$ for some $\tau > 0$;

($\mathcal{A}$2)\,\, $b(t), a(t), m_1(t), m_2(t), \alpha_1(t)$ and $\alpha_2(t)$ have upper and lower bounds for $t\in[0,\tau]$. For convenience, we denote the maximum and minimum values of these functions corresponding to $t$ by marking superscript $M$ and $m$, respectively. For example, $b^M:=\max_{0\leq t\leq \tau}{b(t)}$ and $b^m:=\min_{0\leq t\leq \tau}{b(t)}$.

($\mathcal{A}$3)\,\, The initial value $(J_0,A_0)$ satisfy
\begin{eqnarray}
\left\{
\begin{array}{lll}
J_0(x),A_0(x)\in C^2([-h_0,h_0]), \\[1mm]
J_0(x),A_0(x)>0\,\, \text{for}\,\, x\in(-h_0,h_0), \\[1mm]
J_0(-h_0)=J_0(h_0)=A_0(-h_0)=A_0(h_0)=0. \\[1mm]
\end{array} \right.
\label{a02}
\end{eqnarray}

Different from pervious work with moving boundary, the introduction of harvesting pulse in such an age-structured impulsive system \eqref{a01} in heterogeneous environment is
more close to reality to allow diversity in results, however it will also cause new challenges in theoretical analysis.

Naturally,
we want to establish the existence and uniqueness of global solution to the impulsive problem \eqref{a01}, and explore the existence of the principal eigenvalue and analyse its properties if harvesting pulse happens. Also, can we establish a criteria for spreading-vanishing under harvesting pulse? Weather or not do harvesting pulse from external intervention and expanding capacities form species itself change the population state? Moreover, if pulse occurs, how do the intensity and timing of harvesting affect or alter the dynamical behaviors of solution? In this paper, we will address all of these questions.

The remainder of this paper is arranged as follows. In the next section, the global existence and uniqueness of solution to problem \eqref{a01} with harvesting pulse is given by extending relevant results of problem \eqref{b01} without pulse. Section 3 focuses on the spreading-vanishing dichotomy, in which the principal eigenvalue related to pulse and some results in fixed domain are firstly given. Section 4 is devoted to some sufficient conditions related to pulse about spreading or vanishing.
The numerical approximations in Section 5 will certify our analytical findings and visually illustrate the transitions in the state of species caused by dual effects of harvesting pulse (harvesting efficacy and timing) and expanding capacities.

\section{Global existence and uniqueness}
Before studying the global existence and uniqueness of free boundaries problem \eqref{a01} with harvesting pulse, we first introduce the following auxiliary moving boundaries problem without pulse
\begin{eqnarray}
\left\{
\begin{array}{lll}
J_t-d_1J_{xx}=b(t)A-a(t)J-m_1(t)J-\alpha_1(t)J^2,  & t\in(0,+\infty),\,x \in (g(t),h(t)),\\[1mm]
A_t-d_2A_{xx}=a(t)J-m_2(t)A-\alpha_2(t)A^2, & t\in(0,+\infty),\,x \in (g(t),h(t)),\\[1mm]
J(t,x)=A(t,x)=0, & t\in(0,+\infty),\,x \in \{g(t),h(t)\}, \\[1mm]
h'(t)=-\mu_1J_x(t,h(t))-\mu_2A_x(t,h(t)), & t\in(0,+\infty), \\[1mm]
g'(t)=-\mu_1J_x(t,g(t))-\mu_2A_x(t,g(t)), & t\in(0,+\infty), \\[1mm]
h(0)=h_0, g(0)=-h_0, \\[1mm]
J(0,x)=J_0(x), A(0,x)=A_0(x),& x\in[-h_0,h_0]. \\[1mm]
\end{array} \right.
\label{b01}
\end{eqnarray}

\begin{thm}
Problem \eqref{b01} admits a unique solution $(J,A;(g,h))$ for all $t>0$. In addition,
$$
(J(t,x),A(t,x))\in [C^{(1+\alpha)/2,1+\alpha}([0,+\infty)\times[g(t),h(t)])]^2
$$
and
$$
(g(t),h(t))\in [C^{1+\alpha/2}([0,+\infty)]^2.
$$
\end{thm}
\bpf
Motivated by previous works of \cite{DL}, the proof is divided into three steps and we give the sketches.

\textbf{Step 1} (local existence and uniqueness) We first straighten the free boundaries $g(t)$ and $h(t)$. Considering the following transformation
$$
x=y+\zeta(y)(h(t)-h_0)+\zeta(-y)(g(t)+h_0),\, -\infty<y<+\infty,
$$
where $\zeta(y)\in C^3[0,+\infty)$ satisfying
$$
\zeta(y)=1\,\,\text{if}\,\, |y-h_0|<h_0/4,\,\, \text{and}\,\, \zeta(0)=0 \,\,\text{if}\,\, |y-h_0|>h_0/2.
$$
the free boundary $x=h(t)$ is transformed to the line $y=h_0$, and $x=g(t)$ becomes the line $y=-h_0$.

Next, a complete metric space and a contraction mapping can be established, and then the local existence and uniqueness result can be derived by applying contracting mapping principle. Moreover, for any $\alpha \in (0,1)$, there exists $\tau_0>0$ such that
$$
(J(t,x),A(t,x))\in [C^{(1+\alpha)/2,1+\alpha}([0,\tau_0]\times[g(t),h(t)])]^2
$$
and
$$
(g(t),h(t))\in [C^{1+\alpha/2}([0,\tau_0])]^2.
$$

\textbf{Step 2} (the estimate of solution) Estimations by the comparison principle yield $0<J(t,x),A(t,x)\leq C^*$ and $0<-g'(t)<C_1,\,0<h'(t)<C_1$ for $t\in(0,\tau_0)$, where $C^*$ and $C_1$  are independent of $\tau_0$.

\textbf{Step 3} (extension of $\tau_0$) According to Zorn's lemma, we assume that the solution exists in the maximal time interval $[0,\tau_{max})$ with $\tau_{max}>0$. The estimates of $(J,A;(g,h))$ in Step 2, which together with the standard continuous extension method, yield the global existence and uniqueness of the solution to problem \eqref{b01}.
\epf\\

In the following, the boundedness, global existence and uniqueness of solution to problem \eqref{a01} with harvesting pulse will be established.
\begin{lem}
Assume that $(J,A;(g,h))$ is a solution to problem \eqref{a01} defined for $t\in(0,\tau^*]$ with $\tau^*:=n_0\tau+\tau_0$ for any given nonnegative integer $n_0$ and $0<\tau_0\leq \tau$. Then we acquire
$$
0<J(t,x),A(t,x)\leq C^*:=\max\{\frac{b^M}{\alpha^m _1},\frac{a^M}{\alpha^m_2},\parallel J_0\parallel_\infty,\parallel A_0\parallel_\infty\}
$$
for $t\in(0,\tau^*]$ and $x\in(g(t),h(t))$. Also,
$$
0<-g'(t)<C_1,\,0<h'(t)<C_1
$$
for $t\in(0,\tau^*]$ and $t\neq n\tau$, where positive constants $C^*$ and $C_1$ are independent of $\tau$ and positive integer $n$.
\end{lem}
\bpf
Let
$$\bar J(t,x)=\bar A(t,x)=C^*:=\max\{\frac{b^M}{\alpha^m _1},\frac{a^M}{\alpha^m_2},\parallel J_0\parallel_\infty,\parallel A_0\parallel_\infty\}.$$
One easily checks that
$$
b(t)\bar A-a(t)\bar J-m_1(t)\bar J-\alpha_1(t)\bar J^2 \leq C^*[b(t)-\alpha_1(t)C^*]<0
$$
and
$$
a(t)\bar J-m_2(t)\bar A-\alpha_2(t)\bar A^2\leq C^*[a(t)-\alpha_2(t)C^*]<0
$$
according to the choice of $C^*$.

When $n_0=0$, we get $t\in (0,\tau^*]\subseteq (0,\tau]$. The initial value satisfies $(\bar J(0,x),\bar A(0,,x))\geq (J_0(x),A_0(x))$, and therefore $(\bar J(0^+,x),\bar A(0^+,,x))\geq (J(0^+,x), A_0(0^+,x))$. It follows from the comparison principle that $(J(t,x),A(t,x))\leq (C^*,C^*)$ for $t\in (0,\tau_0]$ and $x\in(g(t),h(t))$.

When $n_0=1$, we obtain $t\in (0,\tau^*]\subseteq (0,\tau+\tau_0]$. Since the interval $t\in (0,\tau_0]$ is discussed in the above, we here fix $t\in ((\tau_0)^+,\tau+\tau_0]$.
It is clear to see from ($\mathcal{H}$1) that
$$
\bar J((\tau_0)^+,x)=C^*=\bar J(\tau_0,x),\,\, \bar A((\tau_0)^+,x)=C^*>H(C^*)= H(\bar A(\tau_0,x)),
$$
and we can apply the comparison principle to conclude that $(J(t,x),A(t,x))\leq (C^*,C^*)$ for $t\in (0,\tau+\tau_0]$ and $x\in(g(t),h(t))$.

Consequently, by induction of $n_0$, $(J(t,x),A(t,x))\leq (C^*,C^*)$ hold for $t\in (0,\tau^*]:=(0,n_0\tau+\tau_0]$ and $x\in(g(t),h(t))$. It then follows from the strong maximum principle that $0<J(t,x),A(t,x)\leq C^*$ for $t\in (0,\tau^*]$ and $x\in(g(t),h(t))$.

Now, we will give the estimates of $g(t)$ and $h(t)$. Thanks to $J(t,g(t))=J(t,h(t))=A(t,g(t))=A(t,h(t))=0$ and Hopf's boundary lemma to the equations of $J$ and $A$, respectively, we derive that
$$
J_x(t,g(t))>0,\,\,J_x(t,h(t))<0
$$
and
$$
A_x(t,g(t))>0,\,\,A_x(t,h(t))<0
$$
for $t\in(0,\tau^*]$. Since the Stefen type conditions
$h'(t)=-\mu_1J_x(t,h(t))-\mu_2A_x(t,h(t))$ and
$g'(t)=-\mu_1J_x(t,g(t))-\mu_2A_x(t,g(t))$,
we acquire $h'(t)>0$ and $-g'(t)>0$ for $t\in(0,\tau^*]$ and $t\neq n\tau$.

The proof ends with $-g'(t)<C_1$ and $h'(t)<C_1$, which are similar to Lemma 2.2 in \cite{DL}. In particular, it is worth noting that $\Omega_M$ and $w(t,x)$ in \cite{DL} should be modified by
$$
\Omega_{1M}:=\{(t,x): 0<t\leq \tau^*,\,h(t)-M^{-1}<x<h(t)\},
$$
$$
w_1(t,x)=C^*[2M(h(t)-x)-M^2(h(t)-x)^2)]
$$
and
$$
\Omega_{2M}:=\{(t,x): 0<t\leq \tau^*,\,g(t)<x<g(t)+M^{-1}\},
$$
$$
w_2(t,x)=C^*[2M(x-g(t))-M^2(x-g(t))^2)]
$$
in the process of proving $h'(t)<C_1$ and $-g'(t)<C_1$.

The key of this proof is to select some appropriate $M$ such that $w_i$ be the upper solution in $\Omega_{iM}$ to problem \eqref{a01}.

Careful calculations assert
$$
C_1=2MC^*(\mu_1+\mu_2)
$$
with
$$
M=\max\{\frac{b^M}{2d_1},\frac{a^M}{2d_2},\frac{4\parallel J_0\parallel_{C^1[0,h_0]}}{3C^*},\frac{4\parallel A_0\parallel_{C^1[0,h_0]}}{3C^*}\}.
$$
This choice of $M$ implies that the initial value satisfying
\begin{eqnarray*}
\left\{
\begin{array}{lll}
w_i(0,x)\geq J(0,x), \\[1mm]
w_i(0,x)\geq A(0,x) \\[1mm]
\end{array} \right.
\end{eqnarray*}
for $x\in[h_0-M^{-1},h_0]$ if $i=1$ and $x\in[g_0,g_0+M^{-1}]$ if $i=2$, respectively. We can also verify by direct calculations that equations meet the requirements of upper solution. It then follows from the comparison principle that $w_i(t,x)$ is the upper solution in $\Omega_{iM}$ with $i=1,2$, and we have $h'(t)<2MC^*(\mu_1+\mu_2)$, $-g'(t)<2MC^*(\mu_1+\mu_2)$ by using the free boundary conditions.
\epf

In the following, we will give the global existence and uniqueness of solution to problem \eqref{a01} with moving boundaries by using Theorem 2.1 and Lemma 2.2.
\begin{thm}
The solution $(J,A;(g,h))$ to problem \eqref{a01} exists and is unique for all $t\in(0,+\infty)$. Moreover,
$$
(J(t,x),A(t,x))\in [C^{1,2}((n\tau,(n+1)\tau]\times[g(t),h(t)])]^2
$$
and
$$
(g(t),h(t))\in [C((0,+\infty))\cap C^1((n\tau,(n+1)\tau])]^2
$$
with $n=0,1,2,\dots.$
\end{thm}
\bpf
(i) When $n=0$, the time interval becomes $t\in(0^+,\tau]$. We regard $(J(0^+,x),A(0^+,x))$ as the initial value of solution $(J(t,x),A(t,x))$ to problem \eqref{a01}. Since $A_0(x)\in C^2([-h_0,h_0])$ and $H(A)\in C^2([0,+\infty))$, we derive that new initial value satisfies that $J(0^+,x)=J(0,x)\in C^2([-h_0,h_0])$ and $A(0^+,x)=H(A_0(x))\in C^2([-h_0,h_0])$. Recalling Theorem 2.1 and Lemma 2.2, we declare that the solution $(J,A;(g,h))$ to problem \eqref{a01} exists and is unique in $t\in(0^+,\tau]$.  Furthermore, $(J(t,x),A(t,x))\in [C^{1,2}((0,\tau]\times[g(t),h(t)])]^2$ and $(g(t),h(t))\in [C((0,+\infty))\cap C^1((0,\tau])]^2$.

(ii) When $n=1$, the time interval is $t\in(\tau^+,2\tau]$. By the same procedure as (i), the new initial value satisfies $J(\tau^+,x)=J(\tau,x)\in C^2([g(\tau),h(\tau)])$ and $A(\tau^+,x)=H(A(\tau,x))\in C^2([g(\tau),h(\tau)])$. One easily checks that the solution $(J,A;(g,h))$ to problem \eqref{a01} exists and is unique in $t\in(\tau^+,2\tau]$. Also, $(J(t,x),A(t,x))\in [C^{1,2}((\tau,2\tau]\times[g(t),h(t)])]^2$ and $(g(t),h(t))\in [C((0,+\infty))\cap C^1((\tau,2\tau])]^2$.

(iii) When $n=2,3,\dots$, the local existence and uniqueness of the solution can be derived by the same process in interval $t\in(2\tau^+,3\tau]$, $t\in(3\tau^+,4\tau]$,\dots, and step by step, we then find a maximal time interval $[0,\tau_{max})$ with a maximal time $\tau_{max}>0$, such that problem \eqref{a01} admits a unique solution in this maximal interval.

(iv) We now claim that $\tau_{max}=+\infty$ and this ends the proof of global existence and uniqueness of solution to problem \eqref{a01}. If not, we can find a positive integer $n_1$ large enough, such that $\tau_{max}\in (n_1\tau,(n_1+1)\tau]$. If $\tau_{max}\in (n_1\tau,(n_1+1)\tau)$, $(J((n_1\tau)^+,x), A((n_1\tau)^+,x))$ can be regard as a new initial value and we derive that $(J(t,x),A(t,x))\in [C^{1,2}(n_1\tau,(n_1+1)\tau]\times[g(t),h(t)])]^2$ according to the process above. This contradicts to the maximal time of $\tau_{max}$. On the other hand, if $\tau_{max}=(n_1+1)\tau$, then $J(t,x), A(t,x)\in C^{1,2}((n_1\tau)^+,(n_1+1)\tau]\times[g(t),h(t)])$, $J(((n_1+1)\tau)^+,x)=J((n_1+1)\tau,x)\in C^2([g((n_1+1)\tau),h((n_1+1)\tau)])$ and $A(((n_1+1)\tau)^+,x)=H(A((n_1+1)\tau,x))\in C^2([g((n_1+1)\tau),h((n_1+1)\tau)])$. Therefore, \eqref{a01} admits a unique solution in $t\in ((n_1+1)\tau,(n_1+2)\tau]$, which contradict to $\tau_{max}=(n_1+1)\tau$. Thus $\tau_{max}=+\infty$ and the solution is global.
\epf

\section{Criteria for spreading or vanishing}
Some preliminaries about the principal eigenvalue and the main results in a fixed domain are firstly given in Subsection 3.1, which is vital to the research of long-time dynamical behaviors of solution to problem \eqref{a01}. Then the spreading-vanishing dichotomy is established in Subsection 3.2.
\subsection{Results in a fixed domain}
In order to investigate the dynamic behaviors of free boundaries problem \eqref{a01} with impulse, we first consider the corresponding problem in a fixed interval $(l_1, l_2)$
\begin{eqnarray}
\left\{
\begin{array}{lll}
J_t=d_1J_{xx}+b(t)A-a(t)J-m_1(t)J-\alpha_1(t)J^2, & t\in((n\tau)^+,(n+1)\tau],\,\,x \in (l_1,l_2),  \\[1mm]
A_t=d_2A_{xx}+a(t)J-m_2(t)A-\alpha_2(t)A^2,& t\in((n\tau)^+,(n+1)\tau],\,\,x \in (l_1,l_2),  \\[1mm]
J((n\tau)^+,x)=J(n\tau,x),& x\in(l_1,l_2),\\[1mm]
A((n\tau)^+,x)=H(A(n\tau,x)),& x\in(l_1,l_2),\\[1mm]
J(t,l_1)=J(t,l_2)=0,A(t,l_1)=A(t,l_2)=0, & t>0,\\[1mm]
J(0,x)=J_0(x),A(0,x)=A_0(x), & x\in[l_1,l_2],\,n=0,1,2,\dots,\\[1mm]
\end{array} \right.
\label{c01}
\end{eqnarray}
and its corresponding periodic eigenvalue problem involving pulse is
\begin{eqnarray}
\left\{
\begin{array}{lll}
\phi_t-d_1\phi_{xx}=b(t)\psi-[a(t)+m_1(t)]\phi+\lambda_1^D\phi,& t\in(0^+,\tau],\,\,x \in (l_1,l_2),  \\[1mm]
\psi_t-d_2 \psi_{xx}=a(t)\phi-m_2(t)\psi+\lambda_1^D\psi, & t\in(0^+,\tau],\,\,x \in (l_1,l_2),  \\[1mm]
\phi(0^+,x)=\phi(0,x),& x\in (l_1,l_2),\\[1mm]
\psi(0^+,x)=H'(0)\psi(0,x),& x\in (l_1,l_2),\\[1mm]
\phi(t,l_1)=\phi(t,l_2)=\psi(t,l_1)=\psi(t,l_2)=0, & t\in[0,\tau],\\[1mm]
\phi(0,x)=\phi(\tau,x),\psi(0,x)=\psi(\tau,x), & x\in [l_1,l_2].
\end{array} \right.
\label{c02}
\end{eqnarray}

The existence of the principal eigenvalue $\lambda_1^D:=\lambda_1^D((l_1,l_2),H'(0))$ to problem \eqref{c02} in heterogeneous environment is similar to that in \cite{XLS}, in which $M^*$ is replaced by $M_1^*:=b^M+a^M\max\{g'(0),1\}+\ln (\max\{1/g'(0),1\})+1$. Some properties about $\lambda_1^D((l_1,l_2),H'(0))$ and dynamic behaviors of solution are also known.

\begin{lem} The following statements hold.

$(i)$ $\lambda_1^D((l_1,l_2),H'(0))=\lambda_1^D((0,l_2-l_1),H'(0))$;

$(ii)$ $\lambda_1^D$ is strictly decreasing with respect to the domain $\Omega$ for any given $H'(0)$, that is,
$$\lambda_1^D(\Omega_2,H'(0))<\lambda_1^D(\Omega_1,H'(0))$$
hold for any $\Omega_1\subseteq \Omega_2\subseteq R^1$ and $\Omega_2 \backslash \Omega_1$ is nonempty. Also, $\lambda_1^D$ is strictly decreasing in $H'(0)$ for any given $\Omega$;

$(iii)$ Suppose that $\lambda_1^D>0$. The solution $(J(t,x),A(t,x))$ to problem \eqref{c01} satisfies
$$
\lim\limits_{t\to+\infty}{(J(t,x),A(t,x))}=(0,0)
$$
uniformly for $x\in [l_1,l_2]$;

$(iv)$ Assume that $\lambda_1^D<0$. Then for nonnegative nontrivial initial value $(J_0(x),A_0(x))$, any solution $(J(t,x),A(t,x))$ to problem \eqref{c01} admits
$$\lim\limits_{n\to+\infty}{(J,A)}(t+n\tau,x)=(u^*,v^*)(t,x),\,t\in[0,+\infty),\,x\in[l_1,l_2],$$
where $(u^*(t,x),v^*(t,x))$ is the unique positive $\tau-$periodic solution satisfying
\begin{eqnarray*}
\left\{
\begin{array}{lll}
U_t=d_1U_{xx}+b(t)V-a(t)U-m_1(t)U-\alpha_1(t)U^2, & t\in(0^+,\tau],\,\,x \in (l_1,l_2),  \\[1mm]
V_t=d_2V_{xx}+a(t)U-m_2(t)V-\alpha_2(t)V^2,& t\in(0^+,\tau],\,\,x \in (l_1,l_2),  \\[1mm]
U(0^+,x)=U(0,x),& x\in(l_1,l_2),\\[1mm]
V(0^+,x)=H(V(\tau,x)),& x\in(l_1,l_2),\\[1mm]
U(t,l_1)=U(t,l_2)=0,\,V(t,l_1)=V(t,l_2)=0, & t>0,\\[1mm]
U(0,x)=U(\tau,x),\,V(0,x)=V(\tau,x), & x\in[l_1,l_2],\,n=0,1,2,\dots.\\[1mm]
\end{array} \right.
\end{eqnarray*}
\label{xuhaiyan}
\end{lem}
\bpf
(i) It follows from the transformation $(t,x)\rightarrow (t,z)$ with $z=x-l_1$ and $l_1<x<l_2$ in \eqref{c02} that $\lambda_1^D((l_1,l_2),H'(0))=\lambda_1^D((0,l_2-l_1),H'(0))$.

(ii) Let
\begin{eqnarray*}
\left\{
\begin{array}{lll}
\phi(t,x)=\alpha(t)\Phi(x),\\[1mm]
\psi(t,x)=\beta(t)\Phi(x)\\[1mm]
\end{array} \right.
\end{eqnarray*}
with $\Phi(x)$ satisfying
\begin{eqnarray*}
\left\{
\begin{array}{lll}
-\Phi_{xx}=\lambda_0\Phi(x),\,&x\in(l_1,l_2),\\[1mm]
\Phi(x)=0,\,\,&x\in \{l_1,l_2\}.\\[1mm]
\end{array} \right.
\end{eqnarray*}
One easily checks that $\lambda_0:=\lambda_0(l_1,l_2)=(\frac{\pi}{l_2-l_1})^2$, and corresponding positive eigenfunction of $\lambda_0$ is $\Phi(x)=\cos[\frac{\pi}{l_2-l_1}(x-\frac{l_1+l_2}{2})]$.

By variable separation approach method, the principal eigenvalue problem \eqref{c02} is turned into
\begin{eqnarray}
\left\{
\begin{array}{lll}
\alpha'(t)+d_1\lambda_0\alpha(t)=b(t)\beta(t)-[a(t)+m_1(t)]\alpha(t)+\lambda_1^D\alpha(t),& t\in(0^+,\tau],\\[1mm]
\beta'(t)+d_2\lambda_0\beta(t)=a(t)\alpha(t)-m_2(t)\beta(t)+\lambda_1^D\beta(t),& t\in(0^+,\tau],\\[1mm]
\alpha(0^+)=\alpha(0),\,\beta(0^+)=H'(0)\beta(0),\\[1mm]
\alpha(0)=\alpha(\tau),\,\beta(0)=\beta(\tau).
\end{array} \right.
\label{c03}
\end{eqnarray}
Recalling that $\lambda_0(l_1,l_2):=\lambda_0(0,l_2-l_1)$ is nonincreasing with respect to domain $\Omega:=[l_1,l_2]$, it then follows from \eqref{c03} yields that $\lambda_1^D((0,l_2-l_1),H'(0))$ is decreasing with respect to $(l_2-l_1)$. Finally, the monotonicity of $\lambda_1^D$ with respect to $H'(0)$ has been proved in \cite{XLS} (Theorem 3.3).

The long-time behaviors (iii) and (iv) in Lemma 3.1 for heterogeneous environment can also be derived by the same manner in \cite{XLS} (Theorems 4.2 and 4.5) and we omit it here.

\epf

\vspace{3mm}
The same as the principal eigenvalue $\lambda_1^D$ in problem \eqref{c01} in a fixed domain, we here introduce the threshold value of free boundaries problem \eqref{a01} by
$$
\lambda_0^F(t):=\lambda_1^D((g(t),h(t)),H'(0)).
$$
Since $-g(t)$ and $h(t)$ are strictly nondecreasing with respect to $t$, it then follows from Lemma 3.1 (ii) that $\lambda_0^F(t)$ is strictly decreasing in $t$. We then can define
$$
\lambda_1^D((-\infty,+\infty),H'(0)):=\lim\limits_{l\to\infty}\lambda_1^D((-l,l),H'(0))
$$
which is related to impulsive harvesting $H'(0)$.

In the following, the comparison principle of mutualistic and impulsive model is firstly set up, and then the criteria for spreading or vanishing of species will be investigated around the principal eigenvalue defined above.
\subsection{Spreading-vanishing dichotomy}
In view of $g'(t)<0$ and $h'(t)>0$ in Lemma 2.2, we declare

$$
\lim\limits_{t\to+\infty}g(t)=g_\infty\,\,\text{and}\,\,\lim\limits_{t\to+\infty}h(t)=h_\infty,
$$
with some $g_\infty\in [-\infty,-h_0)$ and $h_\infty\in (h_0,+\infty]$.

To begin with, we first introduce the comparison principle, which is the basis of the forthcoming theoretical analysis.
\begin{lem}
{\rm (Comparison principle)} Assume that $(\bar J(t,x),\bar A(t,x))\in [C^{1,2}((n\tau,(n+1)\tau]\times[\bar g(t),\bar h(t)])]^2$, $(\bar g(t),\bar h(t))\in [C((0,+\infty))\cap C^1((n\tau,(n+1)\tau])]^2$,
and
\begin{eqnarray}
\left\{
\begin{array}{lll}
\bar J_t-d_1\bar J_{xx}\geq b(t) \bar A-(a(t)+m_1(t))\bar J-\alpha_1(t)\bar J^2,  & t\in((n\tau)^+,(n+1)\tau],\,x \in (\bar g(t),\bar h(t)),\\[1mm]
\bar A_t-d_2\bar A_{xx}\geq a(t)\bar J-m_2(t)\bar A-\alpha_2(t)\bar A^2, & t\in((n\tau)^+,(n+1)\tau],\,x \in (\bar g(t),\bar h(t)),\\[1mm]
\bar J((n\tau)^+,x)\geq \bar J(n\tau,x), & x \in (\bar g(n\tau),\bar h(n\tau)), \\[1mm]
\bar A((n\tau)^+,x)\geq H(\bar A(n\tau,x)), & x \in (\bar g(n\tau),\bar h(n\tau)), \\[1mm]
\bar J(t,x)=\bar A(t,x)=0, & t\in(0,+\infty),\,x \in \{\bar g(t),\bar h(t)\}, \\[1mm]
\bar h'(t)\geq-\mu_1\bar J_x(t,\bar h(t))-\mu_2\bar A_x(t,\bar h(t)), & t\in(n\tau,(n+1)\tau], \\[1mm]
\bar g'(t)\leq-\mu_1\bar J_x(t,\bar g(t))-\mu_2\bar A_x(t,\bar g(t)), & t\in(n\tau,(n+1)\tau], \\[1mm]
\bar g(0)\leq -h_0<h_0\leq\bar h(0), \\[1mm]
\bar J(0,x)\geq J_0(x), \bar A(0,x)\geq A_0(x),& x\in[-h_0,h_0],n=0,1,2,\dots.\\[1mm]
\end{array} \right.
\label{c04}
\end{eqnarray}
Then the solution $(J,A;(g,h))$ to problem \eqref{a01} satisfies
$$(J(t,x),A(t,x))\leq(\bar J(t,x),\bar A(t,x))\,\,\text{for}\,\, t\in(0,+\infty)\,\,\text{and}\,\,x\in[g(t),h(t)],$$
$$\bar g(t)\leq g(t),\,h(t)\leq \bar h(t)\,\,\text{for}\,\, t\in(0,+\infty).$$
\end{lem}

\begin{rmk} We call the pair $(\bar J,\bar A;(\bar g,\bar h))$ in $Lemma \ 3.2$ an upper solution to the free boundaries problem \eqref{a01}. Similarly, we can define a lower solution to problem \eqref{a01} by $(\underline{J},\underline{A};(\underline{g},\underline{h}))$ if inequalities in \eqref{c04} are all reversed.
\end{rmk}

\begin{thm}
Assume that $\lambda_1^D((-\infty,+\infty), H'(0))\geq0$. Then $-\infty<g_\infty<h_\infty<+\infty$ and
$$\lim\limits_{t\to+\infty}{||(J(t,\cdot),A(t,\cdot))||_{C[g(t),h(t)]}}=(0,0).$$
\end{thm}
\bpf An upper-lower solution technique will be used to deduce the desired conclusion. Define
$$
\hat \eta(t)= h\omega(t),\,\,\omega(t):=1+\theta-\frac{\theta}{2}e^{-\gamma t},\,\,t\geq0,
$$
$$
\hat W_1(t,x)=Ce^{-\gamma t}\phi(t,\frac{hx}{\hat \eta(t)}),\,\,\hat W_2(t,x)=Ce^{-\gamma t}\psi(t,\frac{hx}{\hat \eta(t)}),\,\,t\geq0,\,\,x\in[-h,h],
$$
where positive constants $h(\geq h_0)$, $\gamma$, $C$ and $\theta$ to be chosen later. $(\phi(t,x),\psi(t,x))$ is a pair of eigenfunction corresponding to the principal eigenvalue $\lambda_1^D((-h,h),H'(0))$ in \eqref{c02} with $\lambda_1^D((-h,h),H'(0))>0$. Recall that $\phi_x(t,-h)>0$, $\phi_x(t,h)<0$, $\psi_x(t,-h)>0$ and $\psi_x(t,h)<0$, it is easy to verify that $x\phi_x\leq0$ and $x\psi_x\leq0$ in $x\in[-\hat \eta(t),\hat \eta(t)]$. Moreover, in view of Section 3.1, we obtain
$$
\phi_{xx}=-(\frac{\pi}{2h})^2\phi<0\,\,\text{and}\,\,\psi_{xx}=-(\frac{\pi}{2h})^2\psi<0.
$$

Careful calculations in $t\in((n\tau)^+,(n+1)\tau]$ and $x\in (-\hat \eta(t),\hat \eta(t))$ yield
$$\begin{array}{llllll}
&&\hat W_{1t}-d_1\hat W_{1xx}-b(t)\hat W_2+[a(t)+m_1(t)]\hat W_1+\alpha_1(t)\hat W_1^2\\[1mm]
&=&-\gamma Ce^{-\gamma t}\phi+Ce^{-\gamma t}\phi_t-Ce^{-\gamma t}x\phi_x\omega'\omega^{-2}-d_1Ce^{-\gamma t}\phi_{xx}\omega^{-2}\\[1mm]
&-&b(t)Ce^{-\gamma t}\psi+[a(t)+m_1(t)]Ce^{-\gamma t}\phi+\alpha_1(t)(Ce^{-\gamma t})^2(\phi)^2\\[1mm]
&\geq&Ce^{-\gamma t}[\phi_t-d_1\phi_{xx}\omega^{-2}-b(t)\psi+(a(t)+m_1(t))\phi-\gamma \phi]\\[1mm]
&=&Ce^{-\gamma t}[(\phi_t-d_1\phi_{xx})+(d_1\phi_{xx}-d_1\phi_{xx}\omega^{-2})-b(t)\psi+(a(t)+m_1(t))\phi-\gamma \phi]\\[1mm]
&=&Ce^{-\gamma t}[(\lambda_1^D((-h,h),H'(0))-\gamma)\phi+d_1\phi_{xx}(1-\omega^{-2})]\\[1mm]
&=&Ce^{-\gamma t}\phi[(\lambda_1^D((-h,h),H'(0))-\gamma)-d_1(\frac{\pi}{2h})^2(1-\omega^{-2})].
\end{array}$$
Further, a positive constant $\gamma:=\frac{\lambda_1^D((-h,h),H'(0))}{2}$ can be selected such that $\lambda_1^D((-h,h),H'(0))-\gamma=\frac{\lambda_1^D((-h,h),H'(0))}{2}>0$. Since $\lim\limits_{\theta\to 0^+}\omega(t)=\lim\limits_{\theta\to 0^+}(1+\theta-\frac{\theta}{2}e^{-\gamma t})=1$, one can choose $\theta$ small enough such that $d_1(\frac{\pi}{2h})^2(1-\omega^{-2})<\frac{\lambda_1^D((-h,h),H'(0))}{2}$, and then
$$\hat W_{1t}-d_1\hat W_{1xx}-b(t)\hat W_2+[a(t)+m_1(t)]\hat W_1+\alpha_1(t)\hat W_1^2\geq 0.$$

Similarly, It is clear that
$$\begin{array}{llllll}
&&\hat W_{2t}-d_2\hat W_{2xx}-a(t)\hat W_1+m_2(t)\hat W_2+\alpha_2(t)\hat W_2^2\\[1mm]
&\geq&Ce^{-\gamma t}[\psi_t-d_2\psi_{xx}\omega^{-2}-a(t)\phi++m_2(t)\psi-\gamma \psi]\\[1mm]
&=&Ce^{-\gamma t}[(\psi_t-d_2\psi_{xx})+(d_2\psi_{xx}-d_2\psi_{xx}\omega^{-2})-a(t)\phi+m_2(t)\psi-\gamma \psi]\\[1mm]
&=&Ce^{-\gamma t}[(\lambda_1^D((-h,h),H'(0))-\gamma)\psi+d_2\psi_{xx}(1-\omega^{-2})]\\[1mm]
&=&Ce^{-\gamma t}\psi[(\lambda_1^D((-h,h),H'(0))-\gamma)-d_2(\frac{\pi}{2h})^2(1-\omega^{-2})]\\[2mm]
&\geq& 0.
\end{array}$$
Meanwhile, for $x\in[-\hat\eta(n\tau),\hat\eta(n\tau)]$, we get $$\hat W_1((n\tau)^+,x)=Ce^{-\gamma n\tau}\phi((n\tau)^+,\frac{hx}{\hat \eta(n\tau)})=Ce^{-\gamma n\tau}\phi(n\tau,\frac{hx}{\hat \eta(n\tau)})=\hat W_1(n\tau,x)$$ and
$$\begin{array}{llllll}
&&\hat W_2((n\tau)^+,x)=Ce^{-\gamma n\tau}\psi((n\tau)^+,\frac{hx}{\hat \eta(n\tau)})\\[1mm]
&=&Ce^{-\gamma n\tau}H'(0)\psi(n\tau,\frac{hx}{\hat \eta(n\tau)})\\[1mm]
&=&H'(0)\hat W_2(n\tau,x)\\[1mm]
&\geq&H(\hat W_2(n\tau,x))
\end{array}$$
according to ($\mathcal{H}$1). The initial value satisfies $\hat W_1(0,x)=C\phi(0,\frac{x}{1+\frac{\theta}{2}})\geq J_0(x)$ and $\hat W_2(0,x)=C\psi(0,\frac{x}{1+\frac{\theta}{2}})\geq A_0(x)$ in $[-h_0,h_0]$ for sufficiently big $C$. Also, for the boundary conditions, $\hat W_1(t,-\hat \eta(t))=\hat W_1(t,\hat \eta(t))=\hat W_2(t,-\hat \eta(t))=\hat W_2(t,\hat \eta(t))=0$ hold for $t>0$.

Therefore, $(\hat W_1,\hat W_2;(-\hat\eta,\hat\eta))$ is an upper solution to problem \eqref{a01} provided that
\begin{equation}
-\mu_1 \hat W_{1x}(t,\hat\eta(t))-\mu_2\hat W_{2x}(t,\hat\eta(t))\leq \hat\eta'(t)
\label{f01}
\end{equation}
and
\begin{equation}
-\mu_1 \hat W_{1x}(t,-\hat\eta(t))-\mu_2\hat W_{2x}(t,-\hat\eta(t))\geq -\hat\eta'(t)
\label{f02}
\end{equation}
for $t\in (n\tau,(n+1)\tau]$.

Since $\omega^{-1}(t)\leq 2/(2+\theta)$, $\phi_x(t,h)=-\frac{\pi}{2h}\alpha(t)<0$ and $\psi_x(t,h)=-\frac{\pi}{2h}\beta(t)<0$, careful calculations yield
$$\begin{array}{llllll}
&&-\mu_1 \hat W_{1x}(t,\hat\eta(t))-\mu_2\hat W_{2x}(t,\hat\eta(t))\\[1mm]
&=&Ce^{-\gamma t}\omega^{-1}(t)[-\mu_1\phi_x(t,h)-\mu_2\psi_x(t,h)]\\[1mm]
&<&\frac{2Ce^{-\gamma t}}{2+\theta}[-\mu_1\phi_x(t,h)-\mu_2\psi_x(t,h)]\\[1mm]
&\leq&\frac{C\pi e^{-\gamma t}}{(2+\theta)h}[\mu_1 \max_{0\leq t\leq \tau}\alpha(t)+\mu_2 \max_{0\leq t\leq \tau}\beta(t)],
\end{array}$$
where $\alpha(t)$ and $\beta(t)$ satisfy \eqref{c03}. Similarly,
$$
-\mu_1 \hat W_{1x}(t,-\hat\eta(t))-\mu_2\hat W_{2x}(t,-\hat\eta(t))\geq \frac{-C\pi e^{-\gamma t}}{(2+\theta)h}[\mu_1 \max_{0\leq t\leq \tau}\alpha(t)+\mu_2 \max_{0\leq t\leq \tau}\beta(t)].
$$
 we then can choose $h=\max\{h_0, \sqrt{\frac{2C\pi}{\gamma\theta(2+\theta)h}[\mu_1 \max_{0\leq t\leq \tau}\alpha(t)+\mu_2 \max_{0\leq t\leq \tau}\beta(t)]}\}$ such that \eqref{f01} and \eqref{f02} hold.

Therefore, $(\hat W_1,\hat W_2;(-\hat \eta,\hat\eta))$ is an upper solution to problem \eqref{a01}. It follows from a simple comparison principle that $(J(t,x),A(t,x))\leq (\hat W_1,\hat W_2)$ in $(t,x)\in(0,+\infty)\times(g(t),h(t))$, and $h(t)\leq \hat\eta(t)$, $g(t)\geq -\hat\eta(t)$ for $t\in(0,+\infty)$. Since $\lim\limits_{t\to+\infty}{(\hat W_1,\hat W_2)}=(0,0)$, we derive that $\lim\limits_{t\to+\infty}{||(J(t,\cdot),A(t,\cdot))||_{C[g(t),h(t)]}}=(0,0).$ Meanwhile, $h_\infty=\lim\limits_{t\to+\infty}{h(t)}\leq\lim\limits_{t\to+\infty}{\hat\eta(t)}=h(1+\theta)<+\infty$ and
$g_\infty=\lim\limits_{t\to+\infty}{g(t)}\geq\lim\limits_{t\to+\infty}{-\hat\eta(t)}=-h(1+\theta)>-\infty$ can be obtained. This finishes the proof.
\epf

\begin{thm}
If $-\infty<g_\infty<h_\infty<+\infty$. Then $\lambda_1^D((g_\infty,h_\infty),H'(0))\geq0$ and
$$\lim\limits_{t\to+\infty}{||(J(t,\cdot),A(t,\cdot))||_{C[g(t),h(t)]}}=(0,0).$$
\end{thm}
\bpf
\textbf{Step 1}\,\,\,We first prove $\lambda_1^D(g_\infty,h_\infty,H'(0))\geq0$.

Otherwise, we assume that $\lambda_1^D((g_\infty,h_\infty),H'(0))<0$. In view of $\lambda_1^D((g_\infty,h_\infty),H'(0))<0$, $h(t)\rightarrow h_\infty$ and $g(t)\rightarrow g_\infty$ as $t\rightarrow \infty$, a small positive constant $\varepsilon$ could be chosen such that
$$\lambda_1^D((g_\infty+\varepsilon,h_\infty-\varepsilon),H'(0))<0$$
by continuous dependence of $\varepsilon$,
and a large integer $N_0:=N_0(\varepsilon)$ could be chosen such that
$$
g_\infty<g(t)<g_\infty+\varepsilon,\ h_\infty-\varepsilon<h(t)<h_\infty, \forall\,\, t\geq N_0\tau
$$
and

For $n\geq N_0$, let $(\underline{J},\underline{A})(t,x)$ be a solution to the following problem
\begin{eqnarray}
\left\{
\begin{array}{lll}
\underline{J}_t-d_1\underline J_{xx}=b(t)\underline A-(a(t)+m_1(t))\underline J-\alpha_1(t)\underline J^2,  &(t,x)\in \Omega_{\tau,\varepsilon},\\[1mm]
\underline A_t-d_2\underline A_{xx}=a(t)\underline J-m_2(t)\underline A-\alpha_2(t)\underline A^2, & (t,x)\in \Omega_{\tau,\varepsilon},\\[1mm]
\underline J((n\tau)^+,x)=\underline J(n\tau,x),&x \in (g_\infty+\varepsilon,h_\infty-\varepsilon),\\[1mm]
\underline A((n\tau)^+,x)=H(\underline A(n\tau,x)), &x \in (g_\infty+\varepsilon,h_\infty-\varepsilon),\\[1mm]
\underline J(t,g_\infty+\varepsilon)=\underline J(t,h_\infty-\varepsilon)=0, & t\geq N_0\tau,\\[1mm]
\underline A(t,g_\infty+\varepsilon)=\underline A(t,h_\infty-\varepsilon)=0, & t\geq N_0\tau,\\[1mm]
\underline J(N_0\tau,x)=J(N_0\tau,x),\,\underline A(N_0\tau,x)=A(N_0\tau,x),& x\in(g_\infty+\varepsilon,h_\infty-\varepsilon),\\[1mm]
\end{array} \right.
\label{c05}
\end{eqnarray}
where $\Omega_{\tau,\varepsilon}:=\{(t,x):(n\tau)^+<t\leq(n+1)\tau],\,g_\infty+\varepsilon<x<h_\infty-\varepsilon\}$.

A comparison principle indicates
\begin{equation}
(J,A)(t,x)\geq(\underline J,\underline A)(t,x),\,\,\forall t\geq N_0\tau, x\in(g_\infty+\varepsilon,h_\infty-\varepsilon).
\label{c06}
\end{equation}

Recalling $\lambda_1^D((g_\infty+\varepsilon,h_\infty-\varepsilon),H'(0))<0$ and Lemma 3.1 (iv), one can see that $(\underline J,\underline A)$ converges to the unique positive steady state corresponding to problem \eqref{c05} with initial value condition replaced by periodic condition, which is denoted by $(U^\varepsilon,V^\varepsilon)$, and
\begin{equation}
\lim\limits_{n\to+\infty}{(\underline J(t+n\tau,x),\underline A(t+n\tau,x))}=(U^\varepsilon(t,x),V^\varepsilon(t,x)),\,t\in[0,\tau],\,x\in(g_\infty+\varepsilon,h_\infty-\varepsilon).
\label{c07}
\end{equation}
It follows form \eqref{c06} and \eqref{c07} that
\begin{equation}
\lim\limits_{n\to+\infty}\inf{(J(t+n\tau,x),A(t+n\tau,x))}\geq(U^\varepsilon(t,x),V^\varepsilon(t,x)),\,t\in[0,\tau],\, x\in[g_\infty+\varepsilon,h_\infty-\varepsilon].
\label{c071}
\end{equation}

Similarly, let $(\bar J,\bar A)$ be a solution of
\begin{eqnarray}
\left\{
\begin{array}{lll}
\bar{J}_t-d_1\bar J_{xx}=b(t)\bar A-(a(t)+m_1(t))\bar J-\alpha_1(t)\bar J^2,  & t\in((n\tau)^+,(n+1)\tau],\,x \in (g_\infty,h_\infty),\\[1mm]
\bar A_t-d_2\bar A_{xx}=a(t)\bar J-m_2(t)\bar A-\alpha_2(t)\bar A^2, & t\in((n\tau)^+,(n+1)\tau],\,x \in (g_\infty,h_\infty),\\[1mm]
\bar J((n\tau)^+,x)=\bar J(n\tau,x),&x \in (g_\infty,h_\infty),\\[1mm]
\bar A((n\tau)^+,x)=H(\bar A(n\tau,x)), &x \in (g_\infty,h_\infty),\\[1mm]
\bar J(t,g_\infty)=\bar J(t,h_\infty)=0, & t\geq 0, \\[1mm]
\bar A(t,g_\infty)=\bar A(t,h_\infty)=0, & t\geq 0, \\[1mm]
\bar J(0,x)=\tilde J_0(x),\,\bar A(0,x)=\tilde A_0(x),& x\in[g_\infty,h_\infty], \\[1mm]
\end{array} \right.
\label{c08}
\end{eqnarray}
where $\tilde J_0(x)=J_0(x)$ in $[-h_0,h_0]$, $\tilde J_0(x)=0$ in $[g_\infty,-h_0)\bigcup (h_0,h_\infty]$ and $\tilde A_0(x)=A_0(x)$ in $[-h_0,h_0]$, $\tilde A_0(x)=0$ in $[g_\infty,-h_0)\bigcup (h_0,h_\infty]$.

Comparison principle arrives at
$$(J,A)\leq(\bar J,\bar A),\,\,\forall t\geq 0, x\in[g(t),h(t)].$$
By $\lambda_1^D((g_\infty,h_\infty),H'(0))<0$ and Lemma 3.1 (iv) in a fixed domain, we have
$$\lim\limits_{n\to+\infty}{(\bar J(t+n\tau,x),\bar A(t+n\tau,x))}=(U(t,x),V(t,x)),\,t\in[0,\tau],\, x\in[g_\infty,h_\infty],$$
where $(U(t,x),V(t,x))$ is the unique solution to periodic problem corresponding to \eqref{c08},
which means
\begin{equation}
\lim\limits_{n\to+\infty}\sup{(J(t+n\tau,x),A(t+n\tau,x))}\leq(U(t,x),V(t,x)),\,t\in[0,\tau],\,x\in(g_\infty,h_\infty).
\label{c081}
\end{equation}

A standard compactness and regularity of partial differential equation to conclude
$$
(U^\varepsilon(t,x),V^\varepsilon(t,x))\rightarrow (U(t,x),V(t,x))\,\, \text{in}\,\, [C_{loc}^{(1+\theta)/2,1+\theta}([0,\tau]\times (g_\infty,h_\infty))]^2\,\,\text{as}\,\,\varepsilon\rightarrow 0,
$$
which together with \eqref{c071} and \eqref{c081}, yield
$$\lim\limits_{n\to+\infty}{(J(t+n\tau,x),A(t+n\tau,x))}=(U(t,x),V(t,x)),\,t\in[0,\tau],\, x\in(g_\infty,h_\infty).$$

Additionally, it is easily seen that
$$
(J_x(t+n\tau,g(t+n\tau)),A_x(t+n\tau,g(t+n\tau)))\rightarrow (U_x(t,g_\infty),V_x(t,g_\infty))>(0,0) \,\,\text{as}\,\,n\rightarrow\infty
$$
and
$$
(J_x(t+n\tau,h(t+n\tau)),A_x(t+n\tau,h(t+n\tau)))\rightarrow (U_x(t,h_\infty),V_x(t,h_\infty))<(0,0) \,\,\text{as}\,\,n\rightarrow\infty
$$
uniformly on $[0,\tau]$, which together with
$$g'(t+n\tau)=-\mu_1J_x(t+n\tau,g(t+n\tau))-\mu_2A_x(t+n\tau,g(t+n\tau))<0$$
and
$$h'(t+n\tau)=-\mu_1J_x(t+n\tau,h(t+n\tau))-\mu_2A_x(t+n\tau,h(t+n\tau))>0$$
for $t\in(n\tau,(n+1)\tau]$, assert
$\lim\limits_{t\to+\infty}-g(t)=\lim\limits_{t\to+\infty}h(t)=\infty$. This contradicts to $-\infty<g_\infty<h_\infty<+\infty$. Consequently,
$\lambda_1^D((g_\infty,h_\infty),H'(0))\geq0$.

\textbf{Step 2}\,\,\,Recalling $\lambda_1^D((g_\infty,h_\infty),H'(0))\geq0$ and Lemma 3.1(iii) in a fixed domain, we have $$\lim\limits_{t\to+\infty}{||(\bar J(t,\cdot),\bar A(t,\cdot))||_{C[g_\infty,h_\infty]}}=(0,0).$$
Then using a comparison principle in Step 1 yields
$$(J,A)(t,x)\leq(\bar J,\bar A)(t,x)\rightarrow (0,0)\,\,\,\text{uniformly\,\, for}\,\, x\in[g(t),h(t)],\,\,\text{as}\,\,t\rightarrow\infty.$$
The proof is now complete.
\epf

\begin{thm}
Suppose that $\lambda_1^D((-\infty,+\infty),H'(0))<0$ and $-g_\infty=h_\infty=\infty$. Then
\begin{equation}
\lim\limits_{n\to+\infty}{(J(t+n\tau,x),A(t+n\tau,x))}=(J^*(t),A^*(t))
\label{c00}
\end{equation}
locally uniformly in $[0,\tau]\times(-\infty,+\infty)$, where $(J^*(t),A^*(t))$ is the unique positive solution to problem
\begin{eqnarray}
\left\{
\begin{array}{lll}
J_t=b(t)A-(a(t)+m_1(t)) J-\alpha_1(t)J^2,  & t\in(0^+,\tau],\\[1mm]
A_t=a(t)J-m_2(t) A-\alpha_2(t) A^2, & t\in(0^+,\tau],\\[1mm]
J(0^+)=J(0),A(0^+)=H(A(0)),\\[1mm]
J(0)=J(\tau),A(0)=A(\tau). \\[1mm]
\end{array} \right.
\label{c11}
\end{eqnarray}
\end{thm}
\bpf
For clarify, the proof is carried out in three steps.

\textbf{Step 1}\,\,\,The proof of existence and uniqueness of solution $(J^*(t),A^*(t))$ to problem \eqref{c11}.

It is obvious that the upper and lower solution to problem \eqref{c11} exist, and we denote them by $(\bar J(t), \bar A(t))$ and $(\underline{J}(t),\underline{A}(t))$, respectively. Moreover, there exists $J_0(x)$ and $A_0(x)$ such that $\underline{J}(0)\leq J_0(x)\leq \bar J(0)$ and $\underline{A}(0)\leq A_0(x)\leq \bar A(0)$.

Next, iteration sequences $\{(\underline{J}^{(n)},\underline{A}^{(n)})\}$ and $\{(\bar J^{(n)}, \bar A^{(n)})\}$ with initial value $(\underline{J}^{(0)},\underline{A}^{(0)})=(\underline{J},\underline{A})$ and $(\bar{J}^{(0)},\bar{A}^{(0)})=(\bar{J},\bar{A})$ can be constructed, which satisfying
\begin{eqnarray}
\left\{
\begin{array}{lll}
U_t^{(n)}+K^*U^{(n)}=K^* U^{(n-1)}+b(t)V^{(n-1)}-(a(t)+m_1(t))U^{(n-1)}\\[1mm]
\qquad \qquad \quad \quad\quad-\alpha_1(t) (U^{(n-1)})^2,  & t\in(0^+,\tau],\\[1mm]
V_t^{(n)}+K^*V^{(n)}=K^* V^{(n-1)}+a(t)U^{(n-1)}-m_2(t)V^{(n-1)}\\[1mm]
\qquad \qquad \qquad \quad-\alpha_2(t)(V^{(n-1)})^2, & t\in(0^+,\tau],\\[1mm]
U^{(n)}(0^+)=U^{(n-1)}(\tau),\,V^{(n)}(0^+)=H(V^{(n-1)}(\tau)),\\[1mm]
U^{(n)}(0)=U^{(n-1)}(\tau),\,V^{(n)}(0)=\bar V^{(n-1)}(\tau)\\[1mm]
\end{array} \right.
\label{c111}
\end{eqnarray}
with $(U(t),V(t))$ replaced by $(\bar J(t), \bar A(t))$ and $(\underline{J}(t),\underline{A}(t))$, respectively. Here, $K^*$ is chosen large enough such that $K^* J+b(t) A-(a(t)+m_1(t))J-\alpha_1(t) J^2$ and $K^* A+a(t)J-m_2(t) A-\alpha_2(t) A^2$ are nondecreasing with respect to $J$ and $A$, respectively. By \cite{P}, we get
$$\lim\limits_{n\to+\infty}{(\underline{J}^{(n)},\underline{A}^{(n)})(t)}=(\underline{J}^*,\underline{A}^*)(t)\,\,\text{and}\,\,
\lim\limits_{n\to+\infty}{(\bar{J}^{(n)},\bar{A}^{(n)})(t)}=(\bar{J}^*,\bar{A}^*)(t)\,\,\text{for}\,\,t\in[0,\tau]
$$
where $(\underline{J}^*,\underline{A}^*)(t)$ and $(\bar{J}^*,\bar{A}^*)(t)$ are the minimal and maximal positive periodic solutions to corresponding periodic problem, respectively.

Step 1 ends with the uniqueness of solution to mutualist problem \eqref{c11}. To prove the uniqueness, let $(\underline{J}^*,\underline{A}^*)(t)$ defined above and $(J_1, A_1)(t)$ be another solution, then denote
$$
S=\{s\in[0,1],sJ_1(t)\leq{\underline{J}^*(t)},\,sA_1(t)\leq{\underline{A}^*(t)},\,t\in[0,\tau]\}.
$$
We claim that $1\in S$. Otherwise, we assume that $s_0=\sup S<1$.
It follows from careful calculations that
$$\begin{array}{llllll}
&&(\underline{J}^*-s_0J_1)_t=b(t)(\underline{A}^*-s_0A_1)-(a(t)+m_1(t))(\underline{J}^*-s_0J_1)-\alpha_1(t)((\underline{J}^*)^2-s_0J_1^2)\\[1mm]
&\geq&-(a(t)+m_1(t))(\underline{J}^*-s_0J_1)-\alpha_1(t)(\underline{J}^*+s_0J_1)(\underline{J}^*-s_0J_1)\\[1mm]
&=&-K_1^*(t)(\underline{J}^*-s_0J_1),
\end{array}$$
where $-K_1^*(t):=-(a(t)+m_1(t))-\alpha_1(t)(\underline{J}^*+s_0J_1)$.

Similarly,
$$\begin{array}{llllll}
&&(\underline{A}^*-s_0A_1)_t=a(t)(\underline{J}^*-s_0J_1)-m_2(t)(\underline{A}^*-s_0A_1)-\alpha_2(t)((\underline{A}^*)^2-s_0A_1^2)\\[1mm]
&\geq&-m_2(t)(\underline{A}^*-s_0A_1)-\alpha_2(t)(\underline{A}^*+s_0A_1)(\underline{A}^*-s_0A_1)\\[1mm]
&=&-K_2^*(t)(\underline{A}^*-s_0A_1),
\end{array}$$
where $-K_2^*(t):=-m_2(t)-\alpha_2(t)(\underline{A}^*+s_0A_1)$.
Additionally, impulsive condition admits $(\underline{J}^*-s_0J_1)(0^+)=(\underline{J}^*-s_0J_1)(0)\geq0$, and
$$
(\underline{A}^*-s_0A_1)(0^+)=H(\underline{A}^*(0))-s_0H(A_1(0))\geq H(s_0A_1(0))-s_0H(A_1(0))\geq0
$$
according to the monotonicity of $H(A)$ and $H(A)/A$ in ($\mathcal{H}$1).

Regarding $(\underline{J}^*-s_0J_1,\underline{A}^*-s_0A_1)(0^+)$ as the new initial value and using the strong maximum principal in $t\in (0^+,\tau]$, yield
$$\underline{J}^*-s_0J_1>0\,\,{\rm{or}}\,\,\underline{J}^*-s_0J_1\equiv 0,\,t\in\{0^+\}\cup\{(0^+,\tau]\}$$
and
$$\underline{A}^*-s_0A_1>0\,\,{\rm{or}}\,\,\underline{A}^*-s_0A_1\equiv 0,\,t\in\{0^+\}\cup\{(0^+,\tau]\}.$$

Now we have the following four assertions.

$(\romannumeral1)$ Suppose $\underline{J}^*-s_0J_1>0$ and $\underline{A}^*-s_0A_1>0$ for $ t\in\{0^+\}\cup\{(0^+,\tau]\}$. There exist positive constants $\varepsilon_1$ and $\varepsilon_2$ such that $\underline{J}^*-s_0J_1\geq\varepsilon_1 J_1$ and $\underline{A}^*-s_0A_1\geq\varepsilon_2 A_1$. Taking $\varepsilon_0:=\min\{\varepsilon_1,\varepsilon_2\}$, then $\varepsilon_0+s_0\in S$, which contradicts to $s_0=\sup S$.

$(\romannumeral2)$ Assume that $\underline{J}^*-s_0J_1\equiv 0$ and $\underline{A}^*-s_0A_1\equiv 0$ for $t\in\{0^+\}\cup\{(0^+,\tau]\}$. It follows from the second equation of $A$ in problem \eqref{c11} that $s_0(s_0-1)\alpha_2(t)A_1^2=0$. However, the fact $0<s_0<1$ indicates $s_0(s_0-1)\alpha_2(t)A_1^2<0$, which leads a contradiction.

$(\romannumeral3)$ Assume that $\underline{J}^*-s_0J_1\equiv 0$ and $\underline{A}^*-s_0A_1>0$ for $t\in\{0^+\}\cup\{(0^+,\tau]\}$. Since $\underline{J}^*-s_0J_1\equiv 0$, it can be deduced from the first equation in \eqref{c11} that $b(t)(\underline{A}^*-s_0A_1)-\alpha_1(t)((\underline{J}^*)^2-s_0J_1^2)=0$. On the other hand, recalling $0<s_0<1$, we get $-\alpha_1(t)((\underline{J}^*)^2-s_0J_1^2)>-\alpha_1(t)(\underline{J}^*-s_0J_1)(\underline{J}^*+s_0J_1)=0$. Thus $b(t)(\underline{A}^*-s_0A_1)-\alpha_1(t)((\underline{J}^*)^2-s_0J_1^2)>0$ thanks to $\underline{A}^*-s_0A_1>0$, which leads a contradiction.

$(\romannumeral4)$ Assume that $\underline{J}^*-s_0J_1>0$ and $\underline{A}^*-s_0A_1\equiv0$ for $t\in\{0^+\}\cup\{(0^+,\tau]\}$. This is impossible. The concrete analysis is similar to that in $(\romannumeral3)$ and we omit it here.

As a result, $1\in S$, and problem \eqref{c11} admits a unique solution, which is denoted by $(J^*,A^*)(=(\underline J^*,\underline A^*)=(\bar J^*,\bar A^*))$.

\textbf{Step 2}\,\,\,The proof of $\lim\limits_{n\to+\infty}\inf{(J(t+n\tau,x),A(t+n\tau,x))}\geq(J^*(t),A^*(t))$ locally uniformly in $[0,\tau]\times(-\infty,+\infty)$.

Since $\lambda_1^D((-\infty,+\infty),H'(0))<0$, there exists a positive constant $L$ such that $$\lambda_1^D((-L,L),H'(0))<0.$$ Also, by $-g_\infty=h_\infty=\infty$, we can find the positive integer $n_L$ such that $h(t)-g(t)\geq 2L$ for any $t\geq n_LT$.

We first consider the following initial boundary problem
\begin{eqnarray}
\left\{
\begin{array}{lll}
\underline{J}_t-d_1\underline J_{xx}=b(t)\underline A-(a(t)+m_1(t))\underline J-\alpha_1(t)\underline J^2,  & t\in((n\tau)^+,(n+1)\tau],\,x \in (-L,L),\\[1mm]
\underline A_t-d_2\underline A_{xx}=a(t)\underline J-m_2(t)\underline A-\alpha_2(t)\underline A^2, & t\in((n\tau)^+,(n+1)\tau],\,x \in (-L,L),\\[1mm]
\underline J((n\tau)^+,x)=\underline J(n\tau,x),&x \in (-L,L),\,n\geq n_L,\\[1mm]
\underline A((n\tau)^+,x)=H(\underline A(n\tau,x)), &x \in (-L,L),\,n\geq n_L,\\[1mm]
\underline J(t,-L)=\underline J(t,L)=\underline A(t,-L)=\underline A(t,L)=0, & t\geq n_LT, \\[1mm]
\underline J(n_LT,x)=J(n_LT,x),\,\underline A(n_LT,x)=A(n_LT,x),& x\in[-L,L] \\[1mm]
\end{array} \right.
\label{c12}
\end{eqnarray}
in a fixed domain $[-L,L]$. We can apply the comparison principle to conclude that
$$(\underline J(t,x),\underline A(t,x))\leq(J(t,x),A(t,x)),\,\,\forall t\geq n_LT,\,\,-L\leq x\leq L.$$
Recall $\lambda_1^D((-L,L),H'(0))<0$, it follows from Lemma \ref{xuhaiyan} (IV) that
$$\lim\limits_{n\to+\infty}{(\underline J(t+n\tau,x),\underline A(t+n\tau,x))}=(U_L(t,x),V_L(t,x))\,\,\text{uniformly\,\,on\,\,} [0,\tau]\times[-L,L],$$
where $(U_L(t,x),V_L(t,x))$ is the periodic solution corresponding to initial problem \eqref{c12}, and satisfies
\begin{eqnarray*}
\left\{
\begin{array}{lll}
(U_L)_t-d_1(U_L){xx}=b(t)V_L-(a(t)+m_1(t))U_L-\alpha_1(t)U_L^2,  & t\in(0^+,\tau],\,x \in (-L,L),\\[1mm]
(V_L)_t-d_2(V_L)_{xx}=a(t)U_L-m_2(t)V_L-\alpha_2(t)V_L^2, & t\in(0^+,\tau],\,x \in (-L,L),\\[1mm]
U_L(0^+,x)=U_L(0,x),\,V_L(0^+,x)=H(V_L(0,x)), &x \in (-L,L),\\[1mm]
U_L(t,-L)=U_L(t,L)=V_L(t,-L)=V_L(t,L)=0, & t\in [0,\tau], \\[1mm]
U_L(0,x)=U_L(\tau,x),\,V_L(0,x)=V_L(\tau,x),& x\in[-L,L]. \\[1mm]
\end{array} \right.
\end{eqnarray*}
Subsequently, we get
$$\lim\limits_{n\to+\infty}\inf{(J(t+n\tau,x),A(t+n\tau,x))}\geq (U_L(t,x),V_L(t,x))\,\,\text{uniformly\,\,on\,\,}[0,\tau]\times[-L,L],$$
and letting $L\rightarrow\infty$ to conclude
$$
\lim\limits_{n\to+\infty}\inf{(J(t+n\tau,x),A(t+n\tau,x))}\geq (J^*(t),A^*(t))\,\,\text{locally\,\,uniformly\,\,in\,\, }[0,\tau]\times(-\infty,+\infty).
$$

\textbf{Step 3}\,\,\,The proof of $\lim\limits_{n\to+\infty}\sup{(J(t+n\tau,x),A(t+n\tau,x))}\leq(J^*(t),A^*(t))$ locally uniformly in $[0,\tau]\times(-\infty,+\infty)$.

Let $(\bar J(t),\bar A(t))$ be a solution to the following auxiliary problem
\begin{eqnarray}
\left\{
\begin{array}{lll}
\bar J_t=b(t)\bar A-(a(t)+m_1(t))\bar J-\alpha_1(t)\bar J^2,  & t\in((n\tau)^+,(n+)\tau],\\[1mm]
\bar A_t=a(t)\bar J-m_2(t)\bar A-\alpha_2(t)\bar A^2, & t\in((n\tau)^+,(n+)\tau],\\[1mm]
\bar J((n\tau)^+)=\bar J(n\tau),\,\bar A((n\tau)^+)=H(\bar A(n\tau)),\\[1mm]
\bar J(0)=\bar A(0)=C^*,\\[1mm]
\end{array} \right.
\label{c13}
\end{eqnarray}
where $C^*$ is an upper solution to problem \eqref{a01} defined in Lemma 2.2. A comparison argument yields $(J(t,x),A(t,x))\leq(\bar J(t),\bar A(t))$ for $t\in[0,+\infty)$ and $x\in(-\infty,+\infty)$, thus
\begin{equation}
\lim\limits_{n\to+\infty}\sup{(J(t+n\tau,x),A(t+n\tau,x))}\leq\lim\limits_{n\to+\infty}{(\bar J(t+n\tau),\bar A(t+n\tau))}
\label{c14}
\end{equation}
for $t\in[0,+\infty)$ and $x\in(-\infty,+\infty)$.

In the following, iteration sequences $\{\bar U^ {(m)}\}$ and $\{\bar V^ {(m)}\}$ satisfying
\begin{eqnarray*}
\left\{
\begin{array}{lll}
\bar U_t^{(m)}+K^*\bar U^{(m)}=K^* \bar U^{(m-1)}+b(t)\bar V^{(m-1)}-(a(t)+m_1(t))\bar U^{(m-1)}\\[1mm]
\qquad \qquad \quad \quad\quad-\alpha_1(t) (\bar U^{(m-1)})^2,  & t\in((n\tau)^+,(n+1)\tau],\\[1mm]
\bar V_t^{(m)}+K^*\bar V^{(m)}=K^* \bar V^{(m-1)}+a(t)\bar U^{(m-1)}-m_2(t)\bar V^{(m-1)}\\[1mm]
\qquad \qquad \quad \quad\quad-\alpha_2(t)(\bar V^{(m-1)})^2, & t\in((n\tau)^+,(n+1)\tau],\\[1mm]
\bar U^{(m)}((n\tau)^+)=\bar U^{(m-1)}((n+1)\tau),\\[1mm]
\bar V^{(m)}((n\tau)^+)=H(\bar V^{(m-1)}((n+1)\tau)),\\[1mm]
\bar U^{(m)}(n\tau)=\bar U^{(m-1)}((n+1)\tau),\,\bar V^{(m)}(n\tau)=\bar V^{(m-1)}((n+1)\tau)\\[1mm]
\end{array} \right.
\end{eqnarray*}
with initial value $\bar U^ {(0)}=\bar J(0)=C^*$ and  $\bar V^ {(0)}=\bar A(0)=C^*$ can be constructed.

According to \cite{P}, it is obvious that
$$(\{\bar U^ {(m+1)}\},\{\bar V^ {(m+1)}\})\leq (\{\bar U^ {(m)}\},\{\bar V^ {(m)}\})\leq\ldots\leq (\{\bar U^ {(0)}\},\{\bar V^ {(0)}\})=(C^*,C^*).$$
So the limit of iteration sequence exists and satisfies
\begin{equation}
\lim\limits_{m\to+\infty}{(\bar U^ {(m)}(t),\bar V^ {(m)}(t))}=(J^*(t), A^*(t)),
\label{c15}
\end{equation}
where $(J^*(t), A^*(t))$ is the unique solution defined in problem \eqref{c11}.

Recalling \eqref{c14} and \eqref{c15}, to prove $\lim\limits_{n\to+\infty}\sup{(J(t+n\tau,x),A(t+n\tau,x))}\leq(J^*(t),A^*(t))$ locally uniformly in $[0,\tau]\times(-\infty,+\infty)$, we only to verify
$$
(\bar J(t+n\tau),\bar A(t+n\tau))\leq(\bar U^ {(n)}(t),\bar V^ {(n)}(t))
$$
for any given integer $n$ and $t\in[0,\tau]$.
In fact, in view of the definition of $C^*$, one easily checks that
\begin{equation}
(\bar J(t),\bar A(t))\leq (C^*,C^*)=(\bar U^ {(0)}(t),\bar V^ {(0)}(t)),\,\forall t
\in [0,\tau],
\label{c16}
\end{equation}
which derives that
\begin{equation}
(\bar J(\tau),\bar A(\tau))\leq (\bar U^ {(0)}(\tau),\bar V^ {(0)}(\tau))=(\bar U^ {(1)}(0),\bar V^ {(1)}(0)).
\label{c17}
\end{equation}
Recalling the monotonicity of $H(A)$ to $A$ and \eqref{c16}, we acquire
\begin{equation}
(\bar J(\tau^+),\bar A(\tau^+))=(\bar J(\tau), H(\bar A(\tau)))< (\bar U^ {(0)}(\tau),H(\bar V^ {(0)}(\tau)))= (\bar U^{(1)}(0^+),\bar V^{(1)}(0^+)).
\label{c18}
\end{equation}

It follows from \eqref{c17}, \eqref{c18} and the comparison argument that
\begin{equation}
(\bar J(t+\tau),\bar A(t+\tau))\leq (\bar U^ {(1)}(t),\bar V^ {(1)}(t))\,\,\text{for}\,\,t\in[0,\tau].
\label{c19}
\end{equation}

Noticing inequalities \eqref{c16} with $n=0$ and \eqref{c19} with $n=1$, and using the induction for $n$, give that
$$
(\bar J(t+n\tau),\bar A(t+n\tau))\leq(\bar U^ {(n)}(t),\bar V^ {(n)}(t))\,\,\text{for}\,\,t\in[0,\tau].
$$
Therefore, we obtain $$(J^*(t),A^*(t))\leq\lim\limits_{n\to+\infty}\inf{(J,A)(t+n\tau,x)}\leq\lim\limits_{n\to+\infty}\sup{(J,A)(t+n\tau,x)}
\leq(J^*(t),A^*(t))$$
locally uniformly in $[0,\tau]\times(-\infty,+\infty)$ and \eqref{c00} hold. Steps 1-3 finish the proof.
\epf

\begin{lem}
Suppose that $\lambda_1^D((-\infty,+\infty),H'(0))<0$. If $h_\infty<\infty$ or $-\infty<g_\infty$, then $-\infty<g_\infty<h_\infty<+\infty$.
\end{lem}
\bpf Without loss of generality, on the contrary we assume that $h_\infty<\infty$ and $g_\infty=-\infty$. A large positive integer $n_0$ can be
found such that $\lambda_0^F(n_0\tau)<0$, and $\lambda_0^F(t)<\lambda_0^F(n_0\tau)<0$ for $t>n_0\tau$ considering the decreasing of  $\lambda_0^F(t)$ in $t$. We denote $\tau_0:=n_0\tau$ for simplicity unless otherwise specified.

Some preliminary work about the eigenvalue problem is firstly given. We consider the following $\tau-$periodic eigenvalue problem with $\varepsilon$
\begin{eqnarray}
\left\{
\begin{array}{lll}
\phi_t-d_1\phi_{xx}-\varepsilon\phi_x=b(t)\psi-[a(t)+m_1(t)]\phi+\lambda_1^D\phi,& t\in(0^+,\tau],\,\,x \in (g(\tau_0),h(\tau_0)),  \\[1mm]
\psi_t-d_2 \psi_{xx}-\varepsilon\frac{d_2}{d_1}\psi_x=a(t)\phi-m_2(t)\psi+\lambda_1^D\psi, & t\in(0^+,\tau],\,\,x \in (g(\tau_0),h(\tau_0)),  \\[1mm]
\phi(0^+,x)=\phi(0,x),& x\in (g(\tau_0),h(\tau_0)),\\[1mm]
\psi(0^+,x)=H'(0)\psi(0,x),& x\in (g(\tau_0),h(\tau_0)),\\[1mm]
\phi(0,x)=\phi(\tau,x),\psi(0,x)=\psi(\tau,x), & x\in [g(\tau_0),h(\tau_0)],\\[1mm]
\phi(t,x)=\psi(t,x)=0, & t\in[0,\tau], \,\,x=\{g(\tau_0),h(\tau_0)\}
\end{array} \right.
\label{c20}
\end{eqnarray}
in the fixed domain.

If $\varepsilon=0$, problem \eqref{c20} is the eigenvalue problem of  \eqref{c01} with $(l_1,l_2):=(g(\tau_0),h(\tau_0))$, and we denote eigenfunction related to $\lambda_0^F(\tau_0)(:=\lambda_1^D((g(\tau_0),h(\tau_0)),H'(0))$ by $(\phi,\psi)$, which is defined in \eqref{c02}. On the other hand, we choose some $\varepsilon$ sufficiently small, such that the principal eigenvalue $\lambda_0^{F,\varepsilon}(\tau_0)$ of problem \eqref{c20} is negative. Then the corresponding eigenfunction $(\phi^\varepsilon,\psi^\varepsilon)$ can be expressed by
$$
\phi^\varepsilon(t,x)=\alpha(t)e^{-\frac{\varepsilon}{2d_1}x}\cos[\frac{\pi}{h(\tau_0)-g(\tau_0)}(x-\frac{h(\tau_0)+g(\tau_0)}{2})]
$$
and
$$
\psi^\varepsilon(t,x)=\beta(t)e^{-\frac{\varepsilon}{2d_1}x}\cos[\frac{\pi}{h(\tau_0)-g(\tau_0)}(x-\frac{h(\tau_0)+g(\tau_0)}{2})].
$$
Direct calculations yield
\begin{eqnarray*}
\left\{
\begin{array}{lll}
\phi^\varepsilon_{x}=-\frac{\varepsilon}{2d_1}\phi^\varepsilon-\frac{\pi}{h(\tau_0)-g(\tau_0)}e^{-\frac{\varepsilon}{2d_1}x}
\sin(\frac{\pi}{h(\tau_0)-g(\tau_0)}(x-\frac{h(\tau_0)+g(\tau_0)}{2}))\alpha(t),\\[2mm]
\psi^\varepsilon_{x}=-\frac{\varepsilon}{2d_1}\psi^\varepsilon-\frac{\pi}{h(\tau_0)-g(\tau_0)}e^{-\frac{\varepsilon}{2d_1}x}
\sin(\frac{\pi}{h(\tau_0)-g(\tau_0)}(x-\frac{h(\tau_0)+g(\tau_0)}{2}))\beta(t)
\end{array} \right.
\end{eqnarray*}
and
\begin{eqnarray*}
\left\{
\begin{array}{lll}
\phi^\varepsilon_{xx}=[-\frac{\varepsilon^2}{4d_1^2}-(\frac{\pi}{h(\tau_0)-g(\tau_0)})^2]\phi^\varepsilon
-\frac{\varepsilon}{d_1}\phi_x^\varepsilon,\\[2mm]
\psi^\varepsilon_{xx}=[-\frac{\varepsilon^2}{4d_1^2}-(\frac{\pi}{h(\tau_0)-g(\tau_0)})^2]\psi^\varepsilon-\frac{\varepsilon}{d_1}\psi_x^\varepsilon,
\end{array} \right.
\end{eqnarray*}
which will be used in the sequel.

Moreover, it is clear that there exists $x^*\in(g(\tau_0),h(\tau_0))$ such that
$$
\phi^\varepsilon_x(t,x)>0,\psi^\varepsilon_x(t,x)>0\,\,\text{for}\,\,x\in[g(\tau_0),x^*),
$$
$$
\phi^\varepsilon_x(t,x)<0,\psi^\varepsilon_x(t,x)<0\,\,\text{for}\,\,x\in(x^*,h(\tau_0)],
$$
and $\phi^\varepsilon_x(t,x^*)=\psi^\varepsilon_x(t,x^*)=0$. The rest of proof covers two steps.

\textbf{Step 1}\,\,\,We firstly fix $\varepsilon=0$ so as to give a lower bound of $(J(t,x), A(t,x))$ for $t\geq \tau_0$ and $x\in [g(\tau_0),h(\tau_0)]$.

Define
\begin{equation*}
\underline{J}(t,x)=\varepsilon_1\frac{\rho_1}{H'(0)}\phi(t,x),\,\,t\geq0,\,\,x\in [g(\tau_0),h(\tau_0)]
\end{equation*}
and
\begin{eqnarray*}
\underline{A}(t,x)=
\begin{cases}
\displaystyle \varepsilon_1{\psi}(n\tau,x), &t=n\tau, \\[1mm]
\displaystyle \varepsilon_1{\frac{\rho_1}{H'(0)}}\psi((n\tau)^+,x), &t=(n\tau)^+,\\[1mm]
\displaystyle \varepsilon_1{\frac{\rho_1}{H'(0)}}e^{(-\lambda_0^F(\tau_0)-\delta)(t-n\tau)}\psi(t,x),&t\in((n\tau)^+,(n+1)\tau]
\end{cases}
\end{eqnarray*}
for $x\in [g(\tau_0),h(\tau_0)]$. In fact, we fix $\rho_1:=H'(0)e^{(\lambda_0^F(\tau_0)+\delta)\tau}$ such that $\underline{A}(n\tau,x)=\underline{A}((n+1)\tau,x)$. Also, a positive constant $\delta$
can be chosen later to guarantee $\lambda_0^F(\tau_0)+\delta<0$ and $\rho_1<H'(0)$.

Hence, we obtain for $t\in((n\tau)^+,(n+1)\tau]$ and $x\in (g(\tau_0),h(\tau_0))$ that
$$\begin{array}{llllll}
&&\underline{J}_t-d_1\underline{J}_{xx}-b(t)\underline{A}+(a(t)+m_1(t))\underline{J}+\alpha_1(t)
\underline{J}^2\\[1mm]
&=&\varepsilon_1\frac{\rho_1}{H'(0)}\{\phi_t-d_1\phi_{xx}-b(t)e^{(-\lambda_0^F(\tau_0)-\delta)(t-n\tau)}\psi
+(a(t)+m_1(t))\phi+\alpha_1(t)\varepsilon_1\frac{\rho_1}{H'(0)}\phi^2\}\\[1mm]
&=&\varepsilon_1\frac{\rho_1}{H'(0)}[\lambda_0^F(\tau_0)\phi+b(t)\psi(1-e^{(-\lambda_0^F(\tau_0)-\delta)(t-n\tau)})+\alpha_1(t)
\varepsilon_1\frac{\rho_1}{H'(0)}\phi^2]\\[1mm]
&\leq&\varepsilon_1\frac{\rho_1}{H'(0)}\phi[\lambda_0^F(\tau_0)+\alpha^M_1\varepsilon_1\frac{\rho_1}{H'(0)}\phi]\\[1mm]
&\leq&0
\end{array}$$
and
$$\begin{array}{llllll}
\,\,\,\,\,\,\,\underline{A}_t-d_2\underline{A}_{xx}-a(t)\underline{J}+m_2(t)\underline{A}+\alpha_2(t)\underline{A}^2\\[1mm]
=\,\varepsilon_1\frac{\rho_1}{H'(0)}e^{(-\lambda_0^F(\tau_0)-\delta)(t-n\tau)}[-\delta\psi+a(t)\phi(1-e^{(\lambda_0^F(\tau_0)+\delta)(t-n\tau)})
+\alpha_2(t)\varepsilon_1\frac{\rho_1}{H'(0)}e^{(-\lambda_0^F(\tau_0)-\delta)(t-n\tau)}\psi^2]\\[1mm]
\leq\,\varepsilon_1\frac{\rho_1}{H'(0)}e^{(-\lambda_0^F(\tau_0)-\delta)(t-n\tau)}[-\delta\psi+a(t)\phi(1-e^{(\lambda_0^F(\tau_0)+\delta)\tau})
+\alpha_2(t)\varepsilon_1\frac{\rho_1}{H'(0)}e^{(-\lambda_0^F(\tau_0)-\delta)(t-n\tau)}\psi^2]\\[1mm]
\leq\,\varepsilon_1\frac{\rho_1}{H'(0)}e^{(-\lambda_0^F(\tau_0)-\delta)(t-n\tau)}\{\psi[-\delta+a^ML(1-e^{(\lambda_0^F(\tau_0)+\delta)\tau})]
+\alpha^M_2\varepsilon_1\frac{\rho_1}{H'(0)}e^{(-\lambda_0^F(\tau_0)-\delta)\tau}\psi^2\}\\[1mm]
\leq\,0
\end{array}$$
on account of $1/L\leq\frac{\phi}{\psi}\leq L$ for some positive constant $L$ and $-\delta+a^ML(1-e^{(\lambda_0^F(\tau_0)+\delta)\tau})<0$.

We claim that the positive constant $\delta$ satisfying
\begin{eqnarray*}
\left\{
\begin{array}{lll}
-\delta+a^ML(1-e^{(\lambda_0^F(\tau_0)+\delta)\tau})<0,\\[2mm]
\lambda_0^F(\tau_0)+\delta<0
\end{array} \right.
\end{eqnarray*}
exists owing to mathematical analysis by the fact $\delta:=-\lambda_0^F(\tau_0)-\varepsilon^*$ with sufficiently small $\varepsilon^*$, see also in \cite{XLS}(Theorem 4.3) for more details.

Furthermore, for $x\in [g(\tau_0),h(\tau_0)]$, we get
$$
\underline{J}((n\tau)^+,x)=\varepsilon_1{\frac{\rho_1}{H'(0)}}\phi((n\tau)^+,x)=\varepsilon_1{\frac{\rho_1}{H'(0)}}\phi(n\tau,x)=\underline{J}(n\tau,x),
$$
and ($\mathcal{H}$2) implies
$$\begin{array}{llllll}
&&\underline{A}((n\tau)^+,x)-H(\underline{A}(n\tau,x))\\[1mm]
&=&\varepsilon_1\frac{\rho_1}{H'(0)}\psi((n\tau)^+,x)-H(\varepsilon_1\psi(n\tau,x))\\[1mm]
&\leq&\varepsilon_1\rho_1\psi(n\tau,x)-H'(0)\varepsilon_1\psi(n\tau,x)+D[\varepsilon\psi(n\tau,x)]^\gamma\\[1mm]
&=&\varepsilon_1 \psi(n\tau,x)\{\rho_1-H'(0)+D[\varepsilon_1\psi(n\tau,x)]^{\gamma-1}\}\\[1mm]
&\leq&0.
\end{array}$$
Also, $\underline{J}(\tau_0,x)\leq J(\tau_0,x)$ and $\underline{A}(\tau_0,x)\leq A(\tau_0,x)$ hold for $x\in [g(\tau_0),h(\tau_0)]$ with sufficiently small $\varepsilon_1$. A simple comparison principle yields $(\underline{J}(t,x),\underline{A}(t,x))\leq (J(t,x),A(t,x))$  for $t\geq \tau_0$ and $x\in [g(\tau_0),h(\tau_0)]$.
Whence, we obtain
\begin{equation}
\underline{J}(t,x^*)\leq J(t,x^*)\,\,\text{and}\,\,\underline{A}(t,x^*)\leq A(t,x^*)\,\,\text{for}\,\,t\geq \tau_0.
\label{c211}
\end{equation}

\textbf{Step 2}\,\,\,We secondly aim to construct a lower solution to problem \eqref{a01}.

Since $\lim\limits_{t\to+\infty}{h'(t)}=0$, there exists some positive integer $n_1\geq n_0(\tau_1:=n_1\tau\geq n_0\tau)$, such that for $t> n_1\tau$ and $t\neq n\tau$,
$$
h'(t)<\min\{\varepsilon\frac{h(\tau_0)-x^*}{h_\infty-x^*},\varepsilon \frac{d_2}{d_1}\frac{h(\tau_0)-x^*}{h_\infty-x^*}\}
$$
hold.

In the following, we will find the lower bound $(\hat{J},\hat A)(t,x)$ to problem \eqref{a01} for $t\geq \tau_1$ and $x\in[x^*,h(t)]$. Define
\begin{equation*}
\hat{J}(t,x)=\varepsilon_2\frac{\rho_2}{H'(0)}\phi^\varepsilon(t,\hat x),\,\,t\geq0,\,\,x\in [x^*,h(t)]
\end{equation*}
and
\begin{eqnarray*}
\hat {A}(t,x)=
\begin{cases}
\displaystyle \varepsilon_2{\psi^\varepsilon}(n\tau,\hat x), &t=n\tau, \\[3mm]
\displaystyle \varepsilon_2{\frac{\rho_2}{H'(0)}}\psi^\varepsilon((n\tau)^+,\hat x), &t=(n\tau)^+,\\[3mm]
\displaystyle \varepsilon_2{\frac{\rho_2}{H'(0)}}e^{(-\lambda_0^{F,\varepsilon}(\tau_0)-\delta)(t-n\tau)}\psi^\varepsilon(t,\hat x),&t\in((n\tau)^+,(n+1)\tau],
\end{cases}
\end{eqnarray*}
where $(\phi^\varepsilon,\psi^\varepsilon)$ satisfies \eqref{c20} and $\hat x:=x^*+\frac{h(\tau_0)-x^*}{h(t)-x^*}(x-x^*)$ for $x\in [x^*,h(t)]$, so $\hat x\in [x^*,h(\tau_0)]$. Arguing as in the above,
$\rho_2:=H'(0)e^{(\lambda_0^{F,\varepsilon}(\tau_0)+\delta)\tau}<H'(0)$ is obtained.

We shall check that
\begin{eqnarray}
\left\{
\begin{array}{lll}
\hat{J}_t-d_1\hat{J}_{xx}-b(t)\hat{A}+(a(t)+m_1(t))\hat{J}+\alpha_1(t)\hat{J}^2\leq0, &(t,x)\in \Omega_\tau,\\[2mm]
\hat{A}_t-d_2\hat{A}_{xx}-a(t)\hat{J}+m_2(t)\hat{A}+\alpha_2(t)\hat{A}^2\leq0, &(t,x)\in \Omega_\tau,\\[2mm]
\hat{J}((n\tau)^+,x)=\hat J(n\tau,x)),\,\,\hat{A}((n\tau)^+,x)\leq H(\hat{A}(n\tau,x)),&x\in [x^*,h(n\tau)],\\[2mm]
\hat J(t,x^*)\leq J(t,x^*),\,\,\hat A(t,x^*)\leq A(t,x^*),&t\geq \tau_1,\\[2mm]
\hat J(t,h(t))=\hat A(t,h(t))=0,&t\geq \tau_1,\\[2mm]
\hat J(\tau_1,x)\leq J(\tau_1,x),\,\,\hat A(\tau_1,x)\leq A(\tau_1,x),&x\in [x^*,h(\tau_1)],
\end{array} \right.
\label{c21}
\end{eqnarray}
where $\Omega_{\tau}:=\{(t,x):(n\tau)^+<t\leq(n+1)\tau],\,x^*<x<h(t)\}$. Once \eqref{c21} holds, the comparison principle yields $(\hat J,\hat A)\leq(J,A)$ for $t\geq \tau_1$ and $x\in [x^*,h(t)]$, which together with $J(t,h(t))=\hat J(t,h(t))=A(t,h(t))=\hat A(t,h(t))=0$, yields
$$
J_x(t,h(t))\leq\hat J_x(t,h(t))<0,\,\, A_x(t,h(t))\leq\hat A_x(t,h(t))<0\,\,\text{for}\,\,t\geq\tau_1.
$$
On the other hand, using the equation $$J_x(t,h(t))=-h'(t)/\mu_1-\mu_2/\mu_1 A_x(t,h(t))$$
yields $\lim\limits_{t\to\infty}{J_x(t,h(t))}>0$, which leads a contradiction, and then $-\infty<g_\infty<h_\infty<\infty$.

Now we verify \eqref{c21}. The impulsive conditions $\hat{J}((n\tau)^+,x)=\hat J(n\tau,x))$ and $\hat{A}((n\tau)^+,x)\leq H(\hat{A}(n\tau,x))$ can be derived by the same manner in Step 1. For $t\geq \tau_1$, the right boundary satisfies $\hat J(t,h(t))=\hat A(t,h(t))=0$, the left boundary meets $\hat J(t,x^*)\leq \underline J(t,x^*)\leq J(t,x^*)$ and $\hat A(t,x^*)\leq \underline A(t,x^*)\leq A(t,x^*)$ considering the sufficiently small $\varepsilon_2<\varepsilon_1$ and equation \eqref{c211}. Moreover, $\hat J(\tau_1,x)\leq J(\tau_1,x)$ and $\hat A(\tau_1,x)\leq A(\tau_1,x)$ hold for $x\in [x^*,h(\tau_1)]$ and sufficiently small $\varepsilon_2$.

It remains to prove the first two inequalities in \eqref{c21}. In fact,
$$\begin{array}{lllll}
\,\,\,\,\,\,\,\hat{J}_t-d_1\hat{J}_{xx}-b(t)\hat{A}+[a(t)+m_1(t)]\hat{J}+\alpha_1(t)\hat{J}^2\\[1mm]
=\,\varepsilon_2\frac{\rho_2}{H'(0)}\{\phi_t^\varepsilon-\frac{(h(\tau_0)-x^*)(x-x^*)h'(t)}{(h(t)-x^*)^2}\phi_x^\varepsilon-d_1
(\frac{h(\tau_0)-x^*}{h(t)-x^*})^2\phi_{xx}^\varepsilon+[a(t)+m_1(t)]\phi^\varepsilon\\[1mm]
\,\,\,\,\,\,+\alpha_1(t)\varepsilon_2\frac{\rho_2}{H'(0)}(\phi^\varepsilon)^2    -b(t)e^{(-\lambda_0^{F,\varepsilon}(\tau_0)-\delta)(t-n\tau)}\psi^\varepsilon\}\\[1mm]
=\,\varepsilon_2\frac{\rho_2}{H'(0)}\{\lambda_0^{F,\varepsilon}(\tau_0)\phi^\varepsilon+d_1
[1-(\frac{h(\tau_0)-x^*}{h(t)-x^*})^2]\phi_{xx}^\varepsilon+\varepsilon\phi_x^\varepsilon- \frac{(h(\tau_0)-x^*)(x-x^*)h'(t)}{(h(t)-x^*)^2}\phi_x^\varepsilon\\[1mm]
\,\,\,\,\,\,+b(t)(1-e^{(-\lambda_0^{F,\varepsilon}(\tau_0)-\delta)(t-n\tau)})\psi^\varepsilon+   \alpha_1(t)\varepsilon_2\frac{\rho_2}{H'(0)}(\phi^\varepsilon)^2\}\\[1mm]
<\,\varepsilon_2\frac{\rho_2}{H'(0)}\{\lambda_0^{F,\varepsilon}(\tau_0)\phi^\varepsilon-d_1
[1-(\frac{h(\tau_0)-x^*}{h(t)-x^*})^2][\frac{\varepsilon^2}{4d_1^2}+(\frac{\pi}{h(\tau_0)-g(\tau_0)})^2]\phi^\varepsilon\\[1mm]
\,\,\,\,\,\,+[\varepsilon(\frac{h(\tau_0)-x^*}{h(t)-x^*})^2-\frac{(h(\tau_0)-x^*)(x-x^*)h'(t)}{(h(t)-x^*)^2}]\phi_x^\varepsilon+ \alpha^M_1\varepsilon_2\frac{\rho_1}{H'(0)}(\phi^\varepsilon)^2 \}\\[1mm]
<\,\varepsilon_2\frac{\rho_2}{H'(0)}\{\lambda_0^{F,\varepsilon}(\tau_0)\phi^\varepsilon
+[\varepsilon(\frac{h(\tau_0)-x^*}{h(t)-x^*})^2-\frac{(h(\tau_0)-x^*)(x-x^*)h'(t)}{(h(t)-x^*)^2}]\phi_x^\varepsilon+ \alpha^M_1\varepsilon_2\frac{\rho_1}{g'(0)}(\phi^\varepsilon)^2 \}
\end{array}$$
in view of $\frac{h(\tau_0)-x^*}{h(t)-x^*}\in(0,1)$ and $-d_1
[1-(\frac{h(\tau_0)-x^*}{h(t)-x^*})^2][\frac{\varepsilon^2}{4d_1^2}+(\frac{\pi}{h(\tau_0)-g(\tau_0)})^2]\phi^\varepsilon<0$ for $t\geq \tau_1(>\tau_0)$ and $x\in[x^*,h(t)]$. Also,
$$\begin{array}{llllll}
\,\,\,\,\,\,\,\hat{A}_t-d_2\hat{A}_{xx}-a(t)\hat{J}+m_2(t)\hat{A}+\alpha_2(t)\hat{A}^2\\[1mm]
=\,\varepsilon_2\frac{\rho_1}{H'(0)}e^{(-\lambda_0^{F,\varepsilon}-\delta)(t-n\tau)}\{(-\lambda_0^{F,\varepsilon}-\delta)\psi^\varepsilon-
\frac{(h(\tau_0)-x^*)(x-x^*)h'(t)}{(h(t)-x^*)^2}\psi_x^\varepsilon+\psi_t^\varepsilon\\[1mm]
\,\,\,\,\,\,-d_2(\frac{h(\tau_0)-x^*}{h(t)-x^*})^2\psi_{xx}^\varepsilon
-a(t)e^{(\lambda_0^{F,\varepsilon}+\delta)(t-n\tau)}\phi^\varepsilon+m_2(t)\psi^\varepsilon
+\alpha_2(t)\varepsilon_2\frac{\rho_1}{H'(0)}e^{(-\lambda_0^{F,\varepsilon}-\delta)(t-n\tau)}(\psi^\varepsilon)^2\}\\[1mm]
\leq\, \varepsilon_2\frac{\rho_1}{H'(0)}e^{(-\lambda_0^F-\delta)(t-n\tau)}\{-\delta\psi^\varepsilon
+d_2[1-(\frac{h(\tau_0)-x^*}{h(t)-x^*})^2]\psi_{xx}^\varepsilon+[\frac{d_2}{d_1}\varepsilon-\frac{(h(\tau_0)-x^*)(x-x^*)h'(t)}{(h(t)-x^*)^2}]
\psi_x^\varepsilon\\[1mm]
\,\,\,\,\,\,+(1-e^{(\lambda_0^{F,\varepsilon}+\delta)\tau})a(t)\phi^\varepsilon
+\alpha_2(t)\varepsilon_2\frac{\rho_1}{H'(0)}e^{(-\lambda_0^{F,\varepsilon}(\tau_0)-\delta)(t-n\tau)}(\psi^\varepsilon)^2\}\\[1mm]
\leq\, \varepsilon_2\frac{\rho_1}{H'(0)}e^{(-\lambda_0^F-\delta)(t-n\tau)}\{[-\delta-d_2(1-(\frac{h(\tau_0)-x^*}{h(t)-x^*})^2)
(\frac{\varepsilon^2}{4d_1^2}+(\frac{\pi}{h(\tau_0)-g(\tau_0)})^2)]\psi^\varepsilon+[\frac{d_2}{d_1}\varepsilon(\frac{h(\tau_0)-x^*}{h(t)-x^*})^2\\[1mm]
\,\,\,\,\,\,-\frac{(h(\tau_0)-x^*)(x-x^*)h'(t)}{(h(t)-x^*)^2}]\psi_x^\varepsilon
+(1-e^{(\lambda_0^{F,\varepsilon}+\delta)\tau}) a^M\phi^\varepsilon
+\alpha^M_2\varepsilon_2\frac{\rho_1}{H'(0)}e^{(-\lambda_0^{F,\varepsilon}-\delta)(t-n\tau)}(\psi^\varepsilon)^2\}\\[1mm]
\leq\,\varepsilon_2\frac{\rho_1}{H'(0)}e^{(-\lambda_0^{F,\varepsilon}-\delta)(t-n\tau)}\{[-\delta-d_2(1-(\frac{h(\tau_0)-x^*}{h(t)-x^*})^2)
(\frac{\varepsilon^2}{4d_1^2}+(\frac{\pi}{h(\tau_0)-g(\tau_0)})^2)\\[1mm]
\,\,\,\,\,\,+(1-e^{(\lambda_0^{F,\varepsilon}+\delta)\tau})\tilde L a^M]\psi^\varepsilon+[\frac{d_2}{d_1}\varepsilon(\frac{h(\tau_0)-x^*}{h(t)-x^*})^2-
\frac{(h(\tau_0)-x^*)(x-x^*)h'(t)}{(h(t)-x^*)^2}]\psi_x^\varepsilon\\[1mm]
\,\,\,\,\,\,+\alpha^M_2\varepsilon_2\frac{\rho_1}{H'(0)}
e^{(-\lambda_0^{F,\varepsilon}-\delta)(t-n\tau)}(\psi^\varepsilon)^2\}\\[1mm]
\leq\,\varepsilon_2\frac{\rho_1}{H'(0)}e^{(-\lambda_0^{F,\varepsilon}-\delta)(t-n\tau)}\{[-\delta+(1-e^{(\lambda_0^{F,\varepsilon}+\delta)\tau})
\tilde{L}  a^M]\psi^\varepsilon+[\frac{d_2}{d_1}\varepsilon(\frac{h(\tau_0)-x^*}{h(t)-x^*})^2\\[1mm]
\,\,\,\,\,\,-\frac{(h(\tau_0)-x^*)(x-x^*)h'(t)}{(h(t)-x^*)^2}]\psi_x^\varepsilon+\alpha^M_2\varepsilon_2\frac{\rho_1}{H'(0)}
e^{(-\lambda_0^{F,\varepsilon}-\delta)(t-n\tau)}(\psi^\varepsilon)^2\}\\[1mm]
\end{array}$$
in view of $\frac{\phi^\varepsilon}{\psi\varepsilon}\leq \tilde{L}$ for a positive constant $\tilde{L}$, which is independent of $\varepsilon$, and $-d_2(1-(\frac{h(\tau_0)-x^*}{h(t)-x^*})^2)
(\frac{\varepsilon^2}{4d_1^2}+(\frac{\pi}{h(\tau_0)-g(\tau_0)})^2)<0$ for $t\geq \tau_1$ and $x\in[x^*,h(t)]$.

Moreover, recalling that $h'(t)<\min\{\varepsilon\frac{h(\tau_0)-x^*}{h_\infty-x^*},\varepsilon \frac{d_2}{d_1}\frac{h(\tau_0)-x^*}{h_\infty-x^*}\}$, we obtain $$\varepsilon(\frac{h(\tau_0)-x^*}{h(t)-x^*})^2-\frac{(h(\tau_0)-x^*)(x-x^*)h'(t)}{(h(t)-x^*)^2}>0$$ and $$\varepsilon\frac{d_2}{d_1}(\frac{h(\tau_0)-x^*}{h(t)-x^*})^2-\frac{(h(\tau_0)-x^*)(x-x^*)h'(t)}{(h(t)-x^*)^2}>0.$$ Also, it is clear that  $(\phi_x^\varepsilon(t,x),\psi_x^\varepsilon(t,x))\leq (0,0)$ for $t\geq \tau_1$ and $x^*\leq x \leq h(t)$. Therefore,
$$\begin{array}{lllll}
&&\hat{J}_t-d_1\hat{J}_{xx}-b(t)\hat{A}+[a(t)+m_1(t)]\hat{J}+\alpha_1(t)\hat{J}^2\\[1mm]
&\leq&\varepsilon_2\frac{\rho_2}{H'(0)}\{\lambda_0^{F,\varepsilon}(\tau_0)\phi^\varepsilon
+\varepsilon_2\alpha^M_1\frac{\rho_1}{g'(0)}(\phi^\varepsilon)^2 \}\\[1mm]
&\leq0&
\end{array}$$
owning to $\lambda_0^{F,\varepsilon}(\tau_0)<0$ and the arbitrary of $\varepsilon_2$. Similarly,
$$\begin{array}{llllll}
&&\hat{A}_t-d_2\hat{A}_{xx}-a(t)\hat{J}+m_2(t)\hat{A}+\alpha_2(t)\hat{A}^2\\[1mm]
&\leq&\varepsilon_2\frac{\rho_1}{H'(0)}e^{(-\lambda_0^{F,\varepsilon}-\delta)(t-n\tau)}\{[-\delta+(1-e^{(\lambda_0^{F,\varepsilon}+\delta)\tau})
\tilde{L} a^M]\psi^\varepsilon\\[1mm]
&+&\varepsilon_2\alpha^M_2\frac{\rho_1}{H'(0)}
e^{(-\lambda_0^{F,\varepsilon}-\delta)(t-n\tau)}(\psi^\varepsilon)^2\}\\[1mm]
&\leq0&
\end{array}$$
can be derived in view of the existence of  positive constant $\delta$ that satisfying
\begin{eqnarray*}
\left\{
\begin{array}{lll}
-\delta+(1-e^{(\lambda_0^{F,\varepsilon}+\delta)\tau})\tilde{L} a^M)<0,\\[1mm]
\lambda_0^{F,\varepsilon}+\delta<0.
\end{array} \right.
\end{eqnarray*}

The proof is now completed.
\epf

\vspace{3mm}
In conclusion, The following statements hold thanks to Theorems 3.3, 3.5 and Lemma 3.6.

\begin{thm}(Spreading-vanishing dichotomy) Let $(J,A;(g,h))$ be the solution to problem \eqref{a01}, the following statements hold:\\
$(\romannumeral 1)$ if $\lambda_1^D((-\infty,+\infty),H'(0))\geq0$, then $-\infty<g_\infty<h_\infty<+\infty$ and

$\lim\limits_{t\to+\infty}{||(J(t,\cdot),A(t,\cdot))||_{C[g(t),h(t)]}}=(0,0)$;\\
$(\romannumeral 2)$ if $\lambda_1^D((-\infty,+\infty),H'(0))<0$, and

$(a)$ $-g_\infty=h_\infty=\infty$, then $\lim\limits_{m\to+\infty}{(J(t+mT,x),A(t+mT,x))}=(J^*(t),A^*(t))$ locally

uniformly in $[0,\tau]\times(-\infty,+\infty)$; or

$(b)$ $-\infty<g_\infty<h_\infty<+\infty$, then $\lambda_1^D((g_\infty,h_\infty),H'(0))\geq0$ and

$\lim\limits_{t\to+\infty}{||(J(t,\cdot),A(t,\cdot))||_{C[g(t),h(t)]}}=(0,0).$
\end{thm}

\section{Sufficient conditions for spreading or vanishing}
Recalling hypothesises about harvesting pulse introduced in ($\mathcal{H}$1) and results verified in Section 3, we will investigate the sufficient conditions about pulse and expanding capacities for governing spatially spreading or vanishing of juveniles and adults in the following.

For the initial habitat $(-h_0,h_0)$ and the whole area $(-\infty,+\infty)$ of species, let us first denote $g^*:=H_1'(0)$ and $g_{*}:=H_2'(0)$, which satisfy
$$
\lambda_1^D((-h_0,h_0),g^{*})=\lambda_1^D((-h_0,h_0),H_1'(0))=0$$
and
$$\lambda_1^D((-\infty,+\infty),g_{*})=\lambda_1^D((-\infty,+\infty),H_2'(0))=0,$$
respectively. We also claim that $g_*<g^*$. Otherwise, we suppose $g_*\geq g^*$. $\lambda_1^D((-h_0,h_0),g_*)\leq\lambda_1^D((-h_0,h_0),g^{*})=0$ hold according to the strictly decreasing of $\lambda_1^D((-h_0,h_0),H'(0))$ in $H'(0)$. On the other hand, it follows from $\lambda_1^D((-\infty,+\infty),g_*)=0$ and the strictly decreasing of $\lambda_1^D(\Omega,g_*)$ in $\Omega$ that $\lambda_1^D((-h_0,h_0),g_{*})>0$, which leads a contradiction, so $g_*<g^*$. It is worth mentioning that $H_1(u)$ and $H_2(u)$ satisfying ($\mathcal{H}$1) are two impulsive functions with $u\geq0$.

\vspace{2mm}
Let us first generalize Theorem 3.3 for the condition of $H'(0)\leq g_{*}$.

\begin{cor}
If $H'(0)\leq g_{*}$, then $-\infty<g_\infty<h_\infty<+\infty$ and species vanish.
\end{cor}
\bpf
Since $H'(0)\leq g_{*}$ and Lemma 3.1(ii) holds, we obtain $\lambda_1^D((-\infty,+\infty),H'(0))\geq \lambda_1^D((-\infty,+\infty),g_*)$. Therefore, $\lambda_1^D((-\infty,+\infty),H'(0))\geq 0$, which together with Theorem 3.3 to conclude $h_\infty-g_\infty<+\infty$ and species finally vanish.
\epf

\begin{thm}\label{aas1}
If $H'(0)\geq g^{*}$, then $-g_\infty=h_\infty=\infty$, and species spread.
\end{thm}
\bpf
Since $H'(0)\geq g^{*}$, it follows from the monotonic decreasing of $\lambda_1^D((-h_0,h_0),H'(0))$ with respect to $H'(0)$ in Lemma 3.1(ii) that $\lambda_1^D((-h_0,h_0),H'(0))\leq \lambda_1^D((-h_0,h_0),g^{*})=0$. By contradiction, we suppose that $-\infty<g_\infty<h_\infty<\infty$. It can be derived from Theorem 3.4 that $\lambda_1^D((g_\infty,h_\infty),H'(0))\geq0$. So $\lambda_1^D((-h_0,h_0),H'(0))>\lambda_1^D((g_\infty,h_\infty),H'(0))\geq 0$ hold according to the strictly monotonic decreasing of $\lambda_1^D((-g(t),h(t)),H'(0))$ in $t$, which contradicts to $\lambda_1^D((-h_0,h_0),H'(0))\leq0$.
\epf

\begin{thm}
Assume that $H'(0)< g^{*}$. Then $-\infty<g_\infty<h_\infty<+\infty$ for sufficiently small initial value
$||J_0||_{C([-h_0,h_0])}+||A_0||_{C([-h_0,h_0])}$ and species vanish.
\end{thm}
\bpf
The proof is similar as that of Theorem 3.3 by constructing an upper solution to problem \eqref{a01}.
Let
$$
\bar \eta(t)=h_0(1+\theta-\frac{\theta}{2}e^{-\gamma t}):=h_0\omega(t),\,\,t\geq0,
$$
$$
\bar W_1(t,x)=C_1e^{-\gamma t}\phi(t,\frac{hx}{\bar \eta(t)}),\,\,\bar W_2(t,x)=C_1e^{-\gamma t}\psi(t,\frac{hx}{\bar \eta(t)}),\,\,t\geq0,\,\,x\in[-h_0,h_0],
$$
where $(\phi(t,x),\psi(t,x))$ is a pair of eigenfunction of problem \eqref{c02} in $x\in(-h_0,h_0)$ corresponding to the principal eigenvalue $\lambda_1^D((-h_0,h_0),H'(0))$. Since $H'(0)<g^{*}$, it follows from the strictly decreasing of $\lambda_1^D((-h_0,h_0),H'(0))$ in $H'(0)$ in Lemma 3.1(ii) that $\lambda_1^D((-h_0,h_0),H'(0))> \lambda_1^D((-h_0,h_0),g^{*})=0$.
Careful calculations yield
$$
\bar W_{1t}-d_1\bar W_{1xx}-b(t)\bar W_2+[a(t)+m_1(t)]\bar W_1+\alpha_1(t)\bar W_1^2\geq 0
$$
and
$$
\bar W_{2t}-d_2\bar W_{2xx}-a(t)\bar W_1+m_2(t)\bar W_2+\alpha_2(t)\bar W_2^2\geq 0
$$
owing to the choice of $\theta$,$\gamma$ and $\lambda_1^D((-h_0,h_0),H'(0))>0$, for more details, please see proofs in Theorem 3.3. Moreover, impulsive conditions satisfy $\bar W_1((n\tau)^+,x)=\bar W_1(n\tau,x)$ and $\bar W_2((n\tau)^+,x)\geq H(W_2(n\tau,x))$.

The evident difference from proof of Theorem 3.3 is that a positive constant
$$C_1:=\frac{\gamma\theta h_0^2 (2+\theta)}{2\pi[\mu_1\max_{0\leq t\leq \tau }\alpha(t)+\mu_2\max_{0\leq t\leq \tau }\beta(t)]}$$
 can be chosen such that for $t\in(n\tau,(n+1)\tau]$,
$$
-\mu_1 \bar W_{1x}(t,\bar\eta(t))-\mu_2\bar W_{2x}(t,\bar\eta(t))\leq \bar\eta'(t)
$$
and
$$
-\mu_1 \bar W_{1x}(t,-\bar\eta(t))-\mu_2\bar W_{2x}(t,-\bar\eta(t))\geq -\bar\eta'(t)
$$
hold. Also, sufficiently small initial value $(J_0(x),A_0(x))$ guarantees $(J_0,A_0)(x)\leq (\bar W_1,W_2)(0,x)$. Therefore, it follows from the comparison principle that $(\bar W_1,\bar W_2;(-\bar\eta,\bar\eta))$ is an upper solution to \eqref{a01}. The rest proof is also analogous to last part of Theorem 3.3 and we omit it here.
\epf
\vspace{2mm}

Corollary 4.1 and Theorem 4.1 aforementioned show the classifications of $H'(0)\leq g_{*}$ for vanishing and $H'(0)\geq g^{*}$ for spreading of species, respectively. In the following, the impact of expanding capacities $\mu_1$ and $\mu_2$ will be investigated in the case of $g_{*}<H'(0)<g^{*}$, that is, spreading occurs for relatively large expanding capacities while vanishing happens for relatively small expanding capacities. Next, some preliminaries about $\lambda_1^D((g_\infty,h_\infty),H'(0))$ are firstly shown.

\begin{lem}
Suppose $\lambda_1^D((g_\infty,h_\infty),H'(0))<0$, then $-g_\infty=h_\infty=\infty$.
\end{lem}
\bpf
Since $\lambda_1^D((g_\infty,h_\infty),H'(0))<0$, we obtain $\lambda_1^D((-\infty,+\infty),H'(0))<0$. By contradiction, suppose that $h_\infty<+\infty$ or $-g_\infty>-\infty$. It then follows form Lemma 3.6 that $-\infty<g_\infty<h_\infty<\infty$, which together with Theorem 3.4, gives $\lambda_1^D((g_\infty,h_\infty),H'(0))\geq 0$. This contradicts to the fact $\lambda_1^D((g_\infty,h_\infty),H'(0))<0$.
\epf

\begin{lem}
Suppose $g_{*}<H'(0)<g^{*}$. Then there exists a positive constant $2L^*$ such that if $h_\infty-g_\infty>2L^*$, we have $\lambda_1^D((g_\infty,h_\infty),H'(0))<0$.
\end{lem}
\bpf
Recalling that $\lambda_1^D(\Omega,H'(0))$ is strictly decreasing with respect to $H'(0)$ and taking $\Omega:=(-\infty,+\infty)$ and $\Omega:=(-h_0,h_0)$, we have $\lambda_1^D((-\infty,+\infty),H'(0))<\lambda_1^D((-\infty,+\infty),g_{*})=0$ and $\lambda_1^D((-h_0,h_0),H'(0))>\lambda_1^D((-h_0,h_0),g^*)=0$, respectively. Therefore, we derive from $\lambda_1^D((-\infty,+\infty),H'(0))<0$ and $\lambda_1^D((-h_0,h_0),H'(0))>0$ that there exists some time $t_0\in (0,+\infty)$ such that $\lambda_1^D((g(t_0),h(t_0)),H'(0))=0$.

We then define $2L^*:=h(t_0)-g(t_0)$. It follows from Lemma 3.1(i) that $\lambda_1^D((g(t),h(t)),H'(0))\\=\lambda_1^D((0,h(t)-g(t)),H'(0))$ for any $t\geq0$, which together with the strictly decreasing of $\lambda_1^D((0,h(t)-g(t)),H'(0))$ with respect to $t$, yields $\lambda_1^D((0,h_\infty-g_\infty),H'(0))<0$ provided that $h_\infty-g_\infty>2L^*$.
\epf

\vspace{2mm}
We now show the impact of $\mu_1$ and $\mu_2$ on the spreading-vanishing dichotomy when $g_*<H'(0)<g^{*}$.
\begin{thm}\label{aas3}
Assume that $g_{*}<H'(0)<g^{*}$. There exists positive constants $\bar{\mu}$ and $\underline{\mu}$ corresponding to initial value $(J_0(x),A_0(x))$  such that\\
(i) if $\mu_1+\mu_2\geq \bar{\mu}$, then $-g_\infty=h_\infty=\infty$ and spreading occurs;\\
(ii) if $\mu_1+\mu_2\leq \underline{\mu}$, then $-\infty<g_\infty<h_\infty<+\infty$ and vanishing occurs.
\end{thm}
\bpf
(i) Lemma 2.2 gives $0<J(t,x),A(t,x)\leq C^*$ for $t\geq0$ and $g(t)<x<h(t)$, where $C^*$ is independent of $\mu_1$ and $\mu_2$. Further, for such $C^*$, we can find a positive constant $M^*:=\max\{\alpha^MC^*+ m_1^M+ a^M,\alpha_2^MC^*+m_2^M\}$ such that $b(t)A-a(t)J-m_1(t)J-\alpha_1(t)J^2\geq -M^*J$ and $a(t)J-m_2(t)A-\alpha_2(t)A^2\geq-M^*A$ for all $J,A\in[0,C^*]$.

Let us first consider the following auxiliary free boundaries problem
\begin{eqnarray}
\left\{
\begin{array}{lll}
u_t-d_1u_{xx}=-M^*u,  & t\in((n\tau)^+,(n+1)\tau],\,r(t)<x<s(t),\\[1mm]
v_t-d_2v_{xx}=-M^*v, & t\in((n\tau)^+,(n+1)\tau],\,r(t)<x<s(t),\\[1mm]
u((n\tau)^+,x)=u(n\tau,x),& r(n\tau)<x<s(n\tau), \\[1mm]
v((n\tau)^+,x)=H(v(n\tau,x)), & r(n\tau)<x<s(n\tau), \\[1mm]
u(t,x)=v(t,x)=0, & t\in(0,+\infty),\,x \in \{r(t),s(t)\}, \\[1mm]
r'(t)=-\mu_1u_x(t,r(t))-\mu_2v_x(t,r(t)), & t\in(n\tau,(n+1)\tau], \\[1mm]
s'(t)=-\mu_1u_x(t,s(t))-\mu_2v_x(t,s(t)), & t\in(n\tau,(n+1)\tau], \\[1mm]
r(0)=-h_0, s(0)=h_0, \\[1mm]
u(0,x)=u_0(x), v(0,x)=v_0(x),& x\in[-h_0,h_0],n=0,1,2,\dots \\[1mm]
\end{array} \right.
\label{d01}
\end{eqnarray}
Arguing as Lemma 2.2 and Theorem 2.3, one easily checks that problem \eqref{d01} admits a unique solution $(u,v;(r,s))$ with $r'(t)<0$ and $s'(t)>0$ for all $t>0$ and $t\neq n\tau$. Further, to stress the dependence of solution on expanding capability $\mu_i\,(i=1,2)$, we denote $(J,A;(g,h))$ by $(J^{\mu_i},A^{\mu_i};(g^{\mu_i},h^{\mu_i}))$ and $(u,v;(r,s))$ by $(u^{\mu_i},v^{\mu_i};(r^{\mu_i},s^{\mu_i}))$, respectively. The comparison principle asserts
\begin{equation}
(J^{\mu_i},A^{\mu_i})\geq(u^{\mu_i},v^{\mu_i}),\,h^{\mu_i}(t)\geq s^{\mu_i}(t),\,r^{\mu_i}(t)\geq g^{\mu_i}(t)\,\,\text{for}\,\,t\geq0\,\,\text{and}\,\,x\in[r^{\mu_i}(t),s^{\mu_i}(t)].
\label{d02}
\end{equation}

In what follows, we shall to prove that for sufficiently large $\mu_1+\mu_2$, there holds
\begin{equation}
r^{\mu_i}(\tau)\leq-L^*,\, s^{\mu_i}(\tau)\geq L^*.
\label{d03}
\end{equation}

Once \eqref{d03} is established, it then follows from \eqref{d02} that
$g_\infty=\lim\limits_{t\to\infty}{g^{\mu_i}(t)}<g^{\mu_i}(\tau)\leq r^{\mu_i}(\tau)\leq-L^*$ and $h_\infty=\lim\limits_{t\to\infty}{h^{\mu_i}(t)}>h^{\mu_i}(\tau)\geq s^{\mu_i}(\tau)\geq L^*$, so $h_\infty-g_\infty>2L^*$, which together with Lemma 4.4 that $\lambda_1^D((g_\infty,h_\infty),H'(0))<0$. Finally recalling Lemma 4.3, we easily derive that $-g_\infty=h_\infty=\infty$ and spreading happens.

It remains to prove \eqref{d03}. We first take the following initial boundary problem into account
\begin{eqnarray}
\left\{
\begin{array}{lll}
w_t-d_1w_{xx}=-M^*w,  & t\in(0^+,\tau],\,p(t)<x<q(t),\\[1mm]
z_t-d_2z_{xx}=-M^*z, & t\in(0^+,\tau],\,p(t)<x<q(t),\\[1mm]
w(0^+,x)=w(0,x),& p(0)<x<q(0), \\[1mm]
z(0^+,x)=H(z(0,x)), & p(0)<x<q(0), \\[1mm]
w(t,x)=z(t,x)=0, & t\in(0,+\infty),\,x \in\{p(t),q(t)\}, \\[1mm]
p'(t)=-\mu_1^*w_x(t,p(t))-\mu_2^*z_x(t,p(t)), & t\in(0,\tau], \\[1mm]
q'(t)=-\mu_1^*w_x(t,q(t))-\mu_2^*z_x(t,q(t)), & t\in(0,\tau], \\[1mm]
w(0,x)=\underline{u}_0(x), z(0,x)=\underline{v}_0(x),& x\in[-h_0/2,h_0/2], \\[1mm]
\end{array} \right.
\label{d04}
\end{eqnarray}
with $\mu_1\geq\mu_1^*$ and $\mu_2\geq\mu_2^*$. Initial value satisfies
\begin{eqnarray}
\left\{
\begin{array}{lll}
(0,0)<(\underline{u}_0(x),\underline{v}_0(x))<(u_0(x),v_0(x)),\,\,x\in(-h_0/2,h_0/2)\\[1mm]
\underline{u}_0(-h_0/2)=\underline{u}_0(h_0/2)=\underline{v}_0(-h_0/2)=\underline{v}_0(h_0/2)=0
\end{array} \right.
\label{d44}
\end{eqnarray}
with $\underline{u}_0'(-h_0/2)>0,\underline{u}_0'(h_0/2)<0,\underline{v}_0'(-h_0/2)>0$ and $\underline{v}_0'(h_0/2)<0$.  Also,
\begin{equation}
-p(0)=q(0)=h_0/2,\,\,-p(\tau)=q(\tau)=L^*,\,\,-p'(t)>0,\,\, q'(t)>0\,\,\text{for}\,\,t\in(0,\tau].
\label{d05}
\end{equation}
It is clear that the solution $(w,z;(p,q)):=(w^{\mu_i^*},z^{\mu_i^*};(p^{\mu_i^*},q^{\mu_i^*}))$ to problem \eqref{d04} exists and is unique according to the theory of parabolic equations. The Hopf boundary lemma asserts that $w_x(t,p(t))>0,w_x(t,q(t))<0,z_x(t,p(t))>0$ and $z_x(t,q(t))<0$. The choice of $\underline{u}_0,\underline{v}_0,p(t)$ and $q(t)$ implies that there exists some positive constant $\bar\mu:=\mu_1^*+\mu_2^*$ such that
\begin{equation}
q'(t)\leq-\mu_1w_x(t,q(t))-\mu_2z_x(t,q(t))\,\,\text{and}\,\,p'(t)\geq-\mu_1w_x(t,p(t))-\mu_2z_x(t,p(t)).
\label{d06}
\end{equation}
for $t\in[0,\tau]$, $\mu_1\geq\bar\mu$ and $\mu_2\geq\bar\mu$.

In view of \eqref{d01},\eqref{d04},\eqref{d44} and \eqref{d06}, one easily checks by the comparison principle that
$$
(w^{\mu_i^*}(t,x),z^{\mu_i^*}(t,x))\leq(u^{\mu_i}(t,x),v^{\mu_i}(t,x)),\,q^{\mu_i^*}(t)\leq s^{\mu_i}(t),\,p^{\mu_i^*}(t)\geq r^{\mu_i}(t)
$$
for $t\in[0,\tau]$ and $x\in[p(t),q(t)]$, which together with \eqref{d05} yields
$$
 s^{\mu_i}(\tau)\geq q^{\mu_i^*}(\tau)=L^*,\,\,r^{\mu_i}(\tau)\leq p^{\mu_i^*}(\tau)=-L^*
$$
and the desired result \eqref{d03} is deduced.

(ii) To verify the vanishing case in relatively small expanding capacities, we make a minor modification to Theorem 3.3.

In fact, we here choose $h=h_0$ in Theorem 3.3 and keep other notations unchanged. If
$$
0<\mu_1+\mu_2\leq \underline{\mu}:=\frac{\gamma\theta h_0^2 (2+\theta)}{2\pi C\max_{0\leq t\leq \tau}\{\alpha(t)+\beta(t)\}},
$$
then $(\hat W_1,\hat W_2;(-\hat \eta,\hat\eta))$ is an upper solution to \eqref{a01}. Further, we can use it to prove that $-\infty<g_\infty<h_\infty<\infty$ and vanishing occurs, more details can be clearly observed in Theorem 3.3.
\epf

\vspace{3mm}
Combining Theorems \ref{aas1}, \ref{aas3} and Corollary 4.1, we have the following results for two special cases.

\begin{cor}
For any given $\mu_1>0$, there exists $\mu^*\in[0,+\infty]$ depending on initial value $(J_0(x),A_0(x))$ such that when $0<\mu_2\leq\mu^*$, vanishing occurs, while spreading occurs for $\mu_2>\mu^*$. Moreover,

$(i)$ if $H'(0)\leq g_{*}$, then $\mu^*=+\infty$, which implies vanishing occurs for any $\mu_2>0$;

$(ii)$ if $H'(0)\geq g^{*}$, then $\mu^*=0$, which means spreading occurs for any $\mu_2>0$;

$(iii)$ if $g_{*}<H'(0)<g^{*}$, then $0<\mu^*<+\infty$ and a spreading-vanishing dichotomy related to $\mu^*$ exists.
\end{cor}

\begin{cor}
Suppose $\mu_1=\mu_2:=\mu$. We have the spreading-vanishing dichotomy with respect to expanding capacity, that is, there exists  $0\leq\mu^\triangle\leq+\infty$ related to initial value such that for $\mu>\mu^\triangle$, species spread, otherwise, species vanish for $0<\mu\leq\mu^\triangle$.
\end{cor}

\section{Numerical simulation}
In this section, numerical simulations are carried out to show the impact of impulsive intervention, which is induced by human control, then the influence of expanding capacities triggered by individuals are also characterized.

Noticing that the principal eigenvalue $\lambda_1^D$ can not be expressed explicitly, we first present its estimations. For simplicity, variable coefficients in problem \eqref{a01} are replayed by constant coefficients in simulations, that is, $b(t):=b,a(t):=a,m_1(t):=m_1,m_2(t):=m_2,\alpha_1(t)=\alpha_1,\alpha_2(t)=\alpha_2$, unless otherwise specified. Considering the following auxiliary principal eigenvalue problems
\begin{eqnarray}
\left\{
\begin{array}{lll}
\alpha'(t)+d_1\lambda_0\alpha(t)\geq b\beta(t)-(a+m_1)\alpha(t)+\underline{\lambda_1}\alpha(t),& t\in(0^+,\tau],\\[1mm]
\beta'(t)+d_2\lambda_0\beta(t)\geq a\alpha(t)-m_2\beta(t)+\underline{\lambda_1}\beta(t),& t\in(0^+,\tau],\\[1mm]
\alpha(0^+)\geq H'(0)\alpha(0),\,\beta(0^+)\geq H'(0)\beta(0),\\[1mm]
\alpha(0)=\alpha(\tau),\,\beta(0)=\beta(\tau)
\end{array} \right.
\label{e01}
\end{eqnarray}
and
\begin{eqnarray}
\left\{
\begin{array}{lll}
\alpha'(t)+d_1\lambda_0\alpha(t)\leq b\beta(t)-(a+m_1)\alpha(t)+\overline{\lambda_1}\alpha(t),& t\in(0^+,\tau],\\[1mm]
\beta'(t)+d_2\lambda_0\beta(t)\leq a\alpha(t)-m_2\beta(t)+\overline{\lambda_1}\beta(t),& t\in(0^+,\tau],\\[1mm]
\alpha(0^+)\leq\alpha(0),\,\beta(0^+)\leq \beta(0),\\[1mm]
\alpha(0)=\alpha(\tau),\,\beta(0)=\beta(\tau),
\end{array} \right.
\label{e02}
\end{eqnarray}
we choose $\alpha(t)=k\beta(t)$, where $k$ is a positive constant. Then direct calculations for \eqref{e01} and \eqref{e02} yield
\begin{eqnarray*}
\left\{
\begin{array}{lll}
\underline{\lambda_1}\leq \frac{-\ln H'(0)}{\tau}+d_1\lambda_0+a+m_1-\frac{b}{k},\\[1mm]
\underline{\lambda_1}\leq \frac{-\ln H'(0)}{\tau}+d_2\lambda_0-ak+m_2
\end{array} \right.
\end{eqnarray*}
and
\begin{eqnarray*}
\left\{
\begin{array}{lll}
\overline{\lambda_1}\geq d_1\lambda_0+a+m_1-\frac{b}{k},\\[1mm]
\overline{\lambda_1}\geq d_2\lambda_0-ak+m_2,
\end{array} \right.
\end{eqnarray*}
where $\lambda_0$ is defined in Lemma 3.1. It then follows from equations in \eqref{c03} with constant coefficients that
\begin{equation}
\lambda_1^D\geq \min\{-\frac{\ln H'(0)}{\tau}+a+m_1-\frac{b}{k},\, -\frac{\ln H'(0)}{\tau}-ak+m_2\}
\label{e03}
\end{equation}
and
\begin{equation}
\lambda_1^D\leq \max\{d_1\lambda_0+a+m_1-\frac{b}{k},\,d_2\lambda_0-ak+m_2\}.
\label{e04}
\end{equation}

\vspace{3mm}
Now we use numerical simulations and estimations \eqref{e03} and \eqref{e04} to support our theoretical findings about asymptotic behaviors of solution to problem \eqref{a01} aforementioned, and to show the complex impacts of harvesting pulse and expanding capacities. Some parameters are firstly fixed in problem \eqref{a01} with constant coefficients
\begin{equation}
d_1=0.1,\, d_2=0.2,\, b=9.5,\, a=7,\, m_1=3.5,\, m_2=6,\, \alpha_1=0.1,\,\alpha_2=0.5,\,k=0.88,\,h_0=2.
\label{00}
\end{equation}
The impulsive function $H(A)$, expanding capacities $\mu_1$ and $\mu_2$, and harvesting period $\tau$ are chosen later.

In order to investigate the impact of impulsive harvesting on spreading or vanishing of species, we take three cases into consideration: no pulse in Fig. 1, harvesting pulse with a relatively small harvesting rate in Fig. 2 and  harvesting pulse with a relatively large harvesting rate in Fig. 3, respectively.
\begin{figure}[ht]
\centering
\subfigure[]{ {
\includegraphics[width=0.28\textwidth]{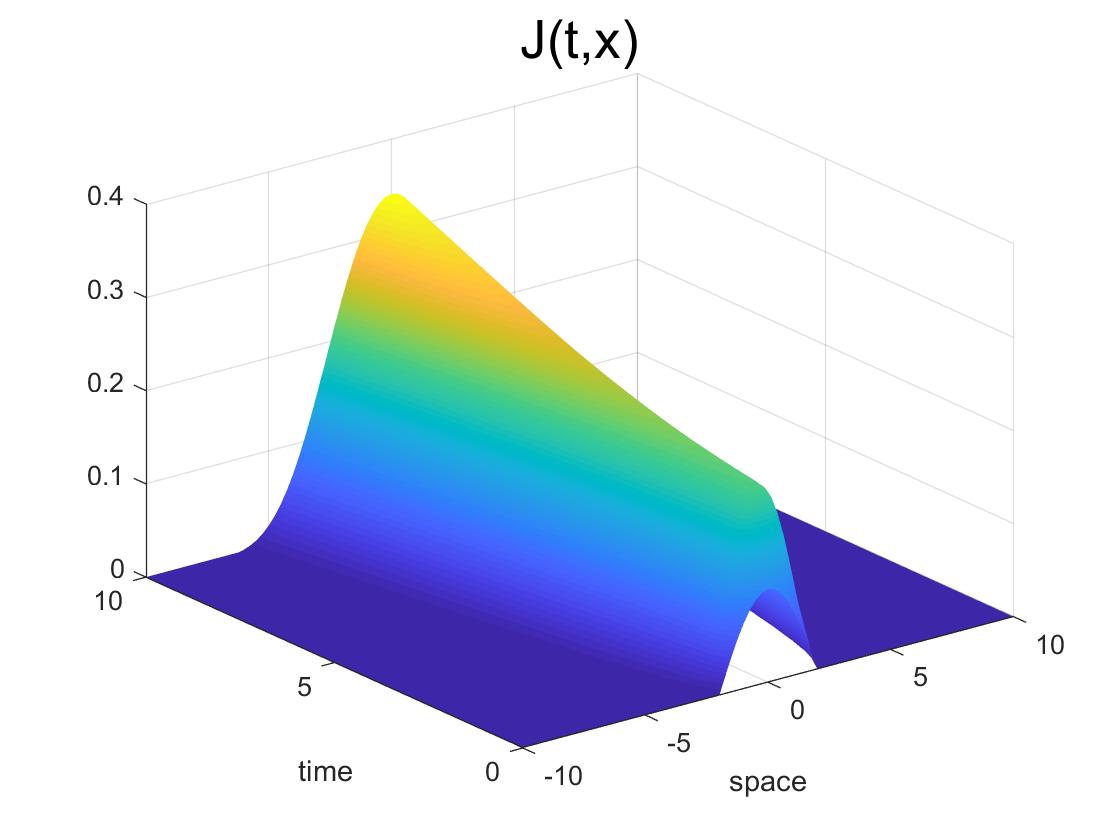}
} }
\subfigure[]{ {
\includegraphics[width=0.28\textwidth]{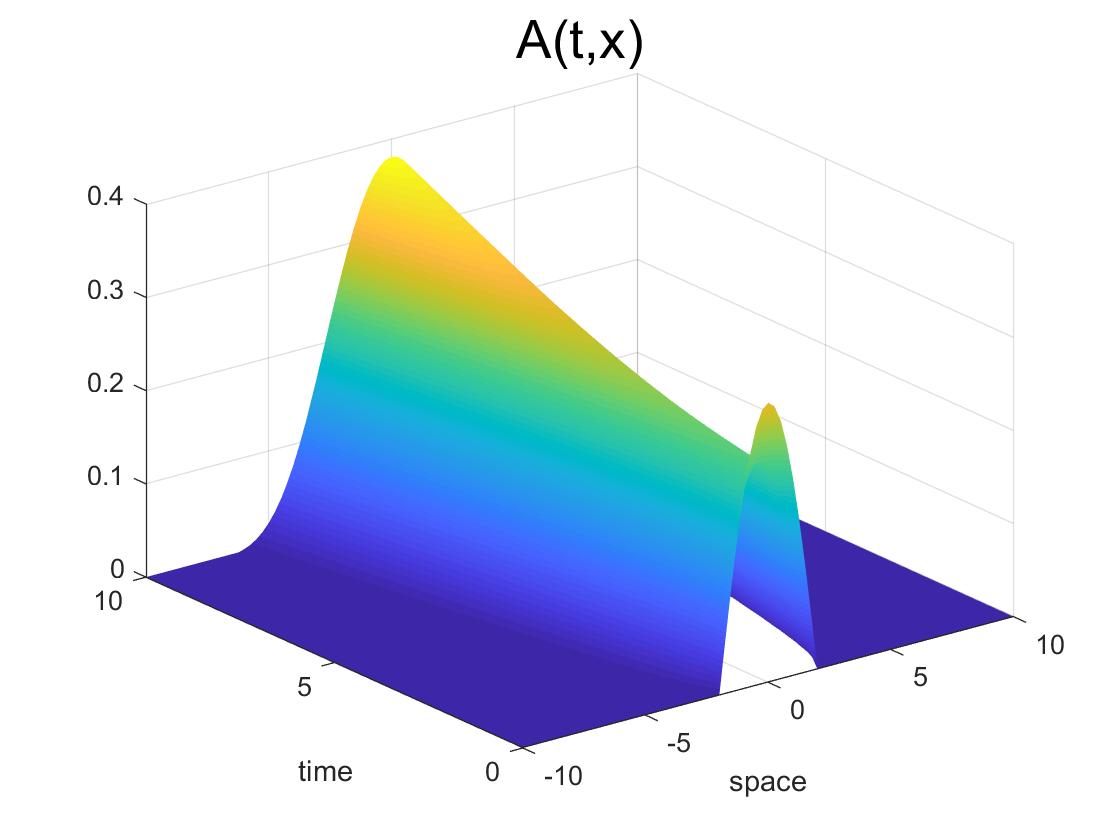}
} }

\subfigure[]{ {
\includegraphics[width=0.28\textwidth]{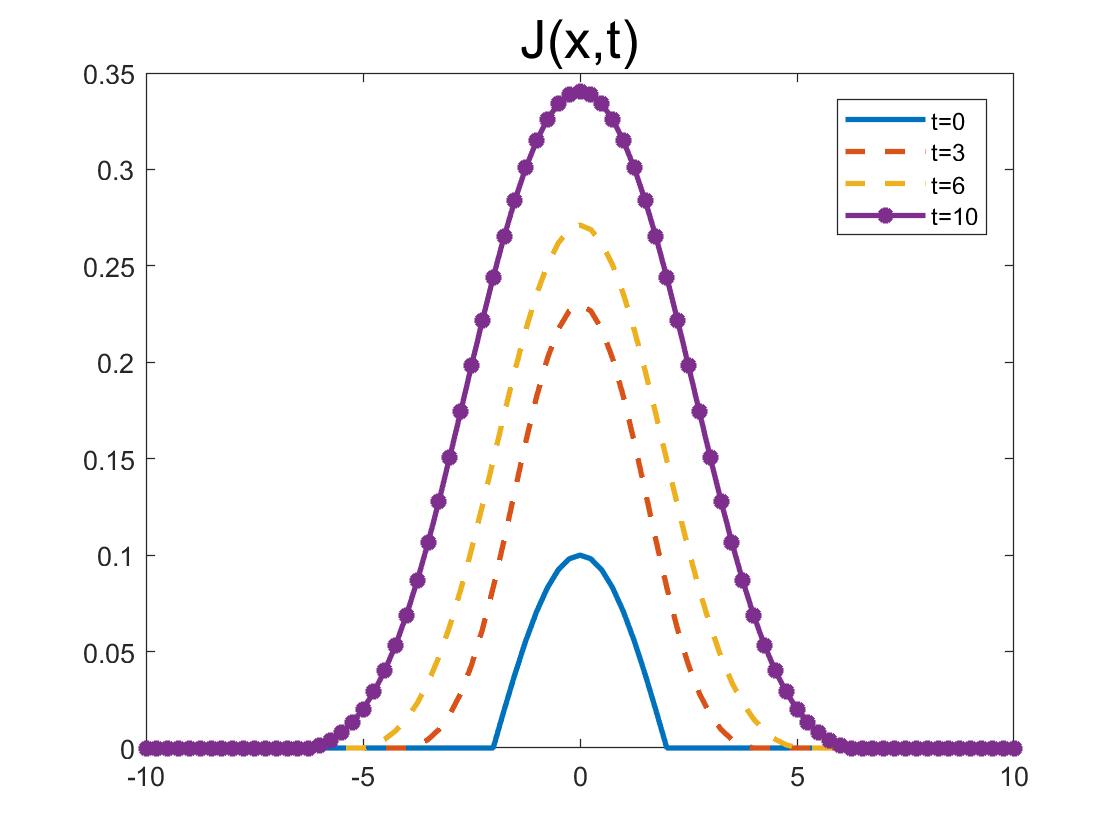}
} }
\subfigure[]{ {
\includegraphics[width=0.28\textwidth]{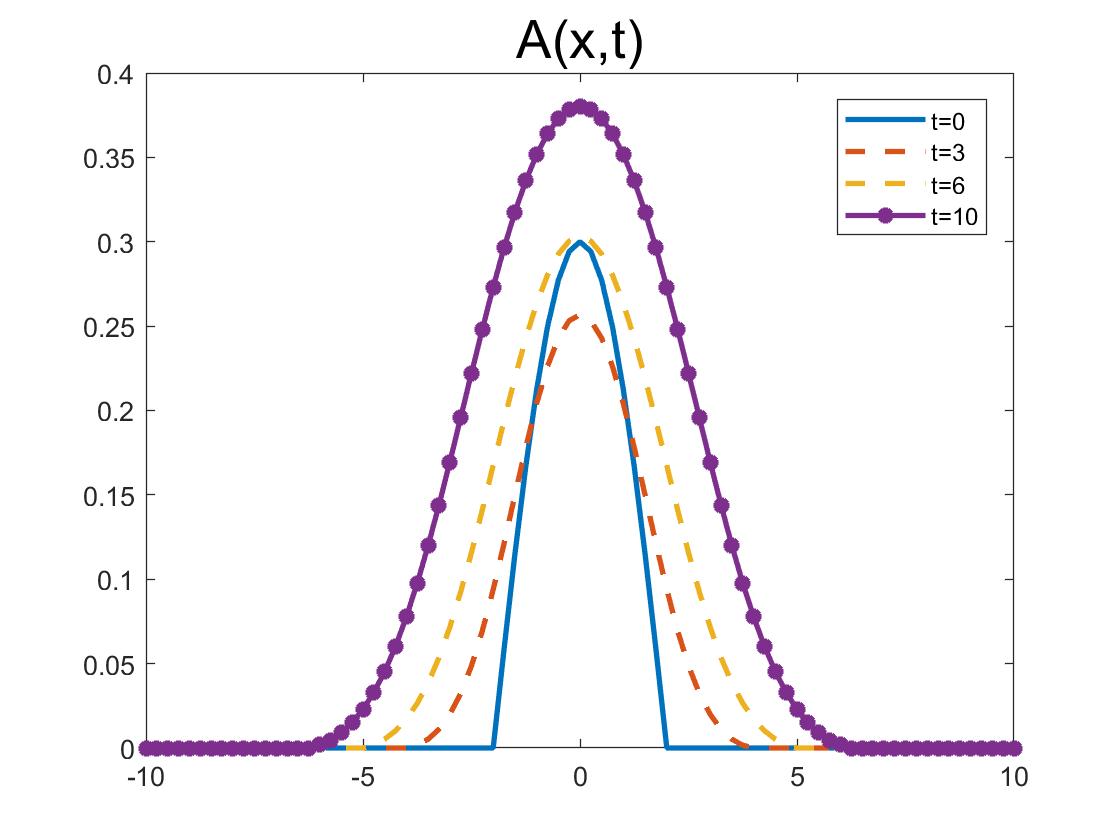}
} }
\caption{\scriptsize Solution of problem \eqref{a01} without harvesting pulse. All parameter values are given by \eqref{00}. Graphs (a) and (b) are stereogram of juveniles and adults, respectively.  Graphs (b) and (d) are sections of (c) and (d), respectively. It shows that juveniles and adults will spatially spread as time increases.
}
\label{tu1}
\end{figure}

One easily checks from Fig. 1 that juveniles $J(t,x)$ and adults $A(t,x)$ are spreading in the habitats without pulse. To investigate the impact of harvesting pulse and different harvesting rates on species spreading, we now take Beverton-Holt impulsive function $$H(A)=mA/(c+A)$$ as an example, where positive constants $0<m\leq c$ are concerned with the intensity of harvesting and to be chosen later, and $H(A)$ satisfies ($\mathcal{H}$1). It is worth noting that harvesting pulse will exert on the adults only.

\begin {exm}
Fix $\mu_1=10,\mu_2=15,\tau=2$ and choose impulsive function $H(A)=7A/(10+A)$ with a relatively small harvesting rate $1-H'(0)=1-m/c=0.3$, a direct calculation in \eqref{e04} yields
$$\begin{array}{llllll}
&&\lambda_1^D((-h_0,h_0),H'(0))\leq \max\{d_1\lambda_0+a+m_1-\frac{b}{k},\,d_2\lambda_0-ak+m_2\}\\[2mm]
&\approx& -0.0366<0.
\end{array}$$
\end {exm}

\begin{figure}
\centering
\subfigure[]{ {
\includegraphics[width=0.28\textwidth]{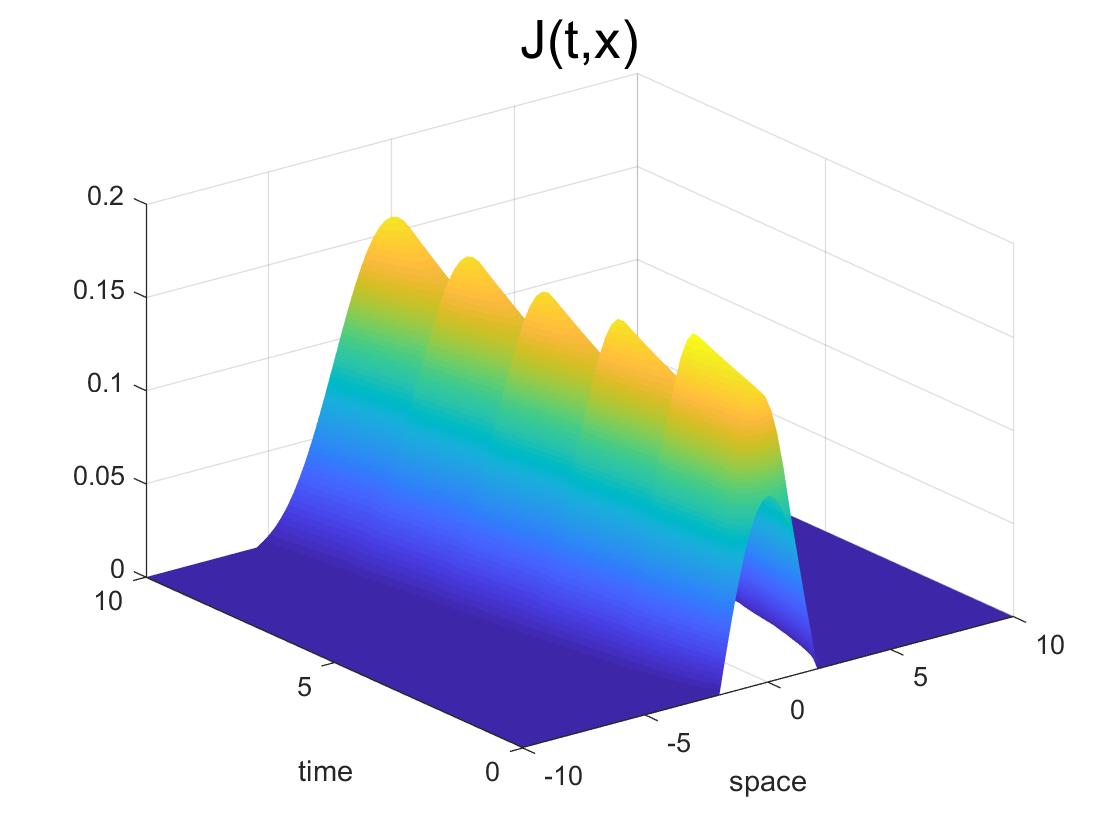}
} }
\subfigure[]{ {
\includegraphics[width=0.28\textwidth]{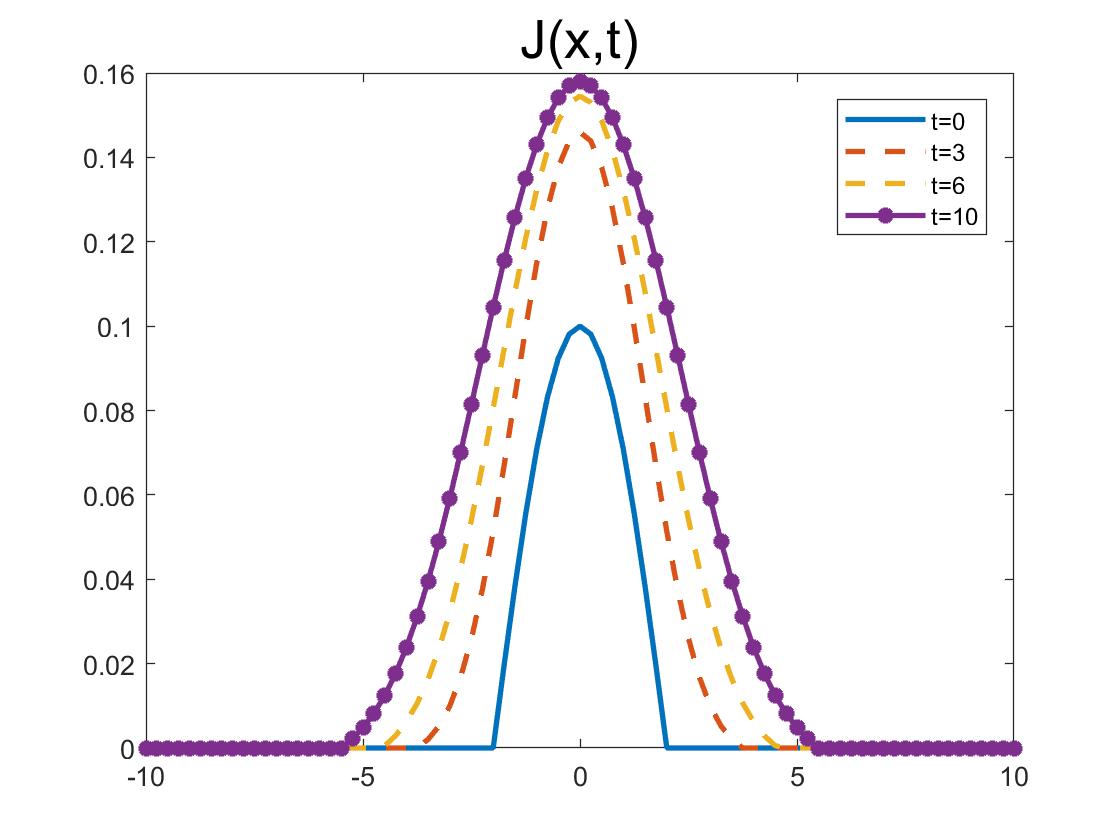}
} }
\subfigure[]{ {
\includegraphics[width=0.28\textwidth]{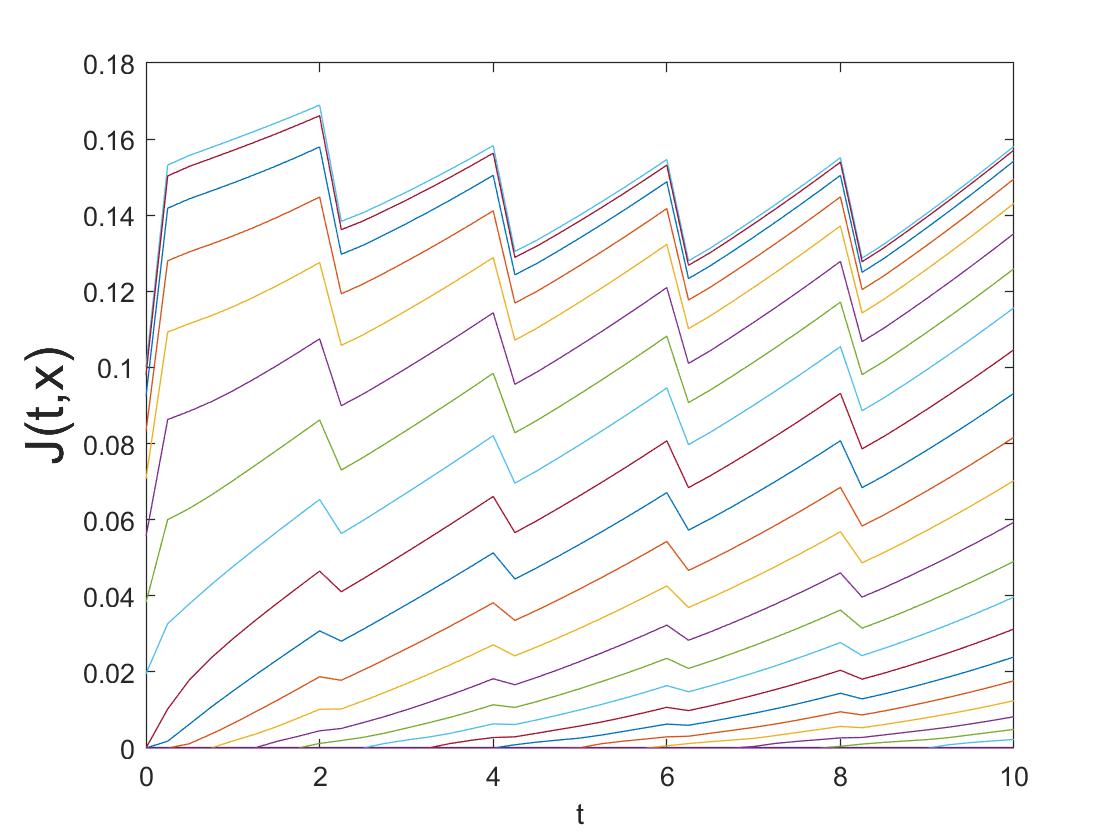}
} }

\subfigure[]{ {
\includegraphics[width=0.28\textwidth]{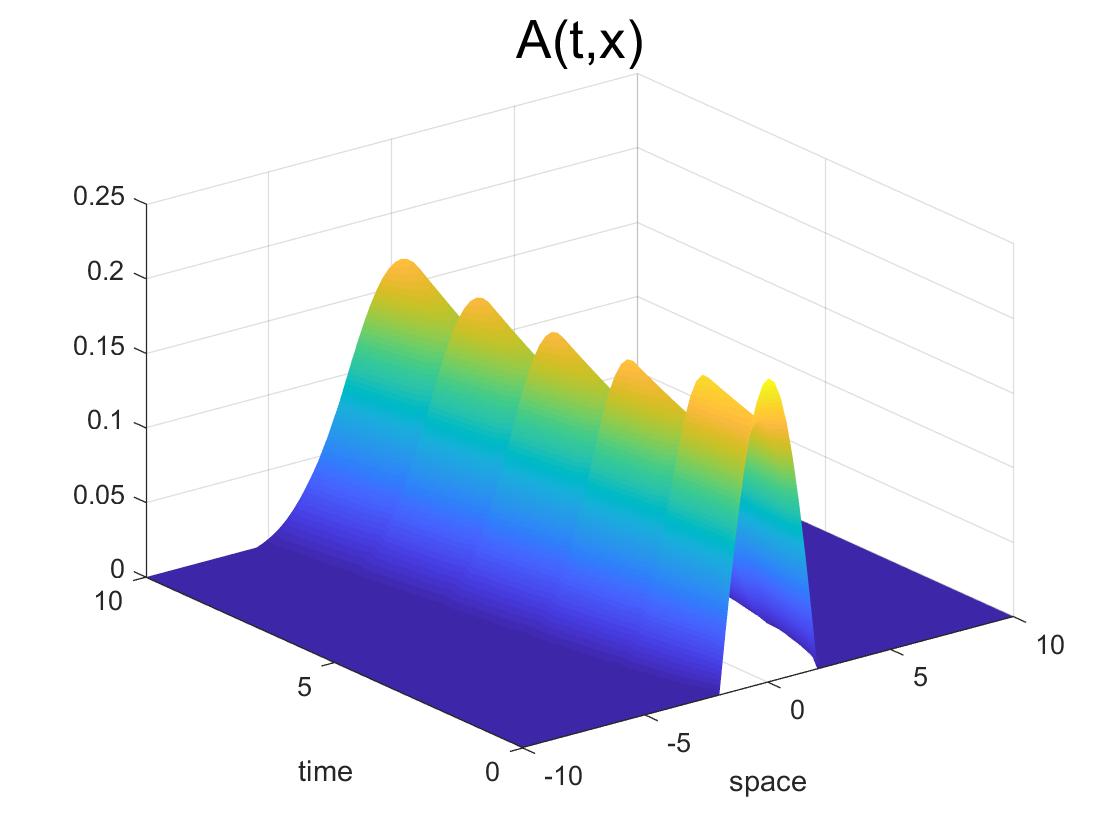}
} }
\subfigure[]{ {
\includegraphics[width=0.28\textwidth]{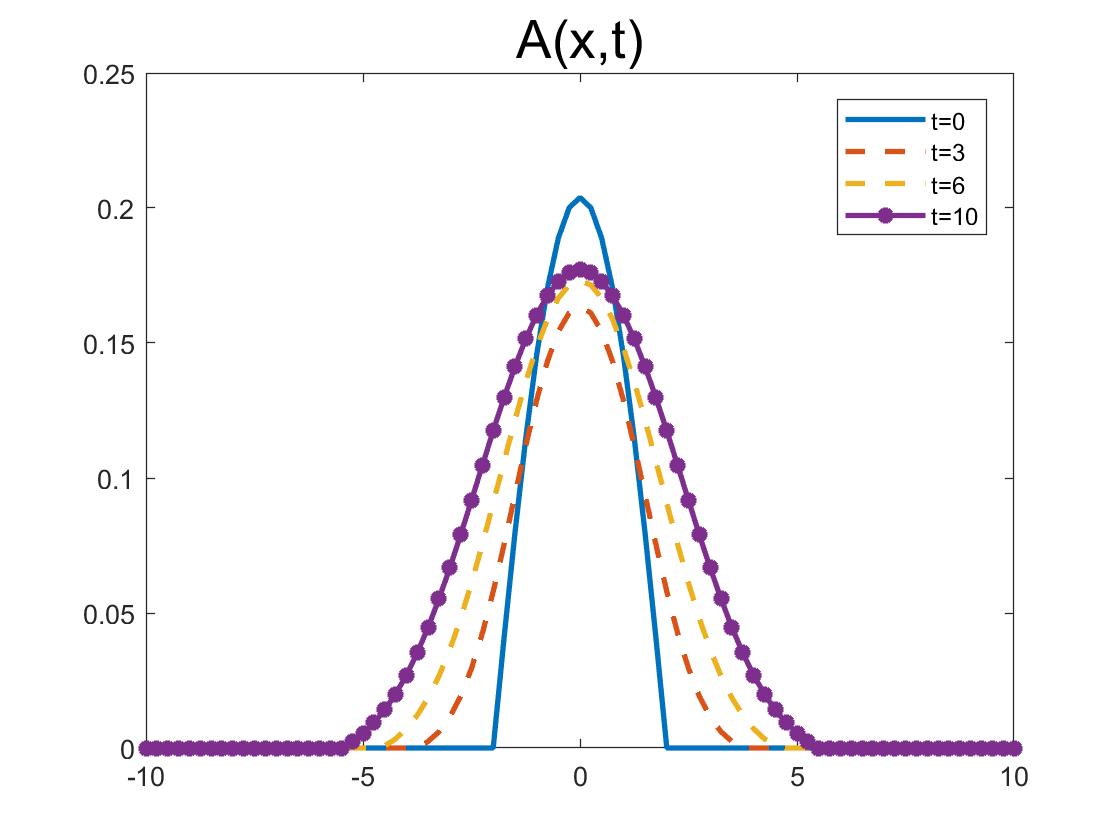}
} }
\subfigure[]{ {
\includegraphics[width=0.28\textwidth]{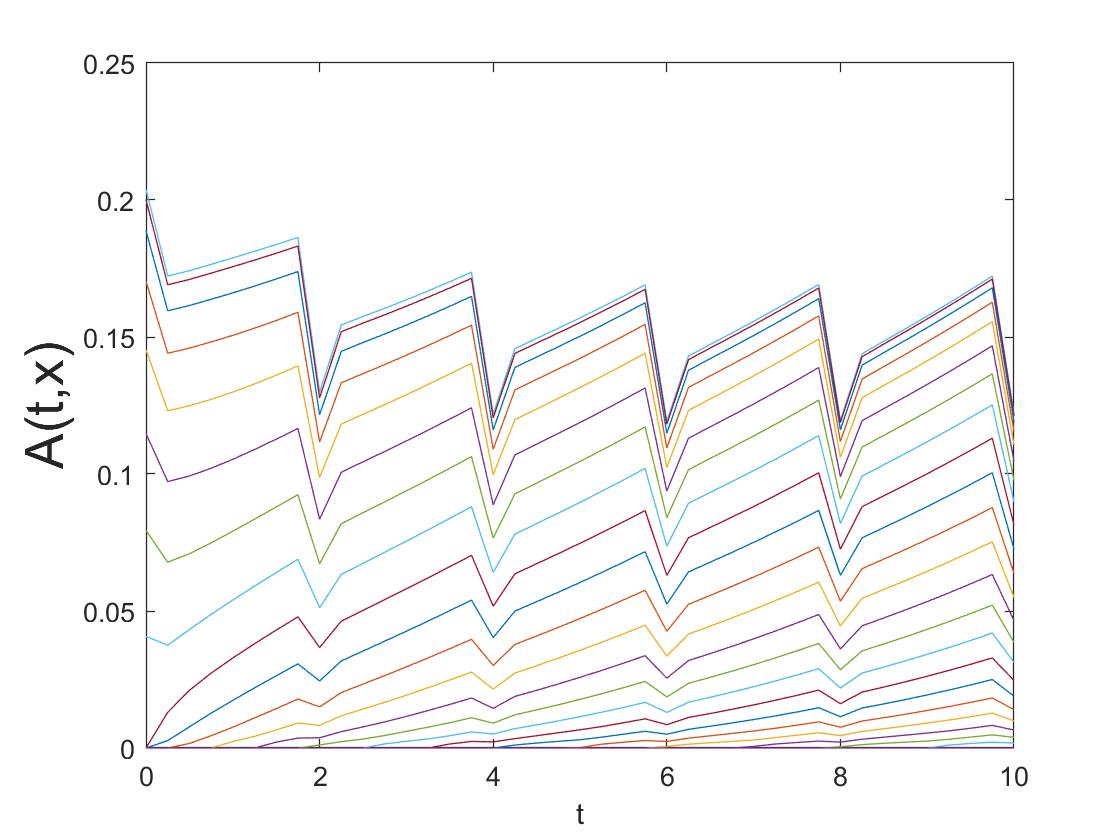}
} }
\caption{\scriptsize A simulation with small harvesting rate ($m=7,c=10$). Graphs (c) and (f) are the projection of (a) and (d) in the right view, respectively. Graphs (a)-(f) exhibit that the species persist eventually with small harvesting rate, and the periodically impulsive harvesting can be vividly observed at every time $t=0,2,4,\dots$ in Graph (f).
}
\label{tu2}
\end{figure}

It follows from $\lambda_1^D((-h_0,h_0),H'(0))\approx-0.0366(<0)$ and Theorem 4.1 that spreading happens in a relatively small harvesting rate, which is presented in Fig. 2. In comparison with Figs 1 and 2, one easily observes that although the species spatially spread in Fig. 2, the density has decreased from Fig. 1 to Fig. 2 owing to the periodic pulse intervention. We wonder whether the state of species will change if a greater harvesting rate exert on the adult? So we will choose pulse function with a relatively big harvesting rate in Example 5.2.

\begin {exm}
Fix $\mu_1=10,\mu_2=15,\tau=2$ and set impulsive function $H(A)=2A/(10+A)$ with a relatively large harvesting rate $1-H'(0)=0.8$, it follows from \eqref{e03} that
$$\begin{array}{llllll}
&&\lambda_1^D((-\infty,+\infty),H'(0))\geq \min\{-\frac{\ln H'(0)}{\tau}+a+m_1-\frac{b}{k},\, -\frac{\ln H'(0)}{\tau}-ak+m_2\}\\[2mm]
&\approx& 0.5093>0.
\end{array}$$
\end {exm}

Theorem 3.3 with $\lambda_1^D((-\infty,+\infty),H'(0))\approx 0.5093(>0)$ shows that the densities of juveniles and adults will decrease to zero, which is in accordance with Fig. 3. Distinctly, Figs. 1 and 3 verify our guess, that is, harvesting pulse helps us restrain the spreading of invasive species, then pulse with large harvesting rate even can change the state of species, from existence to extinction.
\begin{figure}
\centering
\subfigure[]{ {
\includegraphics[width=0.28\textwidth]{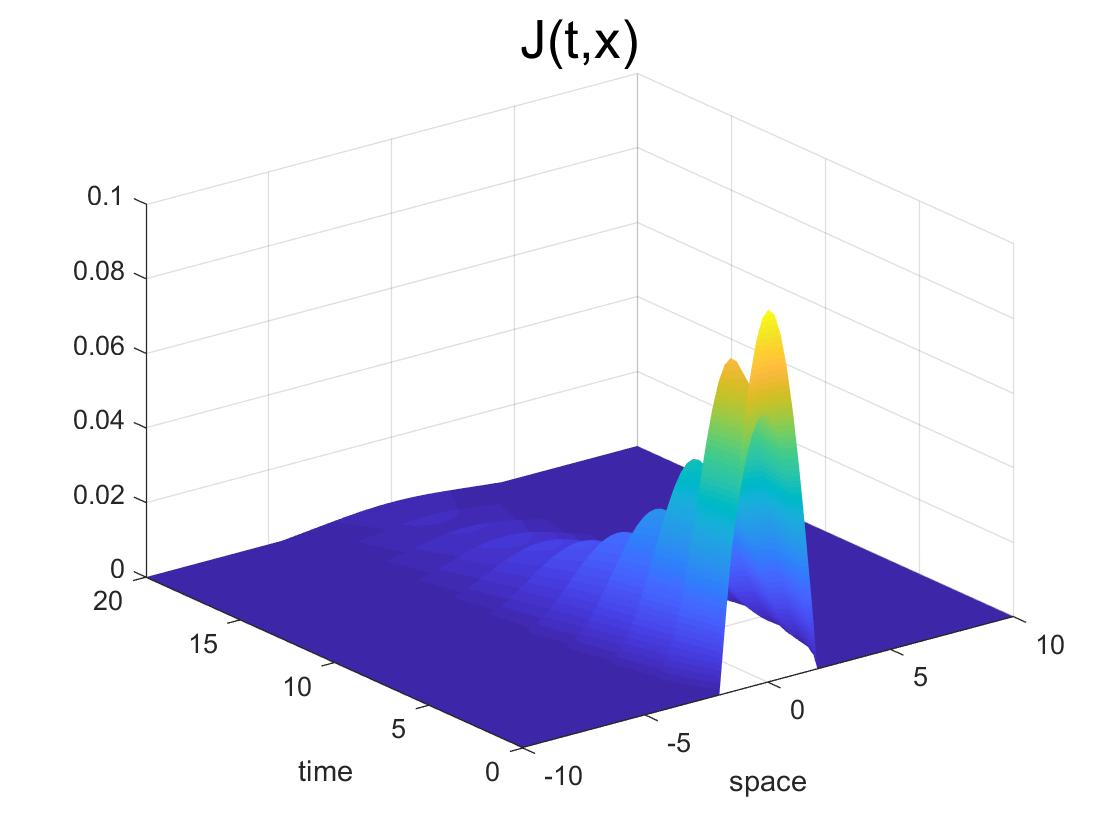}
} }
\subfigure[]{ {
\includegraphics[width=0.28\textwidth]{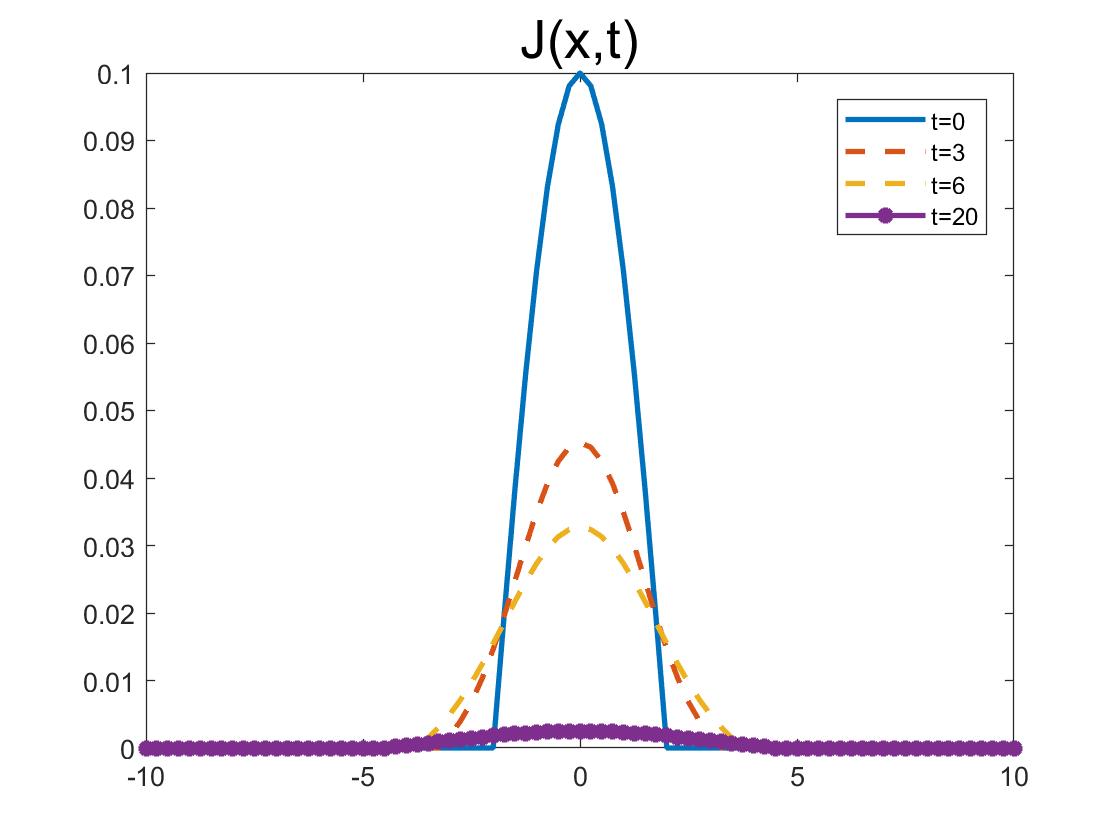}
} }
\subfigure[]{ {
\includegraphics[width=0.28\textwidth]{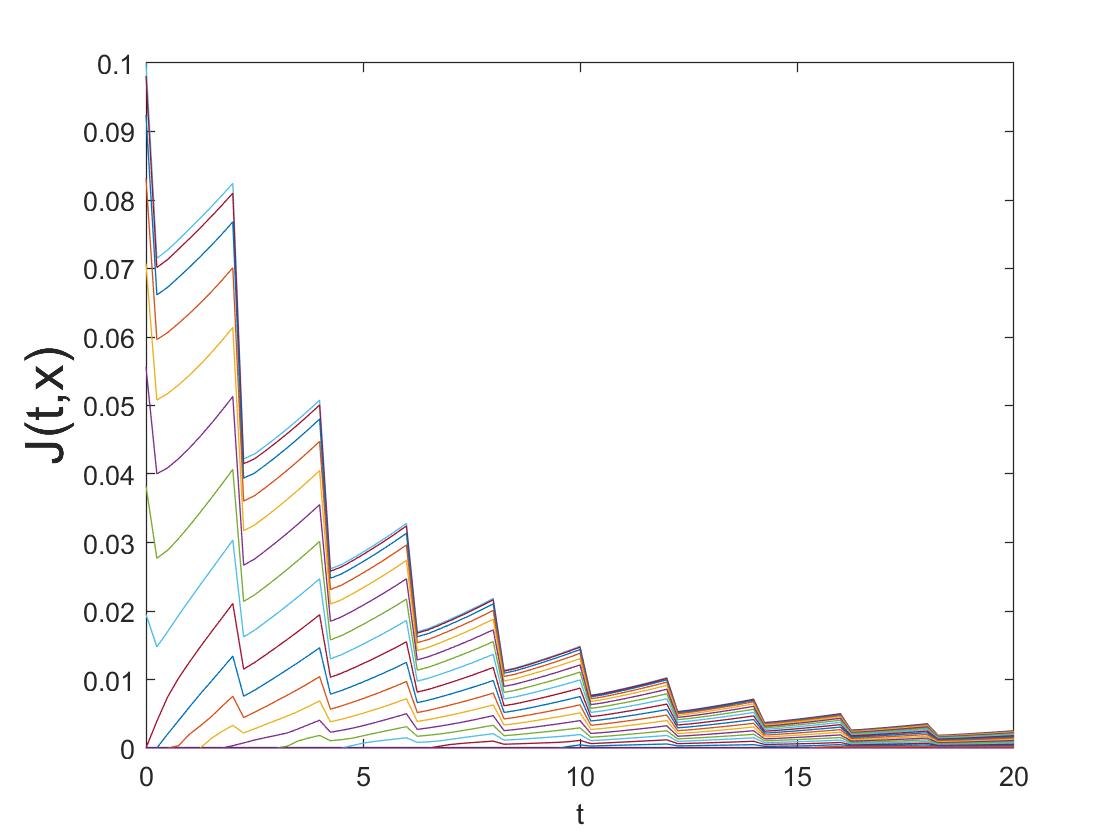}
} }

\subfigure[]{ {
\includegraphics[width=0.28\textwidth]{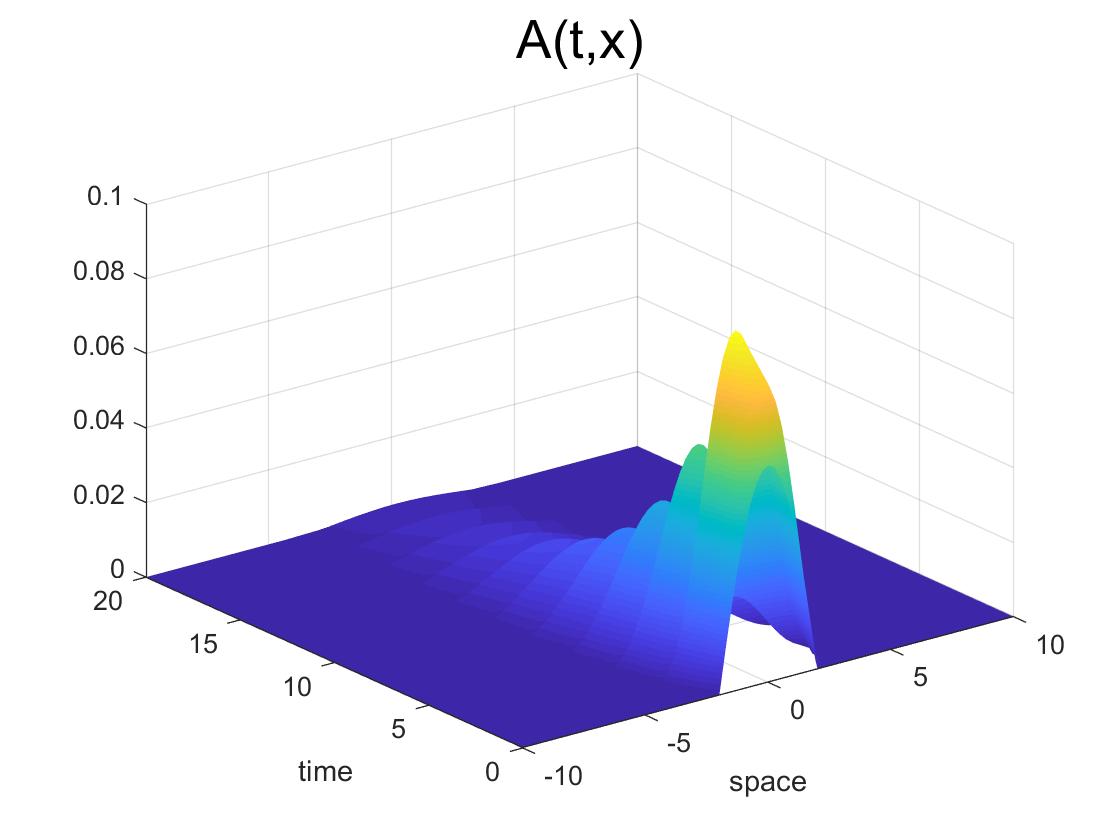}
} }
\subfigure[]{ {
\includegraphics[width=0.28\textwidth]{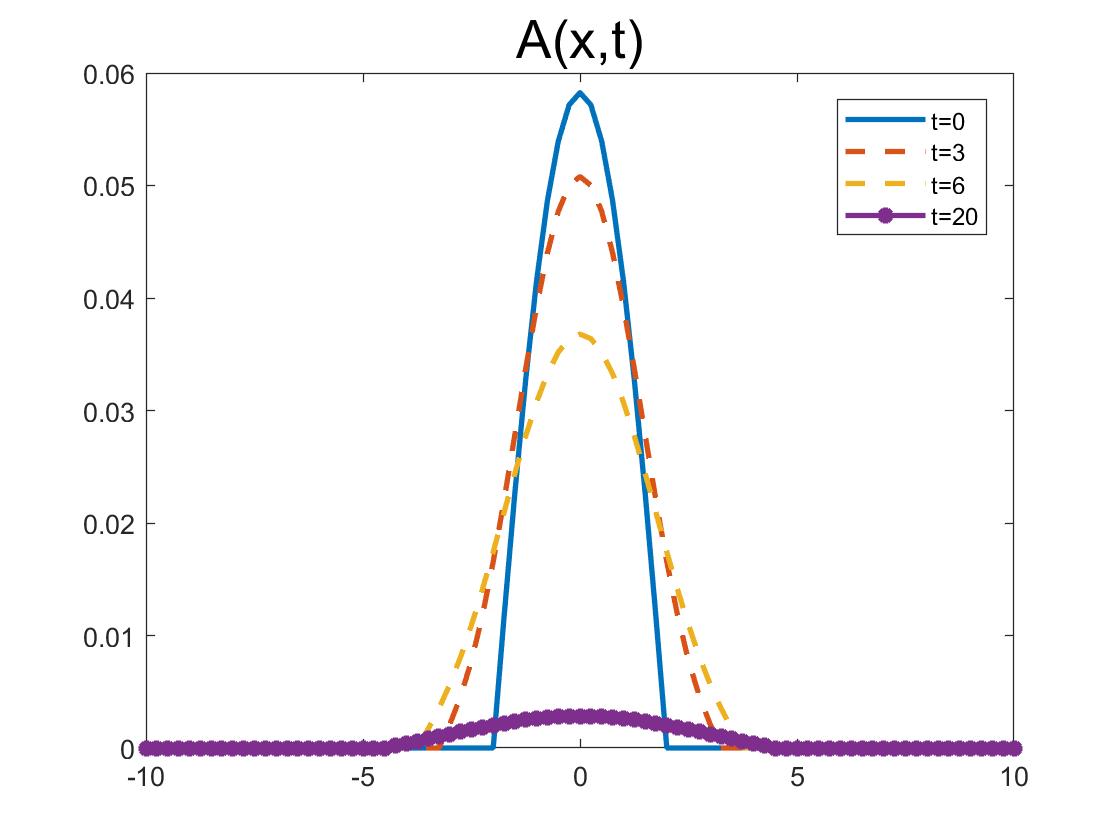}
} }
\subfigure[]{ {
\includegraphics[width=0.28\textwidth]{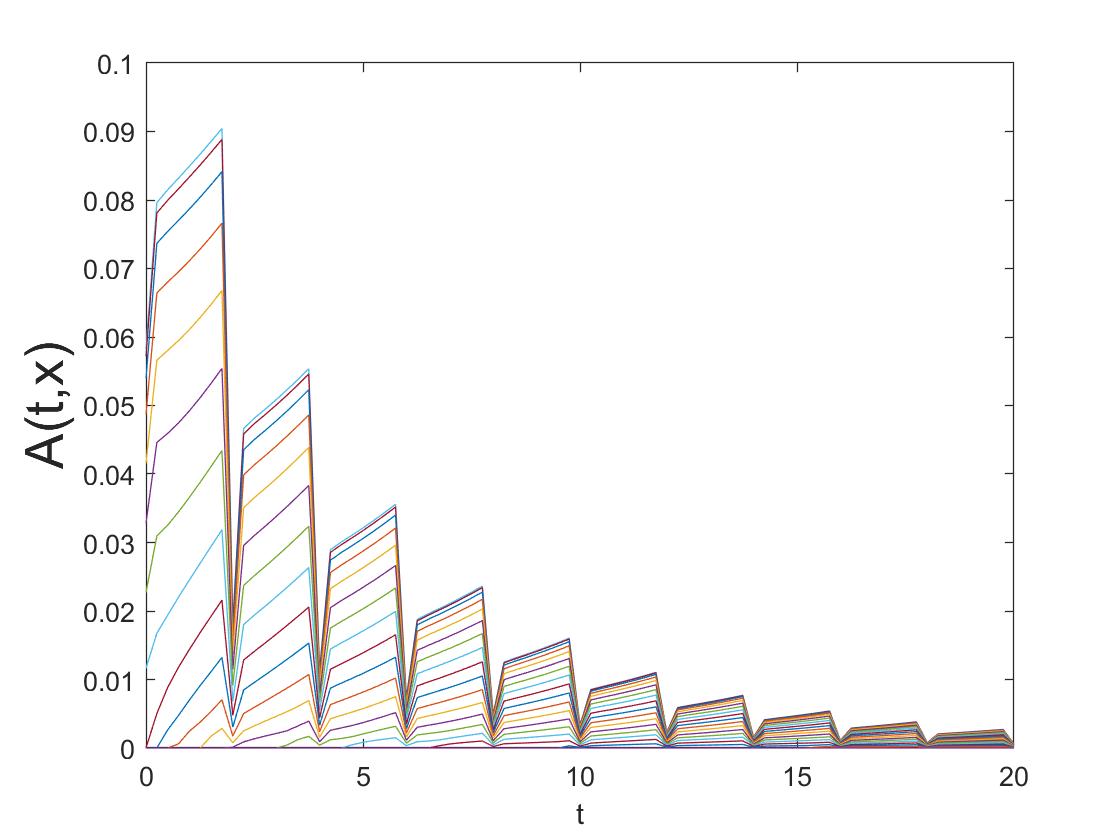}
} }
\caption{\scriptsize A simulation with large harvesting rate($m=2,c=10$). Graph (f) clearly indicates the density of the adult $A(t,x)$ decays to zero with harvesting pulse appearing at every time $t=0,2,4,\dots$, while Graph (c) shows that the juvenile $J(t,x)$ will also tend to vanish since the reciprocity to $A(t,x)$.
}
\label{tu3}
\end{figure}

The spatial distribution of juveniles and adults with pulse are shown in Figs 2 and 3. In view of Fig. 2(f)(or Fig. 3(f)), we clearly see that periodic pulse occurs in adults at every time $t=0,2,4,\dots$, which cause a sudden and sharp drop on density of adults. The cooperative relationship between juveniles and adults also leads to a relative gentle decline of juveniles, see Fig 2.(c)(or Fig 3.(c)). In comparison with Fig. 1, approximations in Fig. 2 indicate that species will still persist under small harvesting rate while eventually extinct under large harvesting rate in Fig. 3.

It follows from Figs 1-3 that pulse control from human intervention can affect the dynamical behaviors of species, that is, harvesting pulse can inhibit species spreading. The larger the harvesting rate, the more obvious the inhibition effect. Moreover, harvesting pulse even leads to the transition of species situations, from spreading (Fig. 1) to vanishing (Fig. 3).

\vspace{2mm}
Analytical result in Theorem 4.5 indicates that the spreading or vanishing of species not only depends on pulse control form human intervention, but also on expanding capacities of species itself. So in order to discuss the impact of expanding capacities on species spread or vanish, we first fix impulsive harvesting function $H(A)=\frac{6A}{10+A}$.

\begin{figure}
\quad\quad
\begin{minipage}{0.25\linewidth}
\centerline{\includegraphics[width=4.5cm]{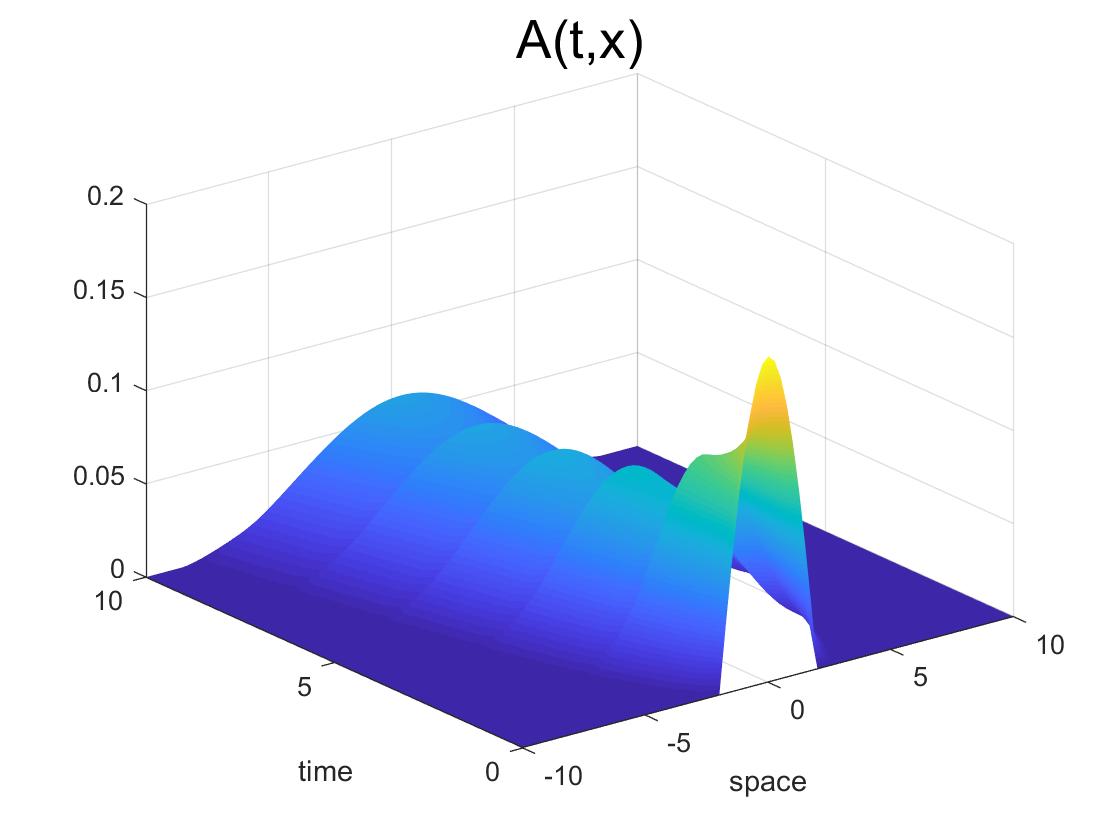}}
\centerline{\small{(a) $\mu_1=30,\mu_2=40$}}
\end{minipage}
\quad\quad\quad
\begin{minipage}{0.25\linewidth}
\centerline{\includegraphics[width=4.5cm]{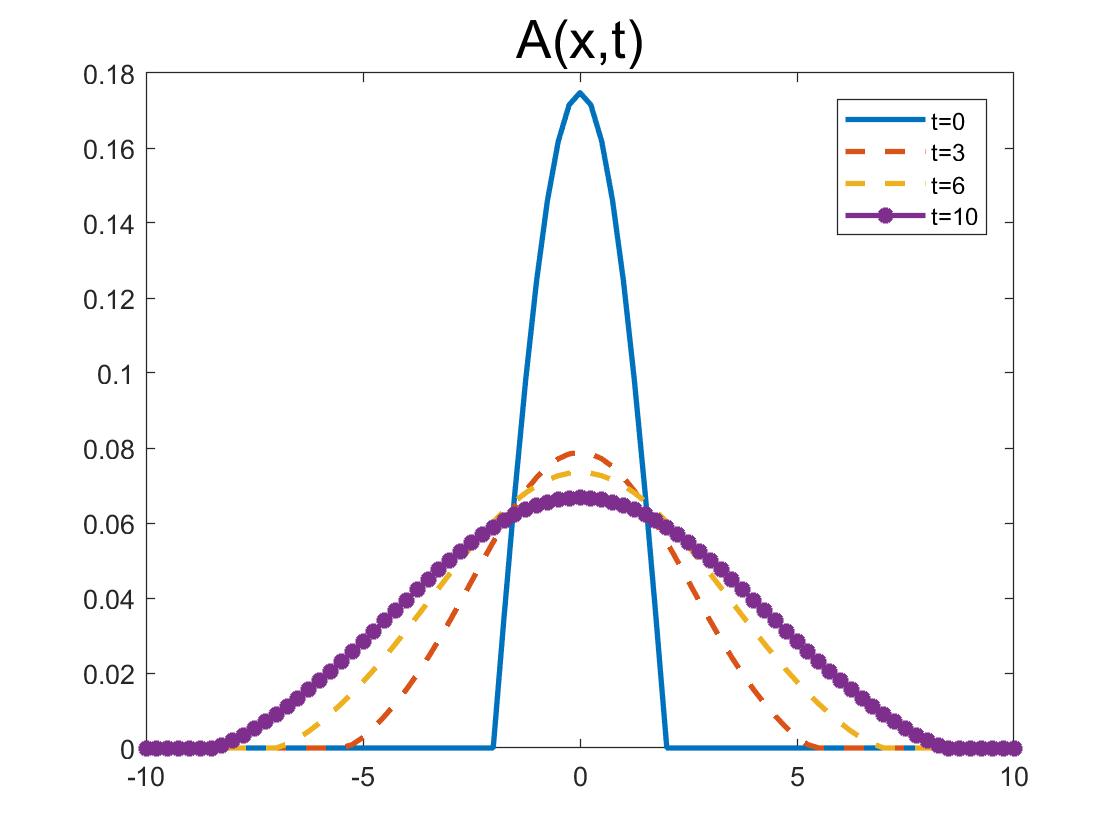}}
\centerline{\small{(b) $\mu_1=30,\mu_2=40$}}
\end{minipage}
\quad\quad\quad
\begin{minipage}{0.25\linewidth}
\centerline{\includegraphics[width=4.5cm]{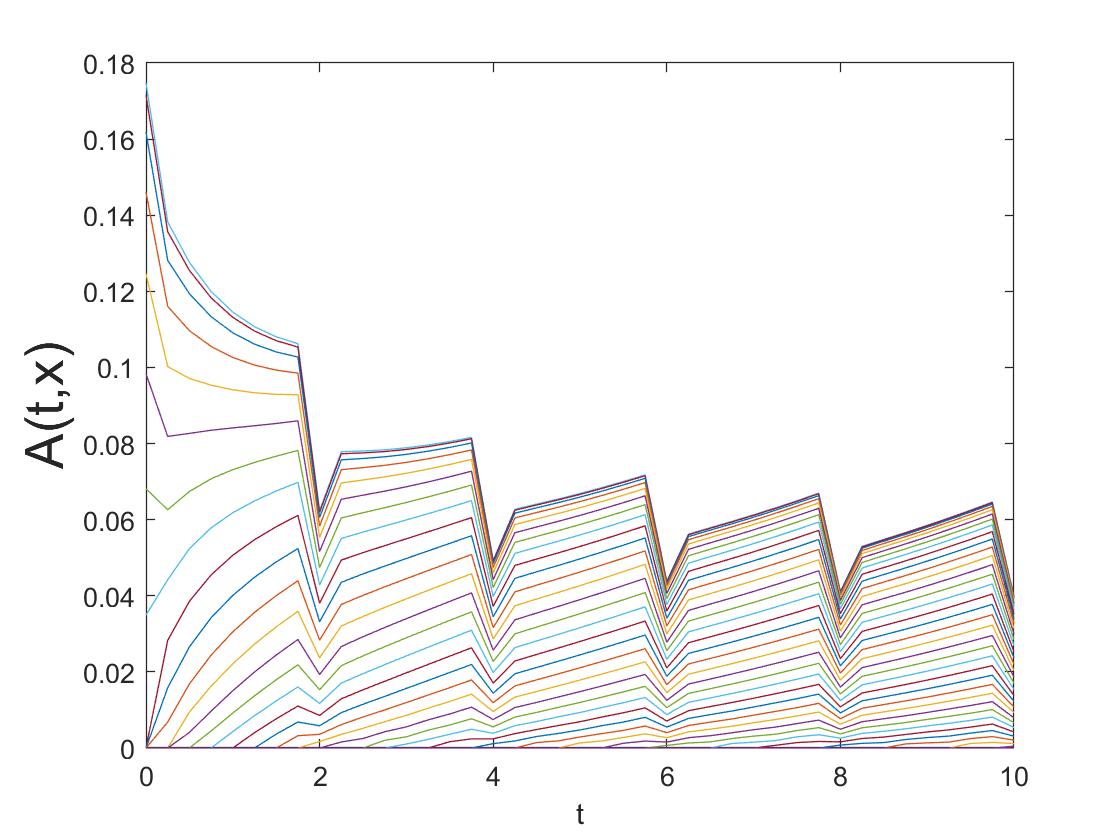}}
\centerline{\small{(c) $\mu_1=30,\mu_2=40$}}
\end{minipage}

\quad\quad
\begin{minipage}{0.25\linewidth}
\centerline{\includegraphics[width=4.5cm]{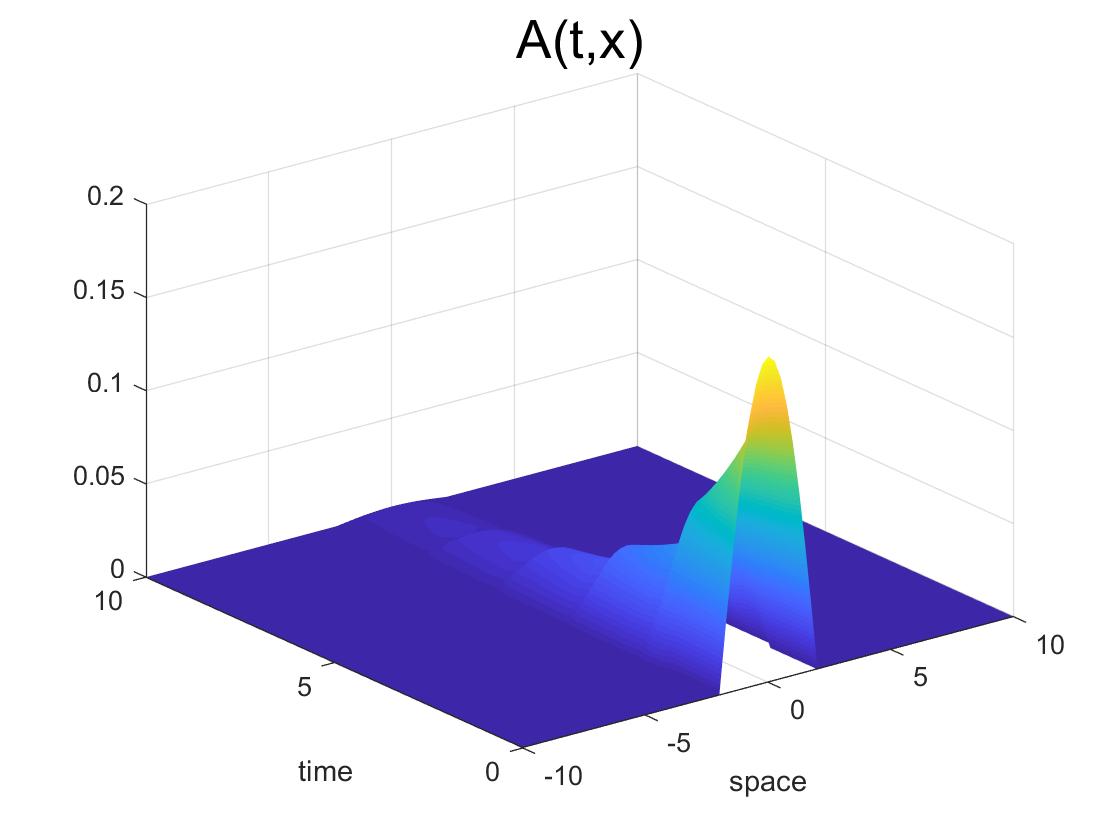}}
\centerline{\small{(d) $\mu_1=0.5,\mu_2=0.6$}}
\end{minipage}
\quad\quad\quad
\begin{minipage}{0.25\linewidth}
\centerline{\includegraphics[width=4.5cm]{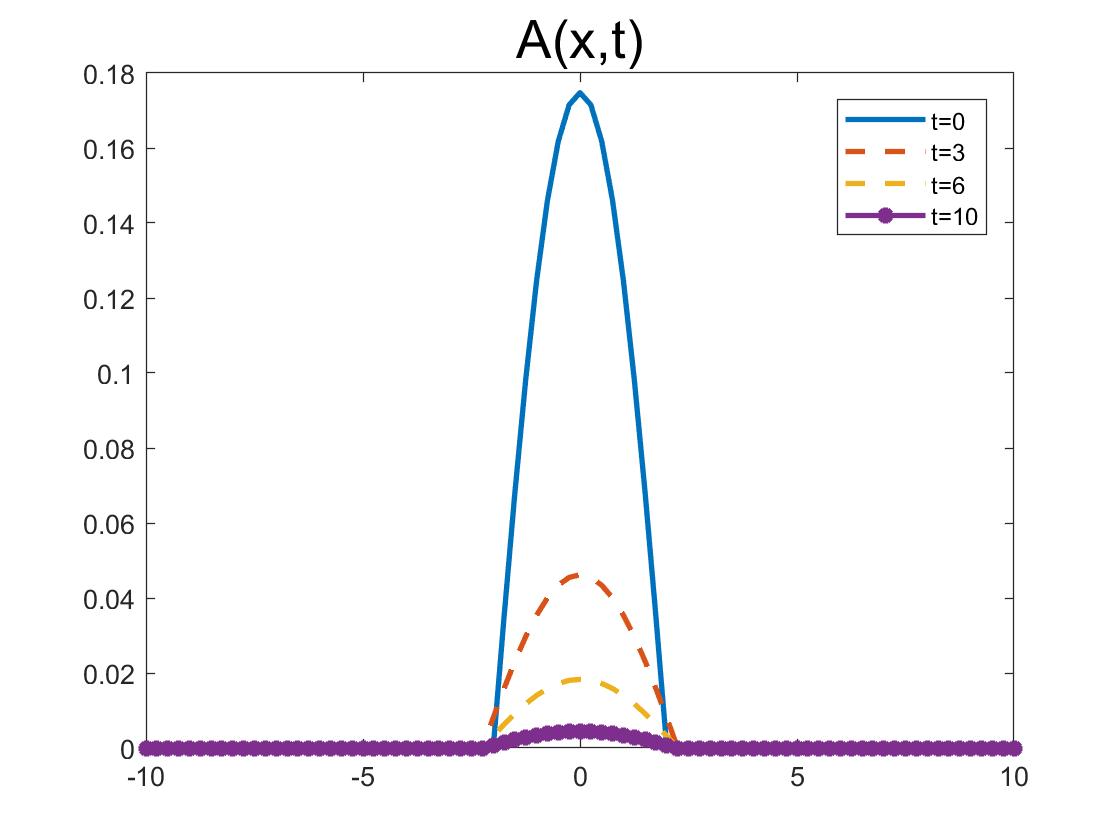}}
\centerline{\small{(e) $\mu_1=0.5,\mu_2=0.6$}}
\end{minipage}
\quad\quad\quad
\begin{minipage}{0.25\linewidth}
\centerline{\includegraphics[width=4.5cm]{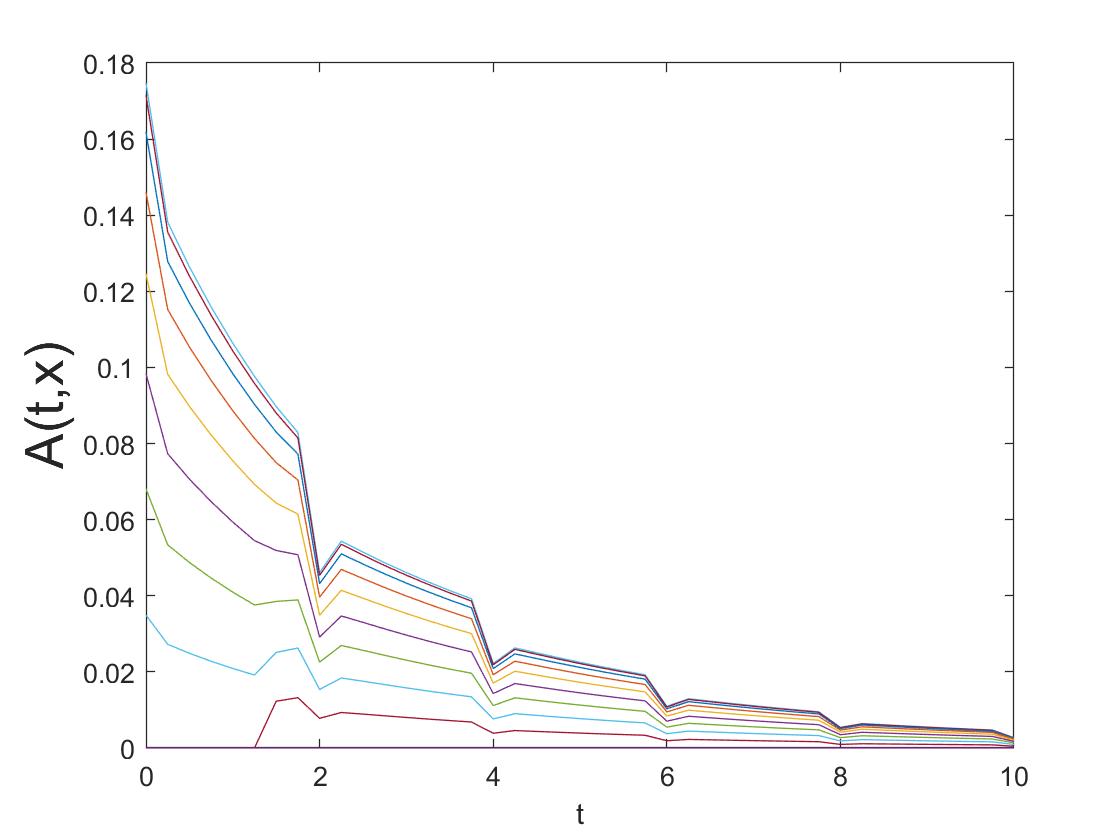}}
\centerline{\small{(f) $\mu_1=0.5,\mu_2=0.6$}}
\end{minipage}
\caption{\scriptsize The density of adults of \eqref{a01} for different choices of $\mu_1$ and $\mu_2$, with parameters $d_1=0.59,\, d_2=1,\, b=9.5,\, a=7,\, m_1=3.5,\, m_2=6,\, \alpha_1=0.1,\,\alpha_2=0.5$ and $\tau=2$. We here take the growing tendency of adults $A(t,x)$ as an example since the cooperative relationship of juveniles and adults. The dynamics of \eqref{a01} varies from spreading in (a)-(c) for relatively big expanding capacities with $\mu_1=30$ and $\mu_2=40$, to vanishing in (d)-(f) for small expanding capacities with $\mu_1=0.5$ and $\mu_2=0.6$.}
\label{tu4}
\end{figure}

One can see from Fig. 4 that when impulsive intervention occurs, the species with large expanding capacities (see (a)-(c) in Fig. 4) will spread while the species become extinct eventually with small expanding capacities (see (d)-(f) in Fig. 4), which enlighten us that when we adopt human control such as impulsive harvesting to ensure sustainable ecological development, the species' own expanding capacities cannot be ignored too. Moreover, a comparison of (b) and (d) in Fig. 4 indicates that the larger the expanding capacities, the more quickly the moving
boundaries expand under the harvesting pulse. On the other hand, the harvesting pulse can slow down the spreading speed of individuals, as is shown in Fig. 1(d) without pulse and Fig. 4(b) with harvesting pulse. It is observed in the process of simulations that the harvesting pulse from human intervention plays a greater role than expanding capacities from species itself concerning the dynamics of species and spreading speeds of moving boundaries.

Besides intensity of harvesting pulse and expanding capacities, we also find that harvesting timing (i.e.
harvesting period $\tau$) can also change the spreading or vanishing of individuals. Compared with $\tau=2$ in the figures above, we will choose $\tau=1$ and $\tau=5$ in Example 5.3 (Fig. 5), respectively.
\begin {exm}
Choose $\tau=1$ and $\tau=5$, respectively. Fix impulsive function and other parameters the same as that in Example 5.2.
\end {exm}
\begin{figure}
\quad\quad\quad\quad
\begin{minipage}{0.4\linewidth}
\centerline{\includegraphics[width=5cm]{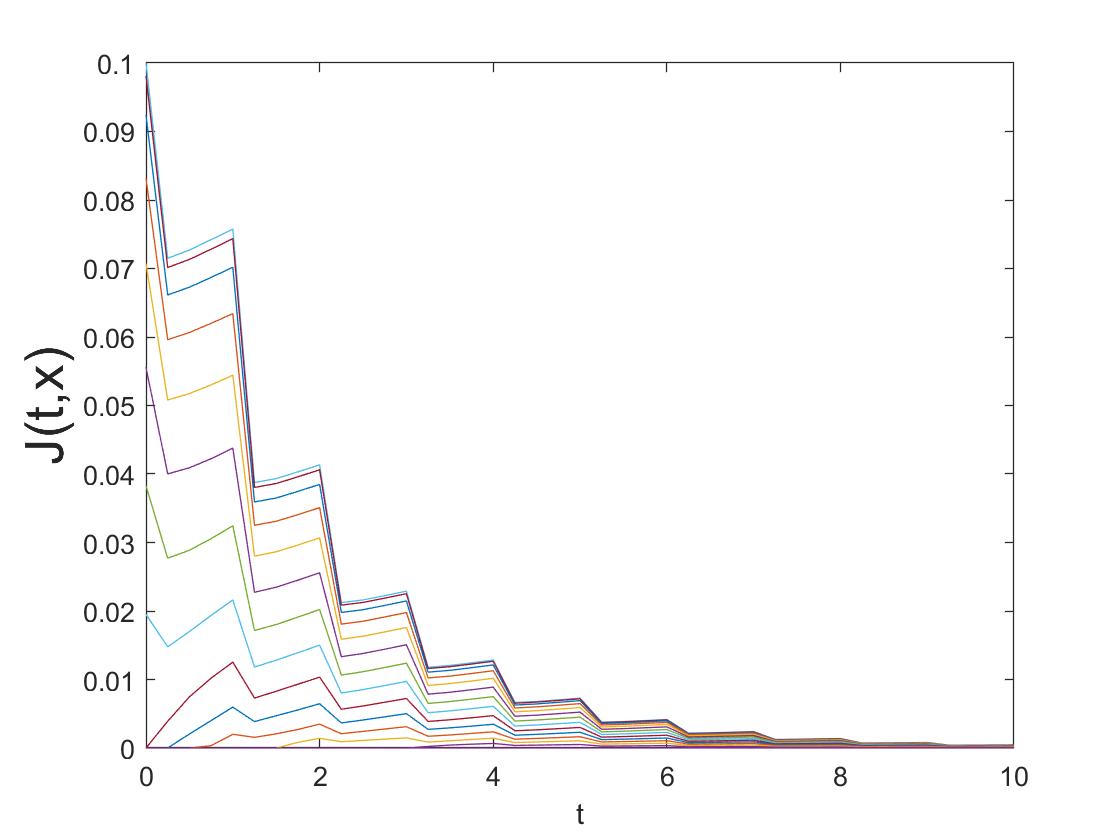}}
\centerline{\small{(a) $\tau=1$}}
\end{minipage}
\begin{minipage}{0.4\linewidth}
\centerline{\includegraphics[width=5cm]{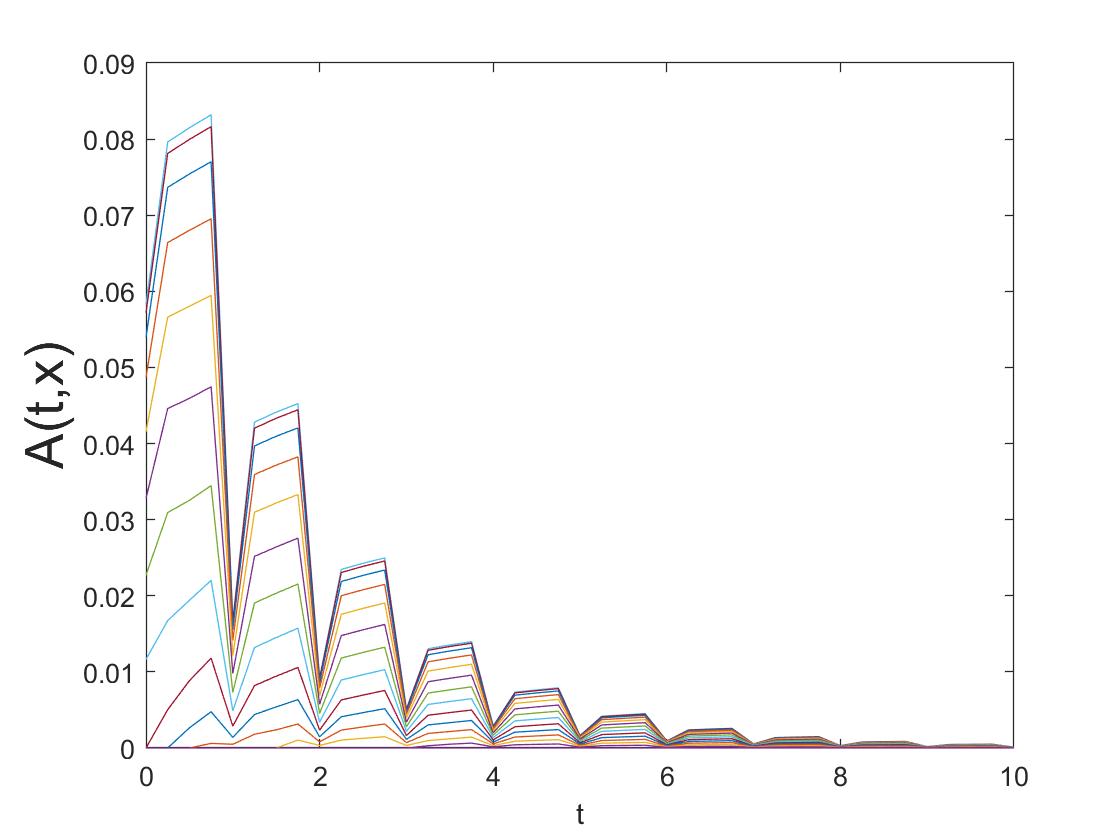}}
\centerline{\small{(b) $\tau=1$}}
\end{minipage}

\quad\quad\quad\quad
\begin{minipage}{0.4\linewidth}
\centerline{\includegraphics[width=5cm]{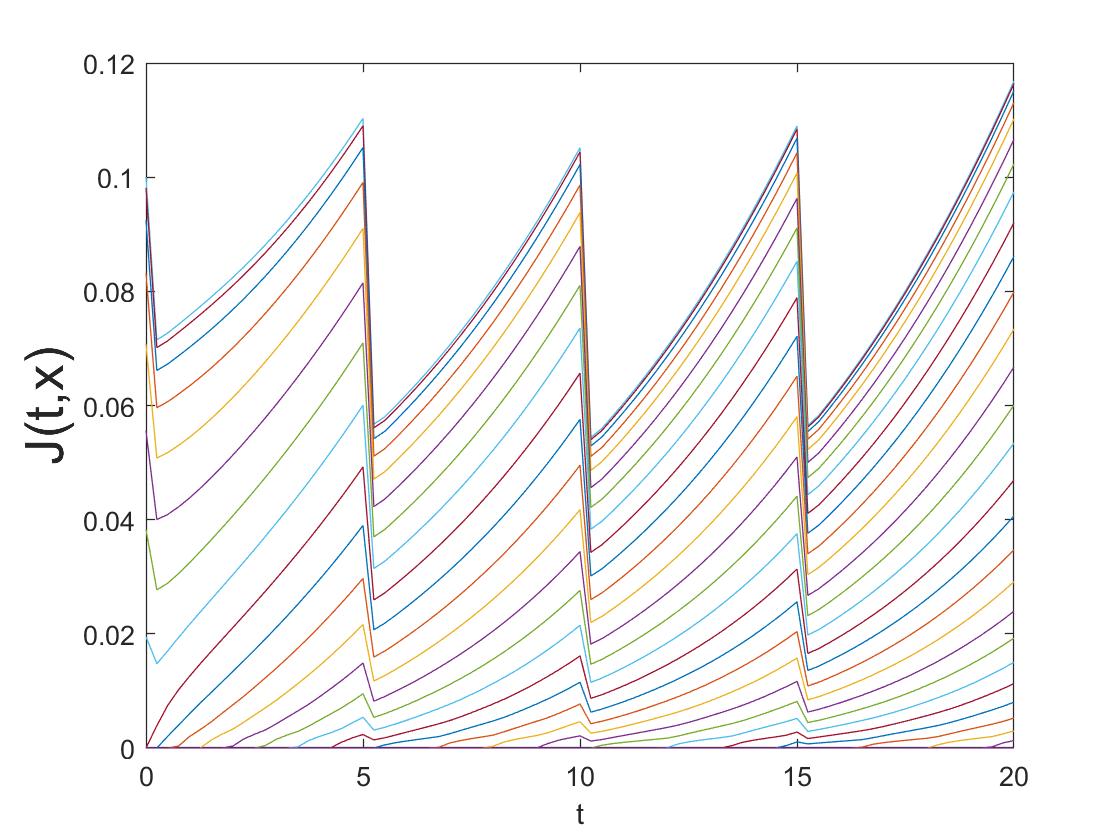}}
\centerline{\small{(c) $\tau=5$}}
\end{minipage}
\begin{minipage}{0.4\linewidth}
\centerline{\includegraphics[width=5cm]{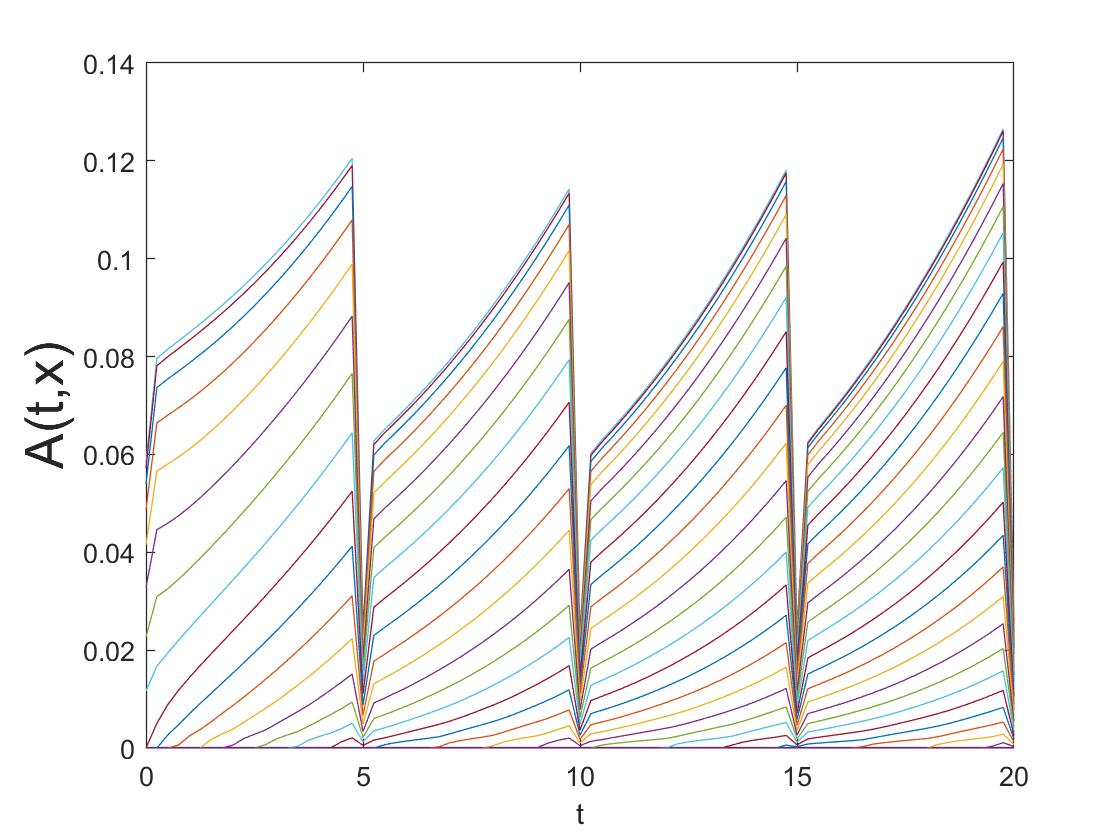}}
\centerline{\small{(d) $\tau=5$}}
\end{minipage}
\caption{\scriptsize Solution of \eqref{a01} for different choices of $\tau$ with big harvesting rate ($m=2,a=10$). All the other parameters are the same as that in Fig. 3. The dynamics of \eqref{a01} vary from vanishing for $\tau=1$ in (a) and (b), to spreading for $\tau=5$ in (c) and (d).}
\label{tu5}
\end{figure}

In contrast to the vanishing of individuals in Fig. 3 with $\tau=2$, we observe from Fig. 5 that when harvesting pulse occurs, short harvesting period with $\tau=1$ in Fig. 5(a-b) will speed up the extinction of species, while long harvesting period with $\tau=5$ in Fig. 5(c-d) will alter the long-time behaviors of solution, from gradual extinction in Fig. 3, to approaching a periodically steady state in Fig. 5(c-d). Analogously, we can get similar conclusions in consideration of Fig. 2 with a spreading case in performing numerical simulations, that is, when harvesting occurs with a small harvesting rate, a longer harvesting period leads to a larger steady state, while a relatively short harvesting period will cause species from expanding to spatially vanishing. Figs 2, 3 and 5 imply that long harvesting periods are conductive to expansion of species, while frequent harvesting pulses lead to the death of species. The above numerical simulations enlighten us that we can adopt external impulsive control covering impulsive intensity and pulse timing on adults to maintain the stability and balance of ecosystem, and finally implement the strategy for sustainable development.

\vspace{2mm}
To sum up, our model generalizes many previous work from physiologically unstructured population with the same dispersal capacity, reproduction rate and mortality rate of all individuals in homogeneous environment to spatially and physiologically structured population in heterogeneous environment. Specifically, we here propose an age-structured juvenile-adult model in moving and heterogeneous environment with periodic harvesting pulse exerting on the adult for describing the dynamical behaviors of species. The global existence and uniqueness of solution to problem \eqref{a01} with harvesting pulse (see Theorem 2.3) is firstly established based on corresponding conclusion of problem \eqref{b01} without pulse (see Theorem 2.1). Then the principal eigenvalue involving periodic harvesting pulse is defined and its properties related to the intensity of harvesting and length of region are investigated in Lemma 3.1. The main goal of the paper is to investigate the criteria for governing spreading or vanishing under the dual effects of harvesting pulse (external factor from human intervention) and expanding capacities (internal factor from species itself). Some sufficient conditions for spreading or vanishing of species are also given, see $H'(0)\leq g_*$ in Corollary 4.1 for vanishing, $H'(0)\geq g^*$ in Theorem 4.1 for spreading, and $g_*<H'(0)< g^*$ in Theorem 4.5 for spreading-vanishing dichotomy. Also, for small initial value, the vanishing case in the condition of $H'(0)< g^*$ is presented in Theorem 4.2. Compared with many other works of reaction-diffusion equations, the introduction of pulse in heterogeneous environment is more natural, and certainly increase the difficulty of theoretical analysis and diversity of complicated outcomes. In fact, the intensity ($1-H'(0)$) and timing ($\tau$) of harvesting pulse can affect or change the dynamical behaviors of solution, see Examples 5.1-5.3 in Section 5. Numerical simulations indicate that small expanding capacities (Fig. 4), large harvesting rate (Fig. 3) and short harvesting period (Fig. 5) on individuals all lead to or accelerate the extinction of the whole population in this juvenile-adult model with harvesting pulse. Moreover, it seems that the intensity and timing of pulse from human intervention will have greater effects on individuals than expanding capacities in the process of performing simulations.

Recently, free boundary problems with different types 
of moving boundaries have been a hot issue and attracted considerable attention. It should be admitted that we only choose a special case of Stefan type boundary condition in this paper, however, there are many other boundary conditions such as \cite{FL} with a nonlocal weighted boundary condition and \cite{PLL} with a nonlocal seasonal type boundary conditions. In addition, the effect of harvesting timing is exhibited here by numerical approximations, and its mathematically proofs may be challenging. Also, the spreading speed in such an age-structured impulsive model has not been solved yet, and we also would like to study whether the long-distant nonlocal effect on boundaries will accelerate or slow down the spreading speed compared with Stefan type boundary condition.
As for how does the harvesting pulse affect the nonlocal spreading, We leave them for further work.

\end{document}